\documentclass[a4paper,12pt,reqno]{article}

\usepackage{amsmath}
\usepackage{amsthm}
\usepackage{amssymb}
\usepackage{latexsym}
\usepackage{graphicx}
\usepackage{stmaryrd}
\usepackage{mathrsfs}
\usepackage{enumerate}
\usepackage{cases}
\usepackage{color}

\newcommand{\BLACK}{\color{black}}

\definecolor{dGREEN}{rgb}{0.0,0.5,0.5}

\newcommand{\glalign}[2]{\lower.6ex\vbox{
\baselineskip\lineskip\ialign{$#1\hfil##\hfil$\crcr#2\crcr=\crcr}}}

\newcommand{\del}{\partial}

\newcommand{\delt}{\partial_{t}}

\newcommand{\dlambda}{\,{\rm d}\lambda}

\newcommand{\dr}{\,{\rm d}r}

\renewcommand{\div}{\mbox{\rm div}\,}

\newcommand{\supp}{\mbox{\rm supp}\,}

\newcommand{\trans}{{}^\top}

\newcommand{\IR}{\mathbb{R}}

\newcommand{\eps}{\varepsilon}

\renewcommand{\Re}{\operatorname{Re}}
\renewcommand{\Im}{\operatorname{Im}}

\def\XXint#1#2#3{{\setbox0=\hbox{$#1{#2#3}{\int}$}
\vcenter{\hbox{$#2#3$}}\kern-.5\wd0}}

\def\eqn#1$$#2$${\begin{equation}\label#1#2\end{equation}}

\numberwithin{equation}{section}

\newtheorem{defi}{Definition}[section]
\newtheorem{thm}[defi]{Theorem}
\newtheorem{cor}[defi]{Corollary}
\newtheorem{prop}[defi]{Proposition}
\newtheorem{lem}[defi]{Lemma}
\newtheorem{rem}[defi]{Remark}
\newtheorem{assumption}[defi]{Assumption}

\def\eqn#1$$#2$${\begin{equation}\label#1#2\end{equation}}

\numberwithin{equation}{section}
\numberwithin{equation}{section}
\allowdisplaybreaks[4]

\begin{document}


\title{
\bf \large Perturbation theory of  the compressible Navier-Stokes equations  and its application\BLACK}

\author{{\normalsize
Kazuyuki Tsuda\footnote{
Kyushu Sangyo University, 3-1 Matsukadai 2-chome,
Higashi-ku, Fukuoka,
813-8503 Japan, \texttt{k-tsuda@ip.}\texttt{kyusan-u.ac.jp}}   }\\[2ex]
}
\date{}
\maketitle

\begin{abstract}
\noindent In this article, a perturbation theory of the compressible Navier-Stokes equations in $\mathbb{R}^n$ $(n \geq 3)$  is studied to investigate decay estimate of solutions around a non-constant state. As a concrete problem, stability is considered for a  perturbation system  from a stationary solution $u_\omega$ belonging to the weak $L^n$ space. 
Decay rates of the perturbation including $L^\infty$ norm are obtained 
which coincide with those of the heat kernel except a bit loss. 
The proof is based on deriving suitable resolvent estimates with perturbation terms in the low frequency part having a parabolic  spectral curve. Our method can be applicable to dispersive hyperbolic systems like wave equations with strong damping.  Indeed, 
a parabolic type decay rate of a solution is obtained  for a  damped wave equation including variable coefficients which satisfy  spatial decay conditions.  
\end{abstract}

\\


\begin{section}{Introduction} 

We consider perturbation theory of    the compressible Navier-Stokes equations (CNS) in $\mathbb{R}^n$  $(n\geq 3)$:
\begin{eqnarray}
\left\{
\begin{array}{ll}
\partial_{t}\rho +\nabla\cdot(\rho v)=0,\\
\rho(\partial_{t}v+(v\cdot\nabla)v)-\nu\Delta v-(\nu+\nu')\nabla(\nabla\cdot v)+{\nabla p(\rho)}=\rho g.
\label{(1.1)}
\end{array}
\right.
\end{eqnarray}
\noindent Here $\rho=\rho(x,t)$ and $v=(v_{1}(x,t),\cdots,v_{n}(x,t))$ denote the unknown density 
and the unknown velocity field, respectively, at time $t\geq 0$ and 
position $x\in\mathbb{R}^n$; 
$p=p(\rho)$ is the pressure that is assumed to be a smooth function of $\rho$ satisfying 
$$
p'(\rho_*)>0
$$
for a given positive constant $\rho_*$; 
$\nu$ and $\nu'$ are the viscosity coefficients that are assumed to be constants satisfying 
$$
\nu>0, \ \ \ 
\frac{2}{n}\nu+\nu'\geq 0;
$$ 
$g=g(x,t)$ is a given external force. 

As a concrete problem, we consider decay estimates of  a solution around a non-constant state, {\rm i.e., } a perturbation system of \eqref{(1.1)} from a solution $u_\omega$.  We assume the following property for the solution. 

\vspace{2ex}

\begin{assumption}\label{ass}
Let $u_\omega(x)=\trans(\phi_\omega(x)+\rho_*, w_\omega(x))$ be a stationary solution to \eqref{(1.1)} around the constant state $\trans(\rho_*, 0)$ satisfying that 
\begin{align}\label{u-omega-esimate}
\begin{aligned}
\|(\phi_\omega, w_{\omega})\|_{L^{n,\infty}} + \|\nabla (\phi_\omega, w_{\omega})\|_{L^{\frac{n}{2}, \infty}}+ \|\nabla (\phi_\omega, w_\omega )\|_{{H}^{s}}
 \leq \epsilon,
\end{aligned}
\end{align}
where $s$ is an integer satisfying $s \geq [n/2]+2$. 
\end{assumption}

{\bf Example.} Shibata and Tanaka  \cite{Shibata-Tanaka}  show existence of a stationary solution  $u_s=\trans(\phi_s-\rho_*, w_s)$ around a given constant state $(\rho_*, 0)$ to \eqref{(1.1)} on $n=3$ satisfying that 
\begin{align}
\begin{aligned}\label{Pr-sol}
\|(1 & +|x|)^2\phi_s\|_{L^\infty} + \|(1+|x|)w_s\|_{L^\infty} \\
& + \sum_{k=1}^4\|(1+|x|)^{k}\nabla^k \phi_s\|_{L^2} + \sum_{k=1}^5\|(1+|x|)^{k-1}\nabla^k w_s\|_{L^2} \ll 1,  
\end{aligned}
\end{align}
under a given stationary force. By the Sobolev embedding,  this solution satisfies \eqref{u-omega-esimate} as $n=3$.  We note that the above spatial decay rates of $L^\infty$ norm, 
{\rm i.e., } $u_s \sim 1/|x|$ and $\nabla u_s \sim 1/|x|$ as $|x| \rightarrow \infty$ 
are the same as those of a physically reasonable solution, which is known as {\em PR} solution in the
sense of Finn \cite{Finn}, constituted by velocity fields decaying  like $|x|^{-1}$ 
with their spatial gradients decaying like $|x|^{-2}$ on exterior domains.  Concerning the exterior domain on CNS $(n=3)$, Novotny and Padula \cite{Novotny-Padula2} show 
existence of a stationary solution having the same spatial decay rates at infinity. 

\vspace{2ex}

\begin{rem} 
{\rm 
In Assumption \ref{ass},  the density has a slower spatial decay than \cite{Shibata-Tanaka}. Indeed, in  the non-isentropic case, the spatial decay of density may be slower than that in the isentropic case as shown by the author \cite{Tsuda-K} for the time periodic solution. Furthermore, 
we  assume larger classes $L^{n,\infty}$ and $L^{n/2,\infty}$ than the weighted $L^\infty$ norms  for  $u_\omega$ and $\nabla u_\omega$ respectively.   
Those classes have well been studied for stability of a  stationary solution to the incompressible Navier-Stokes equations, {\rm e.g.,} Novotny and Padula \cite{NoPa}, 
Borchers and Miyakawa \cite{BoMi95}, Kozono and Yamazaki \cite{Kozono-Yamazaki} and Koba \cite{Koba}. On the other hand, as far as we investigated, there is no result  for the compressible one.  The higher order derivatives only belong to the $L^2$ homogeneous Sobolev space. 
The spatial  decay in Assumption \ref{ass} is slower  than 
the $n$-dimensional Newton potential $O(1/|x|^{n-2})$.  Handling such  slower decaying terms will be reasonable for the application of our method, {\em e.g. } to the damped wave equation. \BLACK
In fact, Assumption \ref{ass-DWE} does not require any weighted   
$L^\infty$ norm, as well as weighted $L^2$ Sobolev norm below. 
}
\end{rem}

\vspace{2ex}

Decay estimates of solutions around a constant state have been studied by many researchers. Matsumura and Nishida \cite{Matsumura-Nishida1, Matsumura-Nishida2} show the global strong solutions for small perturbations from a linear stable constant state $\trans(\rho_*, 0)$ in $\mathbb{R}^3$. They also show a decay rate of solutions in $L^2$ norm for small $L^1$ initial data:  
$$
\|(\rho-\rho_*, v)\|_{L^2} \leq C(1+t)^{-\frac{3}{4}}. 
$$
The decay rate is optimal which coincides with that of the heat kernel. Ponce \cite{Ponce} obtains the decay rates in $L^q$ class:  
$$
\|\nabla^k (\rho-\rho_*, v)\|_{L^q} \leq C(1+t)^{-\frac{n}{2}\Big(1-\frac{1}{q}\Big)-\frac{k}{2}},  \ \ 
2 \leq q \leq \infty,  \ \ 0\leq k \leq 2
$$
for $n=2,3$.
On other unbounded domains, decay estimates of solutions also have been studied. There are results by Kobayashi \cite{Kobayashi} and Kobayashi and Shibata \cite{Kobayashi-Shibata2} for exterior domains and Kagei and Kobayashi \cite{Kagei-Kobayashi1, Kagei-Kobayashi2} for the half space case.  
Hoff and Zumbrun \cite{HZ} investigate the diffusion wave phenomenon \BLACK which affects large time behavior of solutions by the hyperbolic aspect of the system. 
\cite{Kobayashi-Shibata} also studies the diffusion wave phenomenon in details for the linear system. For decay rates of solutions in the critical $L^2$ frame-work, we cite 
Okita \cite{Okita} for $n \geq 3$, and Danchin \cite{Danchin} which includes the two dimensional case.  Danchin and Xu \cite{Danchin-Xu} study  decay rates \BLACK in the critical $L^p$ frame-work. 


\vspace{2ex}



On the other hand, to derive decay rates of solutions around a non-constant state as in our case is more difficult. We consider a perturbation system of $u=\trans(\phi, w)=\trans(\rho-\rho_\omega, v-v_\omega)$ from the solution to \eqref{(1.1)} 
with an initial perturbation 
\begin{eqnarray}
 u|_{t=0}=u_{0}=\trans{(\phi_{0},w_{0})}.\label{(1.2)}
\end{eqnarray}
In \cite{Shibata-Tanaka} the $L^2$ energy stability is also obtained and it is shown that 
$$
\|u\|_{L^\infty} \rightarrow 0  \ \ (t \rightarrow \infty).
$$
As a continued work, Shibata and Tanaka \cite{Shibata-Tanaka2} study the stability again and obtain the following results.   

\vspace{2ex}

\begin{thm}
There exists a solution $u$ globally in time to $\eqref{(1.1)}$ and $\eqref{(1.2)}$ for $u=\trans(\rho-\rho_s, v-v_s)$. Moreover, if 
$$\|u_0\|_{L^{\frac{6}{5}}\cap H^3} \ll 1, $$
it holds that  
\begin{align}\label{decay-est-previous}
\|u(t)\|_{\widehat{H}^1} \leq C(1+t)^{-\frac{1}{2}+\epsilon} 
\end{align}
for $0<t$, where $\widehat{H}^1$ is the homogeneous $L^2$ Sobolev space and 
$\epsilon>0$ is any small positive number.  This together with some space-time integral stated in {\rm \cite{Shibata-Tanaka2}} yields 
\begin{align}\label{Linfty-norm-decay}
\|u(t)\|_{L^\infty}=O(t^{-\frac{1}{2}+\epsilon})  \ \ \mbox{as}  \ \ t\rightarrow \infty. 
\end{align}
\end{thm}

The decay rates \eqref{decay-est-previous}, as well as \eqref{Linfty-norm-decay} seem to be not optimal.   Indeed, \cite{Kobayashi-Shibata}  shows \BLACK that  the low frequency part of solutions to the linearized CNS has the same decay rate as that  of (in) \BLACK the heat kernel. On the other hand, the high frequency has exponential decay due to the spectrum away from the origin. 
Hence it is expected that 
\begin{align}\label{decay-est-original}
\|u(t)\|_{L^\infty} \leq C(1+t)^{-\frac{n}{2p}},  \ \  \|u(t)\|_{\widehat{H}^1} \leq C(1+t)^{-\frac{n}{2}\Big(\frac{1}{2}-\frac{1}{p}\Big)-\frac{1}{2}}  \ \ \mbox{for }  \ \ 0<t,  
\end{align}
for small $u_0 \in L^p$ with $1<p<2$. 
In view of this point, Deguchi  \cite{Deguchi} studies the stability problem on the three dimensional case in homogeneous Besov spaces. He shows that there exists a stationary solution $\trans(\rho_s, v_s )=\trans(\phi_s+\rho_*, v_s)$ with $ \phi_s \in \widehat{B}^{-1/2}_{2,\infty} \cap \widehat{H}^4$ and $v_s \in \widehat{B}^{1/2}_{2,\infty} \cap \widehat{H}^5$
for small steady external force, and the perturbation $u=\trans(\phi, w)=\trans(\rho-\rho_s, v-v_s)$ has the following decay estimate:  
if $u_0$ is small in $L^p \cap \widehat{B}^{1/2}_{2,\infty} \cap \widehat{H}^3$, then 
$$
\|u\|_{\widehat{H}^s} \leq C(1+t)^{-\frac{n}{2}\Big(\frac{1}{p}-\frac{1}{2}\Big)-\frac{s}{2}}\|u_0\|_{L^p \cap H^3}
$$
for $1 \leq p \leq 2$ and $-3/2 < s < 3/2$. 
The crucial point is to apply the time-space $L^1$ integral estimates like the Yamazaki estimates in the weak type Besov space.  

\vspace{2ex} 

\begin{rem}\label{about-solution-space}
{\rm  
Concerning the stationary  solution $u_\omega$ as in Assumption \ref{ass},
even if $\phi_\omega \sim 1/|x|$ as  $|x| \rightarrow \infty$, $\phi_\omega$ will not belong to $\widehat{B}^{-1/2}_{2,\infty}$ according to Danchin \cite[Proposition 2.21]{Danchin}. 
In addition, 
if $v$ satisfies that $(1+|x|) v \in L^\infty$, $(1+|x|)^2 \nabla v \in L^\infty$, 
\BLACK then one can show that $v \in \widehat{B}^{1/2}_{2,\infty}$. 
On the other hand, 
it is unclear  that $w_\omega$  belongs to $\widehat{B}^{1/2}_{2,\infty}$ because the conditions $w_s \in L^{n,\infty}$ and $\nabla w_\omega \in L^{\frac{n}{2}, \infty}$ are weaker than the weighted $L^\infty$ norms.  We can show that 
$L^{3,\infty} \not\subset \widehat{B}^{1/2}_{2,\infty}$ on $\IR^3$. See Appendix for details below.  If $w_s$  decays like $1/|x|$ at infinity, it does not give any control of higher order derivatives  due to small oscillations which show up
under derivatives, but not so much for $w_s$  itself. On the other hand, the function  $1/(1+|
x|)$  is well controlled in the Besov space of order $1/2$ and its
derivatives of higher order and thus $1/|x| \in \widehat{B}^{1/2}_{2,\infty}$ as in \cite[Proposition 2.21]{Danchin}.
}
\BLACK 
\end{rem} 

\vspace{2ex} 

These  results are different from those for \BLACK the incompressible Navier-Stokes equations. Let $\omega=\omega(x)$ be a stationary solution obtained by Borchers and Miyakawa \cite[Theorem 2.3 (ii)]{BoMi95} for given external force field $F=\div G$ with matrix field $G$ being sufficiently small in weighted $L^\infty$ spaces. \BLACK Then there exists a unique solution  $\omega$ of the stationary Navier-Stokes system, 
\begin{align}\label{equ:ns-stat}
		- \nu\Delta \omega + \omega\cdot \nabla \omega +  \nabla p =  F, \quad 
		 \div\omega = 0\; \text{ in }\; \Omega, \quad
	\omega = 0 \; \text{ on }\, \partial \Omega,
\end{align}
\BLACK satisfying the estimate 
\begin{align}\label{PR}
	\begin{aligned}
\||x|^{ n-2}\omega\|_{L^\infty}+ \||x|^{ n-1} \nabla \omega\|_{L^\infty} \leq \epsilon, 
\end{aligned}
\end{align}
where $\epsilon=\eps(\Omega, \nu, G)>0$ is a small constant.  \cite[Theorems 6.3 and  6.8\BLACK]{BoMi95} \BLACK prove asymptotic stability in the  class \eqref{PR} together with the following decay estimates: 
Suppose $u_0 \in L^{n,\infty}_\sigma \cap L^{p,\infty}_\sigma$ with $1<p<n$ and that
$\|(1+|x|)\omega\|_{L^\infty}+\|(1+|x|)^2 \nabla \omega\|_{L^\infty}+ \|u_0\|_{L^{n,\infty}}$ is small. 
Then
\begin{align}\label{decay-est-intro}
	\begin{aligned}
\|u(t)\|_{L^r} & = O\big(t^{-\frac{n}{2}\big(\frac{1}{p}-\frac{1}{r}\big)}\big),  \ \ \|u(t)\|_{L^\infty} = O(t^{-\frac{n}{2p}+\epsilon}),\\	
\|\nabla u(t)\|_{L^q} & = O\big(t^{-\frac{n}{2}\big(\frac{1}{p}-\frac{1}{q}\big)-\frac{1}{2}}\big), 
\|\nabla u(t)\|_{L^n} =O\big(t^{-\frac{n}{2p}+\epsilon}\big), 
\end{aligned}
\end{align}
as $t$ goes to $\infty$, where $1<p<r<\infty$ satisfying  $p \leq q <n$ and $\epsilon$ is any small positive number.  Our aim is to show 
decay estimates like \eqref{decay-est-intro} of the perturbation from the solution satisfying Assumption \ref{ass}.   Since, as shown in \cite{Kobayashi-Shibata}, the linearized system around a constant state has parabolic type decay rates, it is expected that we can get similar decay estimates to \eqref{decay-est-intro} in our case. 

Our main result is stated as follows. 

\vspace{2ex}

\begin{thm}\label{stability}
Let $ n \geq 3$. Suppose that $\|u_0\|_{L^p \cap H^s}\ll 1$ for $1<p < 2$, $s \geq \big[\frac{n}{2}\big]+2 $ and 
$$
 1<  \frac{n}{2}\Big(\frac{1}{p}-\frac{1}{2}\Big)+\dfrac{1}{2},  \ \ \max\Big\{0, 1-\dfrac{n}{4}\Big\}<  \frac{n}{2}\Big(\frac{1}{p}-\frac{1}{2}\Big)< \min\Big\{\dfrac{n}{4}, 1\Big\}. 
$$
In addition, $\dfrac{1}{p}=\dfrac{1}{2}+\dfrac{1}{p_0}$ and $k=\dfrac{n}{2}\Big(\dfrac{1}{2}-\dfrac{1}{p_0}\Big)$ for $2<p_0$  with $2k\leq s$.  
We also suppose that viscosity coefficients are large as in Proposition \ref{resolvent2-each-point} below. 
Then there exists a global solution $u$ to $\eqref{(1.1)}$-$\eqref{(1.2)}$ in $C([0,\infty); H^s)$ satisfying that 
$$
\|u(t)\|_{L^2} \leq C(1+t)^{-\frac{n}{2}\big(\frac{1}{p}-\frac{1}{2}\big)+\epsilon_1},  \ \ 
\|\nabla u(t)\|_{L^2}  \leq C(1+t)^{-\frac{n}{2}\big(\frac{1}{p}-\frac{1}{2}\big)-\frac{1}{2}+\epsilon_1}  \ \ \mbox{for}  \ \ t>0 
$$
and 
$$
\|\nabla u(t)\|_{L^n}+\|u(t)\|_{L^\infty} \leq C(1+t)^{-\delta} \ \ \mbox{for}  \ \ t>0,
$$
where 
$\delta =\min\Big\{\dfrac{n}{2p}-\epsilon_1, 1+ \dfrac{n}{2}\Big(\dfrac{1}{p}-\dfrac{1}{2}\Big)-\epsilon_1\Big\}$ 
for any small positive $\epsilon_1$. 
\end{thm}

\vspace{2ex}

{\rm 
\begin{rem}
Concerning $n=3$, the decay rates  coincide with those of the heat kernel except a bit loss. 
For $4 \leq n$, the decay rates have a restriction 
for the upper bound of decay rate, while we consider  the slower spatial decaying solution in Assumption \ref{ass} than 
the $n$-dimensional Newton potential $O(1/|x|^{n-2})$ as $|x|\rightarrow \infty$.  
\end{rem}
}

\vspace{2ex}
As an application of our perturbation theory, we consider the wave equations with strong damping and perturbation terms as follows. 

\vspace{2ex}

\begin{align}
\begin{aligned}\label{eq-DWE}
    \del_{tt} u - \mu \Delta \del_t u - \mu' \Delta u + B[b]u =0, \\
u(0,x)=u_0,  \ \ \del_t u(0,x)=v_0,
\end{aligned}
\end{align}
where 
$$
B[b]u(t,x)= u \div  b_1(x)  + b_2(x) \cdot \nabla u + b_3 (x) \Delta u. 
$$
We assume smallness and decay property of the coefficients $b_j$ as follows. 

\begin{assumption}\label{ass-DWE}
Let $b_j(x)$ $(j=1,2,3)$ satisfy that 
%
$$
\sum_{j=1,2,3}(\|b_j\|_{L^{n,\infty}} 
+\|\nabla  b_j\|_{L^{\frac{n}{2},\infty}}) +\|(\nabla b_1, b_2, b_3, \nabla b_3)\|_{L^{p_0}}+ \|b_3\|_{W^{1,\infty}} \leq \epsilon, 
$$
where $1/p_0= 1/p-1/2$ for $1<p<2$.   
\end{assumption}

\vspace{2ex}

We have the following a priori estimate of a solution $u$. 

\vspace{2ex}

\begin{thm}\label{decay-DWE}
Let $ n \geq 5$, $u \in C([0,T]; H^1) \cap L^2 (0,T; H^2)$ with $\del_t u \in C([0,T]; L^2)$ 
be a solution to \eqref{eq-DWE} with 
$(u_0,v_0) \in L^p \cap H^2$ for $T>0$ and $1 <p <2$. We suppose that 
$\mu$ is large as in  Proposition \ref{resolvent2-each-point-DEW} and \eqref{mu-condition} below. 
Then $u$ satisfies 
$$
\|u\|_{L^2} \leq C(1+t)^{-\delta+\epsilon_1}(\|v_0\|_{L^p\cap H^1 }+ \|w_0\|_{L^p \cap H^1})
$$
and 
$$
\|\nabla u\|_{L^2} \leq C(1+t)^{-\delta+\epsilon_1-\frac{1}{2}}(\|v_0\|_{L^p\cap H^1 }+ \|w_0\|_{L^p \cap H^1}),  
$$
where $\delta=\frac{n}{2}\Big(\frac{1}{p}-\frac{1}{2}\Big)$ satisfies $\max\{0, 1-\frac{n}{4}\}< \delta< \min\{\frac{n}{4}, 1\}$ and 
$\delta+\frac{1}{2}>1$ for any small $\epsilon_1>0$ and the positive constant $C$ is independent of $T$. 
\end{thm}

\vspace{2ex}



One of the difficulties of this problem is to establish time decay estimates of the solution operator for the perturbed system. 
The spectral curve is a parabola in the low frequency part due to the hyperbolic aspect of the system.  
Indeed, density of fluid is governed by the damped wave equation which has hyperbolic aspect. 
Hence it is difficult to get suitable resolvent estimates, especially, around the origin.  
Let us consider the Oseen equations on the three dimensional case.  
The spectral curve is also a parabola as shown in Farwig and Neustupa \cite{Farwig-Neustupa}.  
In this case, it is shown by Deuring and Varnhorn \cite{ Deuring, Deuring-Varnhorn} that the resolvent estimate have a worse order for the resolvent parameter $\lambda$ near $\lambda=0$ than that of the sectorial  operator. More preciously, a $L^p$ norm of the Oseen resolvent is bounded by $C|\lambda|^{-2}$ with $C>0$ independent of $\lambda$ near $\lambda=0$, and such an estimate cannot hold  if $|\lambda|^{-2}$ is replaced by $|\lambda|^{-\kappa}$ with $\kappa \in [0,3/2)$.  On resolvent estimates of CNS away from the origin, we can see, {\rm e.g.,} Shibata and Tanaka \cite{Shibata-Tanaka3}, Kobayashi, Murata and Saito \cite{KMS} and Huang and Luo \cite{Huang-Luo} in various settings.  Furthermore, the resolvent estimate will not satisfy the Fourier multiplier condition. In fact,  
the linearized solution may grow up in $L^p$ norm on $1 \leq p<2$.  This is called the diffusion wave property which is found in \cite{Kobayashi-Shibata}. See, Remark \ref{no-multiplier} below.   
This situation is the reason why it is different from the incompressible Navier-Stokes problem which is a parabolic one. We have to consider  an integral curve surrounding the parabolic curve in the Cauchy formula of the semigroup.
Another difficulty is that the perturbation terms from the solution satisfying \eqref{u-omega-esimate} belongs to $L^{n,\infty}$ for $u_\omega$ and $L^{n/2,\infty}$ for $\nabla u_\omega$, {\rm i.e.,}  
 has slow decays at spatial infinity, especially, weaker ones than  the inverse-square type potential. As explained in \cite{Deguchi}, the inverse-square type potential terms can
prevent the asymptotic behavior of the solution of the heat equation.  This is the reason why weak-type Besov spaces are used in \cite{Deguchi}.

Our crucial idea is to derive suitable resolvent estimates and decay estimates of a perturbed system in the low frequency part.  Concretely, we consider the cases: (1) $\mbox{Re} \lambda <0$, (2) $\mbox{Re} \lambda >0$,  
where $\lambda$ is the resolvent parameter. 
Estimates on $\mbox{Re} \lambda <0$ near $\lambda=0$ are most difficult. 

On (1) and (2), we estimate by the explicit fundamental solution formula. First we investigate several cases on the parameter $\xi$ of the Fourier transform. 
For the crucial part, {i.e., } $\lambda \sim |\xi|$ with $|\xi| \rightarrow 0$ and $\Re \lambda <0$, we make use of the property of a parabolic curve which we take on the Cauchy integral. We get resolvent estimates with some loss of $|\xi|$ order, which is integrable on $L^2$ in the low frequency part.   Under the resolvent estimates based on (1)-(2) we take some parabola type integral curve surrounding the spectral one in the Cauchy formula of the semigroup. Then we obtain suitable decay estimate of the semigroup. We note that since we take a parabola curve in the integral, 
additional growth decay rate $+1/2$ appears in the time decay estimates. To cancel of this growth term, we further obtain improvement of the resolvent estimates in the Fourier space by combing the precious pointwise resolvent estimates with integrability of $1/|\xi|$ in the low frequency part, which holds on the parabola spectral curve. See Corollary \ref{resovent-estimate-modify-0} below. 

Perturbation terms are included in the semigroup estimate by  perturbed resolvent estimates based on the spatial decays of the terms.  To handle perturbation terms in resolvent estimates,  Lemma \ref{resolvent-est-L1} plays an important role. 
The restriction of exponents appear from Lemma \ref{fgh-est} which is enough to get $L^p$-$L^2$ estimate for $1<p<2$ close to $1$.  The crucial point that the fractional Laplacian  ``passes through" the resolvent with perturbation terms in $L^2$ norm, {\rm i.e., } the fractional Laplacian of the resolvent is bounded by the homogeneous Sobolev norm in $L^2$.  
A similar type estimate is found by Borchers and Miyakawa, \cite[Proposition 3.4, 3.5]{BoMi95} with a first order derivative $\nabla$, {\rm i.e., } $(-\Delta)^{1/2}$.  

On the high frequency part, we use the standard energy method on $L^2$ including perturbation terms. To estimate perturbation terms, the Poincar\'e type inequality is used due to the high frequency.  Our method does not need the conservation form for  the nonlinear estimates due to the perturbation theory. We only apply a product rule \eqref{derivative-form} below. This will have a benefit when we apply our method to the damped wave equations.

This method for the perturbed systems  can be applicable to other  dispersive hyperbolic systems like wave equations with strong damping.  Since the resolvent estimates rely on the explicit solution formula in the Fourier space,  we can apply our method to other systems having solution formulas including parabolic spectral curves with a similar order of $|\xi| \sim 0$ to ours. (Cf., Ponce \cite{Ponce} and Shibata \cite{Shibata-2000}). 
Indeed, we can obtain the $L^p$-$L^2$ decay estimate of solutions to the damped wave equations with perturbation terms in Section 7.  Decay estimates of the wave equations
 with variable coefficients have been  studied by Ikehata, Todorova and Yordanov \cite{ITY}.  In contrast to \cite{ITY}, our result can handle 
 variable coefficients $(\div b_1(x)) u+b_2(x) \cdot \nabla u + b_3 (x)\Delta u $ which satisfy some spatial decay conditions without  divergence form,   see  Assumption \ref{ass-DWE}.  
 Concerning existence of a solution, see Remark \ref{existence-sol-DWE} below.  
 The result is stated on $5 \leq n$ and the lower dimensional case will be handled in a forth-coming paper.  
 
\vspace{2ex}

This paper is organized as follows. 
In section $2$, we introduce notations and auxiliary lemmas used in this paper. In section $3$ we reformulate the problem. 
Section $4$ is devoted to studying the low frequency part 
and we derive the estimates for low frequency part. 
In section $5$,  we establish the  energy estimate 
for the high frequency part. 
In section $6$, we estimate nonlinear terms 
and then give a proof of the stability theorem. 
In section 7, the damped wave with perturbation terms is studied. In section 8, we write details about Remark \ref{about-solution-space}. 
\end{section}

\begin{section}{Preliminaries}
In this section we first introduce some notations which will be used throughout this paper. 
We then introduce some auxiliary lemmata which will be useful in the proof of the main results. 

For a given Banach space $X,$ the norm on $X$ is denoted by $\|\cdot\|_{X}$.

Let $1\leqq p\leqq \infty.$ $L^p$ stands for the usual $L^p$ space over $\mathbb{R}^n$. 
The inner product of $L^2$ is denoted by $(\cdot , \cdot)$. 
For a nonnegative integer $k$, $H^k$ stands for the usual $L^2$-Sobolev space of order $k$. 
As usual, $H^0=L^2$. 

The set of all vector fields 
$w=\trans(w_1,\cdots,w_n)$ on $\mathbb{R}^n$ 
with $w_j\in L^p$ $(j=1,\cdots,n)$, i.e., $(L^p)^n$, is simply denoted by $L^p$; 
and the norm $\|\cdot\|_{(L^p)^n}$ is denoted by $\|\cdot\|_{L^p}$ 
if no confusion will occur. 
Similarly, for a function space $X$, 
the set of all vector fields 
$w=\trans(w_1,\cdots,w_n)$ on $\mathbb{R}^n$ 
with $w_j\in X$ $(j=1,\cdots,n)$, i.e., $X^n$, is simply denoted by $X$; 
and the norm $\|\cdot\|_{X^n}$  is denoted by $\|\cdot\|_{X}$ if no confusion will occur. 
(For example, $(H^k)^n$ is simply denoted by $H^k$ 
and the norm $\|\cdot\|_{(H^k)^n}$ is denoted by $\|\cdot\|_{H^k}$.)

For $u=\trans(\phi,w)$ with $\phi\in H^k$ and $w=\trans(w_1,\cdots,w_n)\in H^m$, 
we define the norm $\|u\|_{H^k\times H^m}$ of $u$ on $H^k\times H^m$ by 
$$
\|u\|_{H^k\times H^m}=\left(\|\phi\|_{H^k}^2+\|w\|_{H^m}^2\right)^{\frac{1}{2}}.
$$
When $m=k$, we simply write $H^k\times (H^k)^n$ as $H^k$, and, also, 
$\|u\|_{H^k\times (H^k)^n}$ as $\|u\|_{H^k}$ 
if no confusion will occur: 
$$
H^k:=H^k\times (H^k)^n, 
\ \ \ 
\|u\|_{H^k}:=\|u\|_{H^k\times (H^k)^n} 
\ \ \ (u=\trans(\phi,w)).
$$
Similarly, when $u=\trans(\phi,w)\in X\times Y$ with $w=\trans(w_1,\cdots,w_n)$ 
for function spaces $X$ and $Y$, 
we denote its norm $\|u\|_{X\times Y}$ by 
$$
\|u\|_{X\times Y}=\left(\|\phi\|_{X}^2+\|w\|_{Y}^2\right)^{\frac{1}{2}}
\ \ \ (u=\trans(\phi,w)).
$$
When $Y=X^n$, 
we simply write $X\times X^n$ as $X$, 
and also its norm $\|u\|_{X\times X^n}$ as $\|u\|_X$: 
$$
X:=X\times X^n, 
\ \ \ 
\|u\|_{X}:=\|u\|_{X\times X^n}
\ \ \ (u=\trans(\phi,w)).
$$

Let $\alpha=(\alpha_{1},\cdots ,\alpha_{n})$ be a multi-index. We use the following notation 
\begin{eqnarray}
\partial^{\alpha}_{x}
=\partial^{\alpha _{1}}_{x_1}\cdots\partial^{\alpha _{n}}_{x_n},
\qquad |\alpha|=\sum_{j=1}^{n}\alpha_{j}.\nonumber
\end{eqnarray}
For any integer $l\geq 0$, $\nabla^{l}f$ denotes $x$-derivatives of order $l$  of a function $f$.

We denote by $\hat{f}$ or $\mathcal{F}[f]$ the Fourier transform of $f$: 
\begin{eqnarray}
\hat{f}(\xi)
=\mathcal{F}[f](\xi)
=\int_{\mathbb{R}^n}f(x)e^{-ix\cdot\xi}dx\quad (\xi\in\mathbb{R}^n).\nonumber
\end{eqnarray}
The inverse Fourier transform of $f$ is denoted by $\mathcal{F}^{- 1}[f]$: 
\begin{eqnarray}
\mathcal{F}^{- 1}[f](x)
=(2\pi)^{-n}\int_{\mathbb{R}^n}f(\xi)e^{i\xi\cdot x}d\xi\quad(x\in\mathbb{R}^n).\nonumber
\end{eqnarray}

We consider to decompose solutions to a low frequency part and a high frequency part of  below. Let $\chi_1\in C^{\infty}_0(\IR^n)$ and $\chi_\infty=1-\chi_1$ be smooth cut-off functions satisfying 
\begin{align}
		 0 \leq  \chi_1 \leq 1, \ \ \chi_1(\xi) =1 \ (\xi \in Q_{r'_1}), \ \ 
		\chi_1(\xi) =0  \ (\xi \in \mathbb{R}^n \setminus Q_{r'_\infty} ), 
\end{align}
where  $Q_r\subset \IR^n$ \BLACK denotes a cube with center at the origin and side length $r$ and $0< r'_1 < r'_\infty$. 
Let us introduce multiplier operators $P_1$ and $P_\infty$  
defined by \BLACK
	\begin{eqnarray}
		P_{j}f=\mathcal{F}^{- 1}\hat\chi_{j}\mathcal{F}[f]\quad(f\in L^p , j=1,\infty). 
	\end{eqnarray}
We denote explanation about the low and high frequency part, respectively. 
\BLACK We define 
$$
L_{(1)}^{p}:= \{u \in L^p; \  \mbox{supp }\hat{u} \subset Q_{r'_\infty}\}, \ \ 1<p \leq 2  
$$
and 
$$
H_{(\infty)}^{k}:=\{u\in H^k; \ \mbox{supp }\hat{u}\subset Q_{r'_1} \}
$$
for a nonnegative integer $k$. Let $r_{1}$ and $r_{\infty}$ be positive constants  satisfying $r_{1}<r_{\infty}$, $r_1= 2r'_1$, $r_\infty =2r'_\infty$. We note that 
for $u \in L^p_{(1)}$ $\supp \hat{u} \subset \{|\xi| \leq r_\infty\}$ and for  $u \in L^2_{(\infty)}$ $\supp \hat{u} \subset \{|\xi| \geq r_1\}$.  In addition,  $H^k \cap L_{(1)}^2=L_{(1)}^2$ for any nonnegative integer $k$. See Lemma \ref{lemP_1} below. 
We note that $L^2_{(1)}$ is a complete metric space with the $L^2$ inner product, {\rm i.e., } the Hilbert space. Hence $(L^2_{(1)})'= L^2_{(1)}$. 
    


Let $s \in \mathbb{R}$, The homogeneous $L^2$ Sobolev space is defined by 
$$
\widehat{H}^s = \{u \in \mathcal{S}'; \|u\|_{\widehat{H}^s}^2=
\displaystyle\int_{\mathbb{R}^n} |\xi|^{2s}|\hat{u}(\xi)| d\xi <+\infty
\}. 
$$
When $s<\frac{d}{2}$, $\widehat{H}^s$ is a Hilbert space by \cite[Proposition 1.34]{BCD}. In addition,  $\widehat{H}^s$ is continuously embedded in $L^{\frac{2n}{n-2s}}$ for $s\in [0, \frac{n}{2})$ by \cite[Theorem 1.38]{BCD}.  Similarly, we define the homogeneous $L^p$ space $\widehat{H}^s_{p}$ for $1\leq p \leq \infty$, see, \cite{Gaudin} for details. 

We define the Lorentz space based on the non-increasing rearrangement as follows. 

\begin{defi}[Distribution Function]
Let $f:\mathbb{R}^n \to \mathbb{C}$ be measurable.  
The \emph{distribution function} of $f$ is defined by
\[
d_f(\lambda) := \big| \{ x \in \mathbb{R}^n : |f(x)| > \lambda \} \big|,
\qquad \lambda > 0,
\]
where $|E|$ denotes the Lebesgue measure of $E$.
\end{defi}

\begin{defi}[Non-increasing Rearrangement]
The \emph{non-increasing rearrangement} of $f$ is
\[
f^*(s) := \inf \{ \lambda > 0 : d_f(\lambda) \leq s \}, \qquad s > 0.
\]
\end{defi}

\begin{defi}[Lorentz Space]
Let $1 \leq p < \infty$ and $1 \leq q \leq \infty$.  
The Lorentz space $L^{p,q}(\mathbb{R}^n)$ consists of measurable functions $f$ such that
\[
\|f\|_{L^{p,q}} =
\begin{cases}
\left( \displaystyle \int_0^\infty \left( t^{1/p} f^*(s) \right)^q \, \frac{ds}{s} \right)^{1/q}, & 1 \leq q < \infty, \\[2ex]
\sup_{s>0} \; t^{1/p} f^*(s), & q = \infty,
\end{cases}
\]
is finite. We note that $L^{p,p} = L^p$ and  $L^{p,\infty}$ is the weak-$L^p$ space.
\end{defi}

For $q=\infty$ and $s \geq 0$,  $L^{\infty}_{s}$ stands for the weighted  $L^\infty$ space over $\mathbb{R}^n$ 
defined by
\begin{eqnarray}
L^{\infty}_s
=
\{f\in L^\infty;\|f\|_{L^\infty_s}:=\|(1+|x|)^s f\|_{L^\infty} < +\infty \}. \nonumber
\end{eqnarray}

We use the following basic lemmata. We begin with the well-known Sobolev type inequality.
	
	\begin{lem}\label{lem2.1.} Let $n\geq 3$ and let $s\geq \left[\frac{n}{2}\right]+1.$ Then there holds the inequality 
		\begin{eqnarray}
			\|f\|_{L^{\infty}}\leq C\|\nabla f\|_{H^{s-1}}\nonumber
		\end{eqnarray}
		for $f\in H^{s}.$
	\end{lem}

\vspace{2ex}

\begin{lem}{\rm (\cite[Lemma 3.1]{FT-lin-half2})}\label{lemP_1}
		{\rm (i)} 
		Let $k$ be a positive number. 
		Then $P_{1}$ is a bounded linear operator from $L^p$ to $H^{k}_p$ for $1<p< \infty$. 
		In fact, it holds that 
		\begin{eqnarray}
			\|(-\Delta)^{k/2}P_{1}f\|_{L^p}\leq C\|f\|_{L^p}\qquad (f\in L^{p}).\nonumber
		\end{eqnarray}
		
		\vspace{1ex}
		{\rm (ii)} 
		Let $k$ be a positive number.
		Then there hold the estimates 
		$$
		\|(-\Delta)^{k/2} f_1\|_{L^p} \leq 
		C\|f_1\|_{L^p} \leq C\|f\|_{L^2}
		\ \ 
		(f\in L^2_{(1)}), 
		$$
		where $2 \leq   p< \infty$.
		
	\end{lem}

\vspace{2ex}

Indeed, {\rm (i)} is directly obtained by \cite[Lemma 3.1]{FT-lin-half2}. On {\rm (ii)}, taking a cut off function $\tilde{\chi}_1\in C^{\infty}_0(\IR^n)$ by 
\begin{align}
		 0 \leq  \tilde{\chi}_1 \leq 1, \ \ \tilde{\chi_1}(\xi) =1 \ (\xi \in Q_{r'_\infty}), \ \ 
		\tilde{\chi}_1(\xi) =0  \ (\xi \in \mathbb{R}^n \setminus Q_{2r'_\infty} )
\end{align}
and letting $\tilde{P}_1 = \mathcal{F}^{-1}\tilde{\chi}_1 \mathcal{F}$ derive that $f = \tilde{P}_1 f$ for $f \in L^p_{(1)}$. Hence a similar idea as in the proof of \cite[Lemma 3.1]{FT-lin-half2} yields  the first inequality in {\rm (ii)}. 
We note that $\mathcal F (-\Delta)^{k/2} f= |\xi|^k \widehat{f}$. The second inequality in {\rm (ii)} is directly obtained by the Sobolev embedding and the  first inequality for $p=2$.

\vspace{2ex}

\begin{lem}{\rm (\cite[Lemma 3.1]{FT-lin-half2})}\label{lemPinfty}
	{\rm (i)} 
	Let $k$ be a nonnegative integer. 
	Then $P_\infty$ is a bounded linear operator on $H^k$. 
	
	\vspace{1ex}
	{\rm (ii)} 
	There hold the inequalities 
	\begin{eqnarray*}
		\|P_\infty f\|_{L^2}
		& \leq &
		C\|\nabla f\|_{L^2} 
		\ \ (f\in H^1),
		\\[2ex]
		\|f_{\infty}\|_{L^2}
		& \leq & 
		C\|\nabla f_{\infty}\|_{L^2} 
		\ \ (f_{\infty}\in H^{1}_{(\infty)}).
	\end{eqnarray*}
\end{lem}

\vspace{2ex}

\begin{lem}{\rm (\cite[Lemma 2.1]{Choi})}\label{composition-est}
Let $f,g \in H^k \cap L^\infty$ with $\nabla f \in L^\infty$ for $k \in \mathbb{N}$. Then it holds that 
$$
\|\nabla^k(fg)- f \nabla^k g\|_{L^2} \leq 
C(\|\nabla f\|_{L^\infty}\|\nabla^{k-1}g\|_{L^2}+ 
\|\nabla^k f\|_{L^2}\|g\|_{L^\infty} ).
$$

\end{lem}

\vspace{2ex}

\begin{lem}[Sobolev inequality]\label{Hardy-Rellich}


 Let $1<p<q<\infty$ and  $\delta=\frac{n}{2}\big(\frac1p-\frac1q\big)$. 
Then for $u\in \widehat {H}^{2\delta}_p(\IR^n)$ there holds the Sobolev embedding estimate
\begin{equation}\label{Sob-ineq}
\|u\|_{L^q} \leq C\|(-\Delta)^\delta u\|_{L^p}.
\end{equation}
\end{lem}

\vspace{2ex}

\begin{lem}[H\"older Inequality on  Lorentz Spaces]\label{Holder-L}
Let $1 < p_1, p_2 < \infty$, $1 \leq q_1, q_2 \leq \infty$, $1\leq p <\infty$ and
\[
\frac{1}{p} = \frac{1}{p_1} + \frac{1}{p_2}, 
\qquad
\frac{1}{q} = \frac{1}{q_1} + \frac{1}{q_2}.
\]
If $f \in L^{p_1,q_1}$ and $g \in L^{p_2,q_2}$, then
\[
\| f g \|_{L^{p,q}} \leq C \, \|f\|_{L^{p_1,q_1}} \, \|g\|_{L^{p_2,q_2}},
\]
for some constant $C>0$ independent of $f,g$.
\end{lem}

\vspace{2ex}

\vspace{2ex}

We will apply a basic lemma which plays an important role. 


\vspace{2ex}

    \begin{lem}\label{fgh-est} 

{\rm (i)} 
 Let $f \in L^{\frac{n}{2},\infty}$, $g \in \widehat{H}^{s_1}$ and $h \in \widehat{H}^{s_2}$ with $0<s_1, s_2 < \frac{n}{2}$ and $s_1+s_2=n-1$.
 Then it holds that 
\begin{equation}
\Big| \displaystyle\int_{\mathbb{R}^n}fgh dx\Big| \leq C\|f\|_{L^{\frac{n}{2},\infty}}\|g\|_{\widehat{H}^{s_1}}\|h\|_{\widehat{H}^{s_1}}. 
\end{equation}  

\vspace{2ex}
{\rm (ii)} 
 Let $f \in L^{{n},\infty}$, $g \in \widehat{H}^{s_1}$ and $h \in \widehat{H}^{s_2}$ with $0<s_1, s_2 < \frac{n}{2}$ and $s_1+s_2=n-2$.
 Then it holds that 
\begin{equation}
\Big| \displaystyle\int_{\mathbb{R}^n}fgh dx\Big| \leq C\|f\|_{L^{n,\infty}}\|g\|_{\widehat{H}^{s_1}}\|h\|_{\widehat{H}^{s_1}}. 
\end{equation}

\end{lem}

\vspace{2ex}
\begin{proof}
{\rm (i)} Since $0<s_1, s_2< \frac{n}{2} $ and $s_1+s_2=n-1$, there exist $2<p,r$ such that $s_1= n\Big(\frac{1}{2}-\frac{1}{p}\Big)$ and $s_2= n\Big(\frac{1}{2}-\frac{1}{r}\Big)$ with $\frac{1}{p}+\frac{1}{r}=\frac{1}{n}$. Then 
the H\"older inequality in the Lorentz spaces, Lemma \ref{Holder-L}  derives that 
\begin{align*}
\Big| \displaystyle\int_{\mathbb{R}^n}fgh dx\Big|
&\leq C\|f\|_{L^{\frac{n}{2}, \infty}} \|g\|_{L^{p, 2}}\|h\|_{L^{r,2}}, 
\end{align*} 
where $1=\frac{2}{n}+\frac{1}{p}+\frac{1}{r}$. Hence  we see from the Sobolev inequality, Lemma \ref{Hardy-Rellich} that 
$$
 \|g\|_{L^{p, 2}}\|h\|_{L^{r,2}} \leq 
 C\|(-\Delta)^\frac{s_1}{2}g\|_{L^{2}} \|((-\Delta)^\frac{s_2}{2}h\|_{L^{2}}.  
$$
This yields Lemma \ref{fgh-est} {\rm (i)}.  \rm{(ii)} is analogous to {\rm (i)}. We take $2<p,r$ such that $s_1= n\Big(\frac{1}{2}-\frac{1}{p}\Big)$ and $s_2= n\Big(\frac{1}{2}-\frac{1}{r}\Big)$ with $\frac{1}{p}+\frac{1}{r}=\frac{2}{n}$. We note that since $s_1+s_2=n-2$, we can take $2<p,r$. 
\end{proof}

\vspace{2ex}

In $\mathbb{R}^n$ $(n\geq 3)$, the  linear problem of CNS without any perturbation terms is written as follows. 
\begin{eqnarray}
\left\{
\begin{array}{lll}
\partial_{t}\rho +\gamma\div  v=0,\\
\partial_{t}v-\alpha \Delta v-\beta\nabla\div v+\gamma \nabla \rho=0,\\
\rho|_{t=0}=\rho_0, \ \ v|_{t=0}=v_0,
\label{cns-linear}
\end{array}
\right.
\end{eqnarray}
where $\alpha, \beta$ and $\gamma$ are positive constants given in \eqref{ns2} and \eqref{(a)} below. 

The resolvent problem of \eqref{cns-linear} is given as follows. 
\begin{eqnarray}\label{resolvent}
\left\{
\begin{array}{lll}
\lambda\rho +\gamma\div  v=f_1,\\
\lambda v-\alpha\Delta v-\beta\nabla\div v+\gamma \nabla \rho=f_2. 
\label{cns-linear-resolvent}
\end{array}
\right.
\end{eqnarray}
Direct computations derive an explicit formula of the solution to \eqref{resolvent}, see \cite[(2.6), (2.9)]{Kobayashi-Shibata}: 
\begin{eqnarray}
\hat{\rho}(t,\xi)&=&\Big(\frac{\lambda+(\alpha+\beta)|\xi|^2}{(\lambda-\lambda_+)(\lambda-\lambda_-)}\Big)\hat{f_1}-\displaystyle\frac{i\gamma \xi\cdot \hat{f}_2}{(\lambda-\lambda_+)(\lambda-\lambda_-)},\nonumber\\
\hat{v}(t,\xi)&=&\displaystyle\frac{1}{\lambda+\alpha|\xi|^2}\Big(I-\frac{\xi\xi\cdot}{|\xi|^2}\Big)\hat{f}_2 -\displaystyle\frac{i\gamma \xi\hat{f}_1}{(\lambda-\lambda_+)(\lambda-\lambda_-)}\nonumber\\ 
&&\quad +\displaystyle\frac{\lambda \xi \xi\cdot \hat{f}_2}{(\lambda-\lambda_+)(\lambda-\lambda_-)|\xi|^2}, 
\label{solution formula1}
\end{eqnarray} 
where $\lambda_{\pm}=-\dfrac{1}{2}(\alpha+\beta \BLACK )|\xi|^{2}\pm\dfrac{1}{2}\sqrt{( \alpha+\beta \BLACK )^{2}|\xi|^{4}-4\gamma^{2}|\xi|^{2}}$, and $|\xi | \neq 0, \dfrac{2\gamma}{\alpha+\beta}$.  We note that due to the hyperbolic aspect of the system, the spectral $\lambda_{\pm}$ curve is the parabolic curve $z= -\dfrac{\alpha+\beta}{2}|\xi|^2 \pm i\gamma |\xi|$, more rigorously, $z= -\dfrac{\alpha+\beta}{2}|\xi|^2 \pm i(\gamma |\xi| - z_1(\xi)) $, where $0< z_1$, $z_1(\xi)=O(|\xi|^3)=o(|\xi|^2)$ as $|\xi| \rightarrow 0$.

We will use the explicit solution formula in the low frequency part.   


\end{section}

\begin{section}{Formulation}
Substituting 
$\phi=\frac{\rho-\rho_{\ast}}{\rho_{\ast}}$ and $w=\frac{v}{\gamma}$ 
with $\gamma=\sqrt{p'(\rho_{\ast})}$ into \eqref{(1.1)}, 
we see that \eqref{(1.1)} is rewritten as 
\eqn{ns2}
$$
\delt u+Au=G(u,g),
$$
where 
\begin{eqnarray}\label{(a)}
A=\begin{pmatrix}
0 &\gamma\mathrm{div}\\
\gamma\nabla &-\alpha\triangle-\beta\nabla\mathrm{div}
\end{pmatrix}
,\quad \alpha =\frac{\nu}{\rho_{\ast}},\quad \beta=\frac{\nu+\nu'}{\rho_{\ast}},
\end{eqnarray}
$u=\trans(\phi, w)$ and
\begin{eqnarray}
G(u,g)&=&\begin{pmatrix}
f^{0}(u)\\
\tilde{f}(u,g)
\end{pmatrix}
,\label{(c)}\\
f^{0}(u)
&=&
-\gamma w \cdot \nabla \phi -\gamma\phi\mathrm{div}w,\label{(d)}
\\
\tilde{f}(u,g)
&=&
-\gamma (w\cdot \nabla w) 
- \nabla (p^{(1)}(\phi)\phi^{2})
- h^{(1)}(\phi) \phi \nabla (p^{(1)}(\phi)\phi^{2}) 
+\frac{g}{\gamma} \notag\\
&&\quad + \frac{\alpha}{\gamma}  h^{(1)}(\phi) \phi \Delta w + \frac{\beta}{\gamma} h^{(1)}(\phi)\phi \nabla \div w -\gamma h^{(1)}(\phi) \phi  \nabla \phi, 
\label{(e)}
\\
p^{(1)}(\phi)
&=&
\frac{\rho_{\ast}}{\gamma}
\int_{0}^{1}(1-\theta)p''(\rho_{\ast}(1+\theta\phi))d\theta,\\ 
h^{(1)}(\phi)
&=&
\int_{0}^{1}h' (1+\theta\phi) d\theta,  \ \ 
h(t)=\dfrac{1}{t}. 
\nonumber
\end{eqnarray}
Let $\bar{\phi}=\phi-\phi_\omega$, $\bar{w}=w-w_\omega$ for $\rho=\rho_* +\rho_* \phi$ and $\rho_\omega=\rho_* +\rho_* \phi_\omega$. Substituting to \eqref{(a)} and replacing notations $\bar{\phi}$ and $\bar{w}$ to $\phi$ and $w$ respectively, we have that 
\eqn{ns3}
$$
\delt u+Au +B[u_\omega]u=\mathcal{G}(u,u_\omega),
$$
where 
\begin{eqnarray}\label{(b2)}
B[u_\omega]u=\gamma\begin{pmatrix}
 b_{11}(u_\omega, u)+ b_{12}(u_\omega, \nabla u) \\
b_{21}(u_\omega, u)+ b_{22}(u_\omega, \nabla u)+ b_{23}(u_\omega, \nabla^2 u)   
\end{pmatrix}
\end{eqnarray}
\begin{align*}
b_{11} (u_\omega, u)&=   w \cdot \nabla  \phi_{\omega} + \phi \div w_{\omega}, \\
b_{12}(u_\omega, \nabla u) &= \phi_{\omega} \div w +w_{\omega} \cdot \nabla \phi, \\
b_{21}(u_\omega, u) &= 
w \cdot \nabla  w_{\omega}  
-  \dfrac{\alpha}{\gamma^2} h^{(1)}(\phi_\omega) \phi \Delta  w_{\omega} - \dfrac{\beta}{\gamma^2} h^{(1)}(\phi_\omega)\phi \nabla \div  w_{\omega}\\
&\quad 
+h^{(1)}(\phi_\omega) \phi \nabla  \phi_{\omega} 
+ 2p^{(1)}(\phi_\omega) \phi \nabla  \phi_{\omega}, \\
b_{22}(u_\omega, \nabla u)&= w_{\omega} \cdot \nabla w+h^{(1)}(\phi_\omega)\phi_{\omega} \nabla \phi + 2p^{(1)}(\phi_\omega) \phi_\omega  \nabla  \phi, \\
b_{23}(u_\omega, \nabla^2 u) & =- \dfrac{\alpha}{\gamma^2} h^{(1)}(\phi_\omega)\phi_{\omega} \Delta w - \dfrac{\beta}{\gamma^2} h^{(1)}(\phi_\omega)\phi_{\omega} \nabla \div w 
\end{align*}
and 
\begin{eqnarray}
\mathcal{G}(u,u_\omega)=G(u, g)- G(u_\omega,g)= \begin{pmatrix}
\mathcal{G}^0(u), \\
\tilde{\mathcal{G}}(u,u_\omega)
\end{pmatrix}
.
\end{eqnarray}
Here 
\begin{align}
\begin{aligned}
   \mathcal{G}^0(u)& = -\gamma \phi \div w -\gamma w \cdot \nabla \phi, \\
    \tilde{\mathcal{G}}(u,u_\omega) &= -\gamma (w\cdot \nabla w) 
- \nabla (p^{(1)}(\phi+\phi_\omega)\phi^{2})
\\
&\quad + \frac{\alpha}{\gamma}  h^{(1)}(\phi+\phi_\omega) \phi \Delta w + \frac{\beta}{\gamma} h^{(1)}(\phi+\phi_\omega)\phi \nabla \div w -\gamma h^{(1)}(\phi+\phi_\omega) \phi  \nabla \phi \\
   &\quad  + \frac{\alpha}{\gamma}h^{(2)}(\phi, \phi_\omega)
    \phi \phi_\omega\Delta w_\omega + \frac{\beta}{\gamma} h^{(2)}(\phi, \phi_\omega)
    \phi \phi_\omega\nabla \div  w_\omega\\
    &\quad - \gamma h^{(2)}(\phi, \phi_\omega) \phi \phi_\omega   \nabla \phi_\omega - 
    \nabla (p^{(2)}(\phi, \phi_\omega) \phi \phi_\omega^2) 
    - 2[\nabla, p^{(1)}(\phi+\phi_\omega)] \phi \phi_\omega\\
    &\quad - ( h^{(1)}(\phi+\phi_\omega) (\phi+\phi_\omega )  \nabla (p^{(1)}(\phi+\phi_\omega)(\phi+\phi_\omega)^{2}) \\
   &\qquad  - h^{(1)}(\phi_\omega) \phi_\omega   \nabla (p^{(1)}(\phi_\omega)\phi_\omega^{2}), \\ 
&\quad +\dfrac{\alpha}{\gamma} h^{(2)}(\phi,\phi_\omega) \phi^2 \Delta  w_{\omega} + \dfrac{\beta}{\gamma} h^{(2)}(\phi, \phi_\omega)\phi^2 \nabla \div  w_{\omega}\\
&\quad  -\gamma h^{(2)}(\phi, \phi_\omega) \phi^2 \nabla  \phi_{\omega} 
- 2 \gamma p^{(2)}(\phi, \phi_\omega) \phi^2 \nabla  \phi_{\omega} \\
&\quad -\gamma h^{(2)}(\phi, \phi_\omega)\phi_{\omega} \phi \nabla \phi - 2\gamma p^{(2)}(\phi, \phi_\omega) \phi_\omega  \phi \nabla  \phi,\\
&\quad +\dfrac{\alpha}{\gamma} h^{(2)}(\phi, \phi_\omega)\phi_{\omega} \phi \Delta w + \dfrac{\beta}{\gamma} h^{(2)}(\phi, \phi_\omega)\phi_{\omega} \phi\nabla \div w,  
\end{aligned}
\end{align}
where 
$$
h^{(2)}(\phi, \phi_\omega)= 
\displaystyle\int_0^1 h^{(1)}(\phi_\omega + \theta \phi) d\theta,  \ \ 
p^{(2)}(\phi, \phi_\omega)= 
\displaystyle\int_0^1 p^{(1)}(\phi_\omega + \theta \phi) d\theta.
$$

Applying $P_j$ $(j=1,\infty)$ to \eqref{ns3}, we decompose the system to the low frequency part and high frequency part respectively. 
The low frequency part $u_1$ satisfies  
\begin{eqnarray}\label{equation-low}
\delt u_1+A u_1 +B_1[u_\omega]u_1= \mathcal{G}_1(u,u_\omega),
\end{eqnarray}
where $u_1=P_1 u$, $u_\infty=P_\infty u$, $u=u_1+u_\infty$, 
\begin{eqnarray}\label{(b2-low)}
B_1[u_\omega]u_1=\gamma P_1 \begin{pmatrix}
 b_{11}(u_\omega, u_1)+ b_{12}(u_\omega, \nabla u_1) \\
b_{21}(u_\omega, u_1)+ b_{22}(u_\omega, \nabla u_1)+ b_{23}(u_\omega, \nabla^2 u_1) 
\end{pmatrix}
\end{eqnarray}
and 
$$
\mathcal{G}_1(u,u_\omega)=P_1\mathcal{G}(u_\omega,  u) -B_1[u_\omega, u_\infty].
$$
The high frequency part $u_\infty$ is a solution to  
\begin{eqnarray}\label{equation-high}
\delt u_\infty+A u_\infty +B_\infty[u_\omega,u]u_\infty= \mathcal{G}_\infty(u,u_\omega),
\end{eqnarray}
where  
\begin{eqnarray*}
B_\infty[u_\omega, \tilde{u}]u_\infty&=&\gamma P_\infty \begin{pmatrix}
 b_{11}(u_\omega, u_\infty)+ b_{12}(u_\omega, \nabla u_\infty) \\
b_{21}(u_\omega, u_\infty)+ b_{22}(u_\omega, \nabla u_\infty)+ b_{23}(u_\omega, \nabla^2 u_\infty) 
\end{pmatrix}
+
\mathcal{B}[\tilde{u}]u_\infty
\end{eqnarray*}
for $$\mathcal{B}[\tilde{u}]u = 
 \gamma P_\infty   \begin{pmatrix}
\tilde{w} \cdot \nabla \phi\\
0
\end{pmatrix}
$$
and 
$$
\mathcal{G}_\infty(u,u_\omega)=P_\infty\mathcal{G}(u_\omega,  u) -B_\infty[u_\omega, u]u_1 - \mathcal{B}[\tilde{u}]u_1. 
$$
Conversely, if there exists a solution $(u_1, u_\infty)$ to \eqref{equation-low}-\eqref{equation-high}, adding \eqref{equation-low} to \eqref{equation-high} and letting $u=u_1+u_\infty$, $u$ is a solution to \eqref{ns2}. Hence we look for a fixed point $u_1$ and $u_\infty$ of \eqref{equation-low} and \eqref{equation-high}. 

\vspace{2ex}

On the low frequency part, the linearlized problem is written as follows. 
\begin{eqnarray}
\left\{
\begin{array}{lll}
\partial_{t}u_1 +L_1 u_1=F_1,\\
u|_{t=0}=u_{01}
\label{cns-linear-low-per2}, 
\end{array}
\right.
\end{eqnarray}
where  $L= A_1+B_1[u_{\omega}]$ and $A_1=A|_{L^2_{(1)}}$.  

On the high frequency part, we consider the following linear system 
\begin{eqnarray}
\left\{
\begin{array}{lll}
\partial_{t}u_\infty +L_\infty u_\infty=F_\infty,\\
u_\infty|_{t=0}=u_{0\infty}
\label{cns-linear-high-per}, 
\end{array}
\right.
\end{eqnarray}
where  $L_\infty= A_\infty+B_\infty[u_{\omega}, \tilde{u}]$,  $A_\infty=A|_{L^2_{(\infty)}}$ 
and 
$\tilde{u}$ is a given function in $C([0,T]; H^s)$, $s\geq \big[\frac{n}{2}\big]+2$.

\end{section}

\begin{section}{Estimates for the low frequency part}
 In this section, we suppose that 
 $|\xi| \leq r_\infty < \frac{2\gamma}{\alpha+\beta}$, where $r_\infty=2r'_\infty$ and $r'_\infty$ is used in $L^p_{(1)}$.  
 We first consider a solution in the low frequency part $u_1=\trans(\phi_1, w_1)$ to the linearized system without perturbation terms. 
\begin{eqnarray}
\left\{
\begin{array}{lll}
\partial_{t}\phi_1 +\gamma\div  w_1=0,\\
\partial_{t}w_1-\alpha \Delta w_1-\beta\nabla\div w_1+\gamma \nabla \phi_1=0,\\
\phi_1|_{t=0}=\phi_{01}, \ \ w_1|_{t=0}=w_{01}.
\label{cns-linear-low}
\end{array}
\right.
\end{eqnarray}

The resolvent problem is given by 
\BLACK

\begin{eqnarray}\label{cns-linear-low-resolvent}
\left\{
\begin{array}{lll}
\lambda\phi_1 +\gamma\div  w_1=f_{11},\\
\lambda w_1-\alpha\Delta w_1-\beta\nabla\div w_1+\gamma \nabla \phi_1=f_{12}. 
\end{array}
\right.
\end{eqnarray}

 Let 
\begin{eqnarray}
\hat{A}_{1}=\begin{pmatrix}
0 &i\gamma\trans{\xi}\\
i\gamma\xi &\alpha|\xi|^{2}I_{n}+\beta\xi \trans{\xi}
\end{pmatrix}
\ \ \ (\xi\in \mathbb{R}^n). 
\nonumber
\end{eqnarray}
Then, ${\cal F}(A_1 u_1)=\hat A_{1}\hat u_1$. 
Hence, if ${\rm supp}\,\hat{u}_{1}\subset Q_{r'_\infty},$ then 
${\rm supp}\,\hat{A}_{1}\hat{u}_{1}\subset Q_{r'_\infty}.$ 
\BLACK 
Furthermore, 
we see  from the property of the low frequency Lemma \ref{lemP_1} that
$$
\|A_1 u_{1}\|_{L^2}\leq C\|u_{1}\|_{L^2},
$$
for $u_1 \in L^2_{(1)}$. 
Therefore, $A_{1}$ is a bounded linear operator on $L^2_{(1)}$. 
It then follows that 
$-A_1$ generates a uniformly continuous semigroup $S(t)=e^{-tA_1}$ 
that is given by 
$$
S(t)u_1=\mathcal{F}^{-1}(e^{-t \hat{A}_{1}}\mathcal{F}u_1)
\ \ \ (u_1\in L^p_{(1)}).
$$
Let $\sigma (-A_1)$ and $\sigma (-\hat{A}_1)$ be spectral sets of $-A_1$ and $-\hat{A}_1$ on $L^2_{(1)}$ respectively. 
We state the following results about the spectrum and the spectral resolution:

\vspace{2ex}

\begin{lem}\label{spectral-resolution} {\rm (\cite{Matsumura-Nishida1})} 
{\rm (i)} 
The set of all eigenvalues of $-\hat{A}_{1}$ consists of $\lambda_{j}(\xi)\,(j=1,\pm)$, 
where 
\BLACK
\begin{eqnarray*}
\left\{
\begin{array}{ll}
\lambda_{1}(\xi)=-\alpha|\xi|^{2},\\
\lambda_{\pm}(\xi)=-\frac{1}{2}( \alpha+\beta \BLACK )|\xi|^{2}\pm\frac{1}{2}\sqrt{( \alpha+\beta \BLACK )^{2}|\xi|^{4}-4\gamma^{2}|\xi|^{2}}.
\end{array}
\right. 
\end{eqnarray*}
If $|\xi|<\frac{2\gamma}{ \alpha+\beta \BLACK }$, then 
$$
{\rm Re}\,\lambda_{\pm}
=-\frac{1}{2}( \alpha+\beta \BLACK)|\xi|^{2}, 
\ \ \ 
{\rm Im}\, \lambda_{\pm}
=\pm\gamma|\xi|\sqrt{1-\frac{( \alpha+\beta \BLACK )^2}{4\gamma^2}|\xi|^2}.
$$

{\rm (ii)} 
If $|\xi|<\frac{2\gamma}{ \alpha+\beta \BLACK }$, then $e^{-t\hat{A}_{1}}$ has the spectral resolution 
\begin{eqnarray*}
e^{-t\hat{A}_{\xi}}=\sum_{j=1,\pm}e^{t\lambda_{j}(\xi)}{\Pi}_{j}(\xi),
\end{eqnarray*}
where ${\Pi}_{j}(\xi)$ are eigenprojections for $\lambda_{j}(\xi)\,(j=1,\pm)$, and $\Pi_{j}(\xi)\,(j=1,\pm)$ satisfy
\[
\Pi_{1}(\xi)=\begin{pmatrix}
0 &0\\
0 &I_{n}-\frac{\xi\trans{\xi}}{|\xi|^{2}}
\end{pmatrix},
\ \ \ 
\Pi_{\pm}(\xi)
=\pm\frac{1}{\lambda_{+}-\lambda_{-}}\begin{pmatrix}
-\lambda_{\mp} &-i\gamma\trans{\xi}\\
-i\gamma\xi &\lambda_{\pm}\frac{\xi\trans{\xi}}{|\xi|^{2}}
\end{pmatrix}
.
\]
Furthermore, for $0<r_\infty<\frac{2\gamma}{ \alpha+\beta \BLACK }$, 
there exists a constant $C= C(r_\infty)\BLACK>0$ such that the estimates
\begin{eqnarray}
\|\Pi_{j}(\xi)\|\leq C\quad (j=1,\pm) 
\label{Pi}
\end{eqnarray}
hold for $|\xi|\leq r_{\infty}$.
\end{lem}

\vspace{2ex}





\vspace{2ex}






\vspace{2ex}

We have the following resolvent estimate on each $\xi$ which plays  an important role.  

\vspace{2ex}

\begin{prop}\label{resolvent2-each-point}
Let $c_0>0$ be a fixed positive number. 
We suppose that $\alpha+\beta \geq 5\gamma^2$ with $\gamma \geq 1$. 
Let $\lambda \in \mathbb{C}\setminus \sigma (-\hat{A}_1)$.  Let 
 resolvent sets \BLACK $\mathcal R_{1,\pm}$ be defined as 
\begin{align*}
\mathcal{R}_{1, \pm } & =\{\lambda \in \mathbb{C}\setminus \sigma (-\hat{A}_1); 
\Re \lambda <0,  \   \lambda = -a^2 \pm (a + c_0) i \mbox{ with } a>0 
\}, \\
%
\mathcal{R}_2 & = \{\lambda \in \mathbb{C}\setminus \sigma (-\hat{A}_1); \Re \lambda >0\}.
\end{align*}
%
%
\BLACK 
Then 
 for $\lambda\in \mathcal R_{1,\pm}\cup \mathcal R_2$ \BLACK
the solution $\hat{u}_1(\xi)=\trans(\hat{\phi}_1, \hat{w}_1) \in L^2_{(1)}$ to \eqref{cns-linear-low-resolvent} satisfies 
for $F_1 =\trans(f_{11}, f_{12} )$ the estimate 
$$
|\hat{\phi}_1| + |\hat{w}_1| \leq \dfrac{C}{|\lambda||\xi|}(|\hat{f}_{11}|+|\hat{f}_{12}|)
$$
and 
$$
||\xi|^2 \hat{\phi}_1| + ||\xi|^2\hat{w}_1| \leq C(|\hat{f}_{11}|+|\hat{f}_{12}|), 
$$
where positive constants $C$ are independent of $c_0$, $\lambda$ and $\xi$.  
\end{prop}

\vspace{2ex}

\begin{proof} The proof  is split into two cases:\BLACK 

Case $(a)$: $|\lambda| < \dfrac{\gamma}{\sqrt{2}} |\xi| \ \textrm{ or } \  \sqrt{2}\gamma|\xi|< |\lambda|,$
\vspace{1ex}

Case $(b)$: $\dfrac{\gamma}{\sqrt{2}}|\xi|\leq |\lambda| \leq  \sqrt{2}\gamma|\xi|$. 

{\bf Case (a):}  In this case we consider $\lambda\in  \mathcal R_{1,\pm}\cup \mathcal R_2$ simultaneously \BLACK 
 \BLACK  and use the explicit formula \eqref{solution formula1}. 
Here $(\lambda-\lambda_+)(\lambda-\lambda_-)= \lambda^2+(\alpha+\beta)\lambda |\xi|^2 +\gamma^{ 2\BLACK} |\xi|^2$. 

 First let $\lambda$ satisfy $|\lambda| < \dfrac{\gamma}{\sqrt{2}}|\xi|$. \BLACK We note that 
\begin{align}\label{resolvent-est-hard1}
\lambda^2+(\alpha+\beta)\lambda |\xi|^2 +\gamma^2 |\xi|^2 \geq 
C|\xi|^2,  
\end{align}
as $|\xi| \rightarrow 0$. 
 Indeed, when $|\lambda| < \dfrac{\gamma}{\sqrt{2}}|\xi|$, 
 then  $|\lambda^2 + \gamma^2 |\xi|^2| \geq \dfrac{\gamma^2}{2}|\xi|^2$, 
 and the remainder $(\alpha+\beta)\lambda |\xi|^2 = O(|\xi|^3) = o(|\xi|^2)$ 
 can be ignored in (4.5). 
\BLACK 
In this case we get 
that 
\begin{align}\label{resolvent-est-hard-10}
\Big|\frac{ \lambda}{(\lambda-\lambda_+)(\lambda-\lambda_-)}\Big| \leq \dfrac{C}{|\lambda|},  \ \  \Big|\frac{ |\xi|}{(\lambda-\lambda_+)(\lambda-\lambda_-)}\Big| \leq \dfrac{C}{|\lambda|}
\end{align}
and 
\begin{align}\label{resolvent-est-hard-11}
\Big|\frac{ \lambda|\xi|^2}{(\lambda-\lambda_+)(\lambda-\lambda_-)}\Big| \leq C,  \ \  \Big|\frac{ |\xi|^3}{(\lambda-\lambda_+)(\lambda-\lambda_-)}\Big| \leq C
\end{align}
uniformly for $|\xi|\ll1$ 
\BLACK. 

Second, on $\sqrt{2}\gamma|\xi|<|\lambda|$ for 
$$
\Big|\frac{ \lambda|\xi|^2}{(\lambda-\lambda_+)(\lambda-\lambda_-)}\Big|
= \Big|\frac{ |\xi|^2}{\lambda + (\alpha+\beta)|\xi|^2 + \gamma^2 \frac{|\xi|^2}{\lambda}}\Big|, 
$$
we note that 
$$
\dfrac{\gamma^2 |\xi|^2}{|\lambda|} \leq \dfrac{\gamma}{\sqrt{2}}|\xi|. 
$$
Hence for  $|\xi|\ll 1$ we see that 
$$
\Big|\lambda + (\alpha+\beta)|\xi|^2 + \gamma^{ 2\BLACK} \frac{|\xi|^2}{\lambda}\Big| \geq \sqrt{2}\gamma |\xi|- (\alpha+\beta)|\xi|^2 -\dfrac{\gamma}{\sqrt{2}}|\xi| \geq 
c|\xi|.
$$
Therefore it holds that 
\begin{align}\label{resolvent-est-hard-22-1}
\Bigg|\frac{ |\xi|^2}{\lambda + (\alpha+\beta)|\xi|^2 + \gamma^2 \frac{|\xi|^2}{\lambda}}\Bigg| \leq C, 
\end{align}
as well as 
\begin{align}\label{resolvent-est-hard-22-2}
\Big|\frac{ |\xi|^3}{\lambda^2 + (\alpha+\beta)\lambda|\xi|^2 + \gamma^2 |\xi|^2}\Big| \leq C
\end{align}
by the condition $\sqrt{2}\gamma|\xi|<|\lambda|$. Moreover, since $|\xi|^2 \leq  \dfrac{r_\infty |\lambda|}{\sqrt{2}\gamma} \leq \dfrac{|\lambda|}{3(\alpha+\beta)}$ for $|\xi| \leq r_\infty \ll 1$, 
it holds that 
$$
|\lambda^2 +  (\alpha+\beta)\lambda |\xi|^2 + \gamma^2 |\xi|^2| \geq |\lambda( \lambda+(\alpha+\beta) |\xi|^2)| - \dfrac{1}{2} |\lambda|^2 \geq \dfrac{|\lambda|^2}{6}. 
$$
Hence we see that 
\begin{align}\label{resolvent-est-hard-22-01}
\Big|\frac{ \lambda}{(\lambda-\lambda_+)(\lambda-\lambda_-)}\Big| \leq \dfrac{C}{|\lambda|}, 
\end{align}
as well as
\begin{align}\label{resolvent-est-hard-22-02}
\Big|\frac{ |\xi|}{(\lambda-\lambda_+)(\lambda-\lambda_-)}\Big| \leq \dfrac{C}{|\lambda|} 
\end{align}
by the condition $\sqrt{2}\gamma|\xi|<|\lambda|$.  Especially, we note that when $|\lambda| \geq 1$, by the above proof we do not need  the form $\lambda = -a^2 \pm (a+c_0) i$ for $\Re \lambda <0$. \BLACK

{\bf Case (b):} Let $\dfrac{\gamma}{\sqrt{2}}|\xi|\leq |\lambda| \leq  \sqrt{2}\gamma|\xi|$. We shall show for the sets $\mathcal R_{1,\pm}$ and $\mathcal R_2$ separately \BLACK that   
\begin{align}\label{resolvent-est-hard2}
\Big|\frac{ \lambda}{(\lambda-\lambda_+)(\lambda-\lambda_-)}\Big| \leq  \dfrac{C}{|\lambda||\xi|},  \ \ \Big|\frac{ \lambda|\xi|^2}{(\lambda-\lambda_+)(\lambda-\lambda_-)}\Big| \leq  C ,
\end{align}
and 
\begin{align}\label{resolvent-est-hard3}
\Big|\frac{|\xi|}{(\lambda-\lambda_+)(\lambda-\lambda_-)}\Big| \leq  \dfrac{C}{|\lambda||\xi|},  \ \ \Big|\frac{ |\xi|^3}{(\lambda-\lambda_+)(\lambda-\lambda_-)}\Big| \leq C
\end{align}
where $C$ are independent of $\lambda$ and $|\xi|\ll 1$. 
Since $|\lambda|\sim|\xi|$, it suffices to prove $\eqref{resolvent-est-hard2}_1$.

{\bf Set $\mathcal R_{\bf 1,\pm}$:} \BLACK We note that by taking $ r_\infty \leq \dfrac{1}{2\sqrt{2}\gamma}$, we see from  the low frequency part condition $|\xi| \leq r_\infty$ that 
$$
|\Im \lambda | \leq |\lambda| \leq 
 \sqrt{2}\gamma |\xi| \leq \dfrac{1}{2}. 
$$ 
Since $\Re \lambda < 0$, $\lambda = -a^2 \pm (a+c_0) i$ for $a>0$ and $c_0>0$, it holds that 
$$
a \leq |\Im \lambda|. 
$$
\BLACK 
Then 
$
|\Re \lambda|=a^2 \ll |\Im \lambda|, 
$
especially, $
|\Re \lambda| <  \frac{1}{2}|\Im \lambda|, 
$
$$
\dfrac{\gamma}{\sqrt{2}}|\xi| \leq |\Im \lambda|+ |\Re \lambda|  \leq  2 |\Im \lambda|.
$$
\BLACK
Hence we see that 
\begin{align}
\dfrac{\gamma}{2\sqrt{2}}|\xi| \leq |\Im \lambda| \leq \sqrt{2}\gamma|\xi|. 
\end{align}
In addition, 
 \begin{align}\label{lambdacondition}
 |\Re \lambda| \leq 2\gamma^2|\xi|^{ 2\BLACK}. 
\end{align}
\BLACK We rewrite  $\dfrac{ |\xi|}{(\lambda-\lambda_+)(\lambda-\lambda_-)}$ to 
\begin{align}\label{rewrite-resolvent}
\frac{ |\xi|}{(\lambda-\lambda_+)(\lambda-\lambda_-)} = \frac{ |\xi|}{\lambda_+-\lambda_-}\Big\{\dfrac{1}{\lambda-\lambda_+}-\dfrac{1}{\lambda-\lambda_-}\Big\}. 
\end{align}
Here  $\lambda_+-\lambda_- = \frac{1}{2}\sqrt{(\alpha+\beta)^2|\xi|^4  - 4\gamma^2|\xi|^2} \sim \gamma |\xi|i$ as $|\xi| \ll 1$. Since 
$$
\lambda -\lambda_+ =  \Re \lambda + (\Im \lambda) i + \frac{\alpha+\beta}{2}|\xi|^2   + i\gamma|\xi|\sqrt{1-\frac{( \alpha+\beta  )^2}{4\gamma^2}|\xi|^2} \,, 
$$
\BLACK
this together with \eqref{lambdacondition} verifies that if $\alpha+\beta \geq 5\gamma^2$ 
\begin{align*}
    |\lambda -\lambda_+| 
    &\geq |\Re (\lambda -\lambda_+)| \geq 
    \dfrac{(\alpha+\beta)}{2}|\xi|^2- |\Re \lambda| \\
    & \geq \dfrac{(\alpha+\beta)}{2}|\xi|^2 - 2\gamma^2 |\xi|^2 \geq C_2|\xi|^2, 
\end{align*}
 where $C_2$ is independent of $c_0$, $\lambda$ and $\xi$. \BLACK 
Similarly, 
$
|\lambda -\lambda_-| \geq C_2|\xi|^2.
$
Substituting these inequalities to \eqref{rewrite-resolvent} we obtain that 
\begin{align}\label{rewrite-resolvent2}
\Big|\frac{ |\xi|}{(\lambda-\lambda_+)(\lambda-\lambda_-)}\Big| \leq C\frac{1}{|\xi|^{ 3}} \leq C \frac{1}{|\lambda||\xi|} 
\end{align}
uniformly for  $c_0,$  $\lambda$ and $\xi$. \BLACK  Hence \eqref{resolvent-est-hard2} and \eqref{resolvent-est-hard3} are completely proved. \BLACK 

\BLACK

{\bf Set $\mathcal R_{\bf 2}$:} \BLACK For $\Re \lambda >0$,  we  rewrite  $\dfrac{ \lambda }{(\lambda-\lambda_+)(\lambda-\lambda_-)}$ to 
\begin{align}\label{rewrite-resolvent2}
\frac{ \lambda}{(\lambda-\lambda_+)(\lambda-\lambda_-)} = \frac{ 1}{\lambda + (\alpha+\beta)|\xi|^2 + \dfrac{\gamma^{2\BLACK} |\xi|^2}{\lambda}}. 
\end{align}
\BLACK We note that 
$$
\lambda + (\alpha+\beta)|\xi|^2 + \dfrac{\gamma^{ 2\BLACK} |\xi|^2}{\lambda} = 
\Re \lambda + (\Im \lambda)i + (\alpha+\beta)|\xi|^2 + \dfrac{\gamma^2|\xi|^2(\Re \lambda- (\Im \lambda)i)}{|\lambda|^2}.
$$
Hence 
$$
\Big|\lambda + (\alpha+\beta)|\xi|^2 + \dfrac{\gamma^{ 2\BLACK} |\xi|^2}{\lambda} \Big| \geq  
\Big|\Re\Big( \lambda + (\alpha+\beta)|\xi|^2 + \dfrac{\gamma^2 |\xi|^2}{\lambda} \Big)\Big| \geq (\alpha+\beta)|\xi|^2. 
$$
\BLACK 
Therefore we have that 
\begin{align}\label{rewrite-resolvent3}
\Bigg| \frac{1}{\lambda + (\alpha+\beta)|\xi|^2 + \dfrac{\gamma^2 |\xi|^2}{ \lambda\BLACK}}\Bigg| \leq \dfrac{C}{|\xi|^2} 
\leq \dfrac{C}{|\lambda||\xi|} 
\end{align}
and thus 
\begin{align}\label{rewrite-resolvent4}
\Big|\frac{ |\xi|}{(\lambda-\lambda_+)(\lambda-\lambda_-)} \Big| \leq \dfrac{C}{|\lambda||\xi|},  \ \ 
\Big|\frac{ |\xi|^3}{(\lambda-\lambda_+)(\lambda-\lambda_-)} \Big| \leq C,  
\end{align}
by the condition $\dfrac{\gamma}{\sqrt{2}}|\xi|\leq |\lambda| \leq  \sqrt{2}\gamma|\xi|$. 

In addition, the property of resolvent of the Stokes equation, cf. \cite[(22)]{ST}, derives that  for any $0<\delta<\pi$ there exists a constant $C_\delta>0$ such that \BLACK
\begin{align}\label{resolvent-est-hard4}
\Big| \frac{1}{\lambda+|\xi|^2}\Big| \leq C_\delta\frac{1}{|\lambda|+|\xi|^2} \quad\textrm{ for all } \lambda\in \Sigma_{\pi-\delta}. \BLACK
\end{align}
%


Therefore, \eqref{resolvent-est-hard-10}, \eqref{resolvent-est-hard-11}, \eqref{resolvent-est-hard-22-1}, \eqref{resolvent-est-hard-22-2}, \eqref{resolvent-est-hard-22-01}, \eqref{resolvent-est-hard-22-02}, 
\eqref{resolvent-est-hard2}, \eqref{resolvent-est-hard3} and \eqref{resolvent-est-hard4}  derive Proposition \ref{resolvent2-each-point}. 
\end{proof}

\vspace{2ex}

{\rm \begin{rem}\label{remark-resolvent2} 
Under the conditions $\alpha+\beta \geq 5 \gamma^2$ with $\gamma \geq 1$, any complex number $\lambda= -r^2 \pm (r+c_0)i$ for $r>0$ does not touch the spectral curve $z=-\frac{\alpha+\beta}{2}|\xi|^2 \pm \gamma |\xi|i$, 
more rigorously, $z= -\frac{\alpha+\beta}{2}|\xi|^2 \pm i(\gamma |\xi| - z_1(\xi) )$, where $0 < z_1$, $z_1(\xi)=O(|\xi|^3)=o(|\xi|^2)$ as $|\xi| \rightarrow 0$. 

\end{rem}}

\vspace{2ex}

In \eqref{solution formula1}, we set 

\begin{eqnarray*}
\hat{\rho}(t,\xi)&=&\Big(\frac{\lambda+(\alpha+\beta)|\xi|^2}{(\lambda-\lambda_+)(\lambda-\lambda_-)}\Big)\hat{f_1}-\displaystyle\frac{i\gamma \xi\cdot \hat{f}_2}{(\lambda-\lambda_+)(\lambda-\lambda_-)},\\
&=:& I_{\rho,1} + I_{\rho,2}, \\
\hat{v}(t,\xi)&=&\displaystyle\frac{1}{\lambda+\alpha|\xi|^2}\Big(I-\frac{\xi\xi\cdot}{|\xi|^2}\Big)\hat{f}_2 -\displaystyle\frac{i\gamma \xi\hat{f}_1}{(\lambda-\lambda_+)(\lambda-\lambda_-)}\nonumber\\ 
&&\quad +\displaystyle\frac{\lambda \xi \xi\cdot \hat{f}_2}{(\lambda-\lambda_+)(\lambda-\lambda_-)|\xi|^2}\\
 &=:& I_{v,1} + I_{v,2}. 
\end{eqnarray*} 

\vspace{2ex}

\begin{cor}\label{resovent-estimate-modify-0}
Let $c_0>0$ be a fixed positive number. 
We suppose that $\alpha+\beta \geq 5\gamma^2$ with $\gamma \geq 1$. 
Let $\lambda \in \mathbb{C}\setminus \sigma (-\hat{A}_1)$.  Let 
 resolvent sets \BLACK $\mathcal R_{1,\pm}$ be defined as in Proposition \ref{resolvent2-each-point}. 
 
{\rm (i)} Then 
 for $\lambda\in \mathcal R_{1,\pm}\cup \mathcal R_2$ \BLACK
 the following estimates hold for any small epsilon $\epsilon$ independent of $c_0$, $\lambda$ and $\xi$. 
$$
|I_{\rho, 1}| + |I_{\rho, 2}|+ |I_{v, 1}| \leq \dfrac{C}{|\lambda|^{1/2 +\epsilon}|\xi|^{3/2-\epsilon}}(|\hat{f}_{11}|+|\hat{f}_{12}|)
$$
and 
$$
||\xi|I_{\rho, 1}| + ||\xi|I_{\rho, 2}|+ ||\xi|I_{v, 1}|\leq C(|\hat{f}_{11}|+|\hat{f}_{12}|), 
$$
where positive constants $C$ are independent of $c_0$, $\lambda$ and $\xi$.  

{\rm (ii)} For $\lambda\in \mathcal R_{1,\pm}$, it holds that 
$$
 |I_{v, 2}| \leq \dfrac{C}{|\lambda|^{1/2 +\epsilon}|\xi|^{3/2-\epsilon}} |\hat{f}_{12}|
$$
and 
$$
||\xi| I_{v, 2}|\leq \dfrac{C}{|\xi|}|\hat{f}_{12}|.  
$$
\end{cor}

\vspace{2ex}

\begin{proof}
{\rm (i)} When $|\lambda| < \frac{\gamma}{\sqrt{2}}|\xi|$, we see from \eqref{resolvent-est-hard-10} that 
$$
|I_{\rho, 1}| + |I_{\rho, 2}|+ |I_{v, 1}| \leq \frac{C}{|\xi|} \leq \dfrac{C|\xi|}{|\lambda|^{1/2 +\epsilon}|\xi|^{3/2-\epsilon}}(|\hat{f}_{11}|+|\hat{f}_{12}|). 
$$
Hence since $|\xi| \leq r_0$, 
$$
|I_{\rho, 1}| + |I_{\rho, 2}|+ |I_{v, 1}|   \leq \dfrac{C}{|\lambda|^{1/2 +\epsilon}|\xi|^{3/2-\epsilon}}(|\hat{f}_{11}|+|\hat{f}_{12}|)
$$
as well as 
$$
||\xi|I_{\rho, 1}| + ||\xi|I_{\rho, 2}|+ ||\xi|I_{v, 1}|\leq C(|\hat{f}_{11}|+|\hat{f}_{12}|). 
$$
In the case $\sqrt{2}\gamma|\xi|<|\lambda|$, the estimates holds similarly by \eqref{resolvent-est-hard-22-01}-\eqref{resolvent-est-hard-22-02}.  Estimates in the case $\frac{\gamma}{\sqrt{2}} \leq |\lambda| \leq  \sqrt{2}\gamma|\xi|$ directly holds by \eqref{resolvent-est-hard2} and \eqref{resolvent-est-hard3}.

{\rm (ii)} Let $\lambda=-r^2 \pm r i$ for $r<0$. When $|\lambda| \leq \frac{1}{\sqrt{2}}|\xi|$, 
we note that $|\Im \lambda| \leq \frac{1}{\sqrt{2}}|\xi|$ and thus $|\Re \lambda | \leq \frac{1}{2}|\xi|^2$. 
We see from \eqref{resolvent-est-hard4} that 
\begin{eqnarray*}
\frac{1}{|\lambda +|\xi|^2|} &\leq& \frac{1}{|\Re (\lambda+|\xi|^2)|} \\
&\leq & \frac{1}{|\xi|^2 - |\Re \lambda|} \\
&\leq &\frac{C}{|\xi|^2} \leq \frac{C}{|\lambda|^{1/2 +\epsilon}|\xi|^{3/2-\epsilon}}. 
\end{eqnarray*}
When  $|\lambda| \geq \frac{1}{\sqrt{2}}|\xi|$, it follows from \eqref{resolvent-est-hard4} that 
\begin{eqnarray*}
\frac{1}{|\lambda +|\xi|^2|} \leq \frac{C}{|\lambda|} &=& \frac{C|\xi|^{1/2-\epsilon}}{|\lambda||\xi|^{1/2-\epsilon}}\\
&\leq  &\frac{C|\lambda |^{1/2-\epsilon}}{|\lambda||\xi|^{3/2-\epsilon}} =  \frac{C}{|\lambda|^{1/2 +\epsilon}|\xi|^{3/2-\epsilon}}.
\end{eqnarray*}
The estimate $||\xi| I_{v, 2}|$ is analogous. 
\end{proof}

Proposition \ref{resolvent2-each-point} verifies the following resolvent estimate in $L^2_{(1)}$. 

\vspace{2ex}

\begin{prop}\label{resolvent2}
Let $\lambda \in \mathbb{C}\setminus \sigma (-{A}_1)$. We suppose that $\alpha+\beta \geq 5\gamma^2$ with $\gamma \geq 1$.  Let 
resolvent sets \BLACK $\mathcal R_{1,\pm}$ be defined as 
\begin{align*}
\mathcal{R}_{1, \pm } & =\{\lambda \in \mathbb{C}\setminus \sigma (-{A}_1); 
\Re \lambda <0,  \   \lambda = -a^2 \pm (a + c_0) i \mbox{ with } a>0 
\}, \\
%
\mathcal{R}_2 & = \{\lambda \in \mathbb{C}\setminus \sigma (-{A}_1); \Re \lambda >0\}.
\end{align*}
%
%
\BLACK 
Then 
 for $\lambda\in \mathcal R_{1,\pm}\cup \mathcal R_2$ \BLACK the solution $u_1=\trans(\phi_1,w_1) \in L^2_{(1)}$ to \eqref{cns-linear-low-resolvent} satisfies 
or $F_1=\trans(f_{11}, f_{12} )\in L^2_{(1)}$ the estimate 
$$
\| (-\Delta)^{\delta} u_1\|_{L^2} 
=\| (-\Delta)^{\delta} (\lambda+A_1)^{-1} F_1\|_{L^2} 
\leq \dfrac{C_{n}}{|\lambda|^{1-\delta }}\{\|f_{11}\|_{L^2}+\|f_{12} \|_{L^2}\}
$$
for $0 \leq \delta \leq 1$, 
 where $C_n$ is independent of $c_0$ and $\lambda$ 
. \BLACK 
\end{prop}

\vspace{2ex}

\begin{proof}
We see from Proposition \ref{resolvent2-each-point} and the Plancherel theorem 
that  
$$
\| u_1\|_{L^2}\leq \dfrac{C_{n}}{|\lambda|}\{\|f_{11}\|_{L^2}+\|f_{12} \|_{L^2}\}
$$
and 
$$
\| \nabla^2 u_1\|_{L^2}\leq C_{n}\{\|f_{11}\|_{L^2}+\|f_{12} \|_{L^2}\}. 
$$
We note that  $1/|\xi|$ is integrable on $L^2_{(1)}$ for $3\leq n$ due to the cut-off property, {\rm i.e., }
$$
\Big\|\frac{1}{|\xi|}\Big\|_{L^2_{(1)}} \leq C \Big\|\frac{1}{|\xi|}\Big\|_{L^2(|\xi| \leq r_\infty)} <+\infty.  
$$

The interpolation thus derives that 
$$
\| (-\Delta)^{\delta} u_1\|_{L^2}\leq \dfrac{C_{n}}{|\lambda|^{1-\delta }}\{\|f_{11}\|_{L^2}+\|f_{12} \|_{L^2}\}
$$
for $0 \leq \delta \leq 1$. 
\end{proof}

\vspace{2ex}

{\rm \begin{rem}\label{no-multiplier}
 We cannot expect that 
the above uniform resolvent estimate 
for $\| \nabla^2 (\lambda +A_1)^{-1} f \|_{L^p}$ holds  with $1<p< 2$ because when 
$1<p<2$, the semigroup has the growth order rate in time  as shown in \cite{Kobayashi-Shibata}. Indeed,  the resolvent estimates have worse order terms for the lower exponent $1<p<2$ in Proposition \ref{resolvent2-each-point}.  In addition, Proposition \ref{resolvent2-each-point} shows that resolvent estimates in each $\xi$ have worse order, however, the $L^2$ integrals have suitable 
decay order of $\lambda$ due to the low frequency part.  

  \end{rem}
}

\vspace{2ex}

\begin{cor}\label{resolvent2-cor}
Let $\lambda \in \mathbb{C}\setminus  \sigma (-{A}_1)$ as in Proposition \ref{resolvent2}.   Then 
for $F_1=\trans(f_{11}, f_{12} )\in L^2_{(1)}$ 
$$
\| (-\Delta)^{\frac{k}{2}}(\lambda +A_1)^{-1} (-\Delta)^{1-\frac{k}{2}} F_1 \|_{L^2}\leq C_{n}\{\|f_{11}\|_{L^2}+\|f_{12} \|_{L^2}\}, 
$$
where $0 \leq k \leq 2$. 
\end{cor}

\vspace{2ex}

The estimate is analogous to  
the explicit formula \eqref{solution formula1} and the proof of Proposition \ref{resolvent2}.  We note that $\mathcal F (-\Delta)^{k/2} f= |\xi|^k \widehat{f}$. 

\vspace{2ex}

By the specific resolvent estimates Corollary \ref{resovent-estimate-modify-0} in the Fourier space, we get the corresponding resolvent estimates as follows. 

\vspace{2ex}

\begin{cor}\label{resolvent2-cor-2-modify}
 Let $\lambda\in \mathcal R_{1,\pm}$  as in Proposition \ref{resolvent2}.   Then 
for $F_1=\trans(f_{11}, f_{12} )\in L^2_{(1)}$  
$$
\| (-\Delta)^{\delta}(\lambda +A_1)^{-1}F_1|
\leq \dfrac{C_{n}}{|\lambda|^{1/2+\epsilon-\delta }}\{\|f_{11}\|_{L^2}+\|f_{12} \|_{L^2}\}
$$
for $0 \leq \delta < \frac{1}{2}$ and any small $\epsilon>0$ independent of $c_0$ and $\lambda$. In addition, 
$$ 
\| (-\Delta)^{1/2}(\lambda +A_1)^{-1}F_1|
\leq C_{n}\{\|f_{11}\|_{L^2}+\|f_{12} \|_{L^2}\}. 
$$

\end{cor}

\vspace{2ex}

Furthermore, the duality argument and the Plancherel theorem directly verify that

 \vspace{2ex}
\begin{cor}\label{resolvent2-cor-2}
Let $\lambda \in \mathbb{C}\setminus  \sigma (-{A}_1)$ as in Proposition \ref{resolvent2}.   Then 
for $F_1=\trans(f_{11}, f_{12} )\in L^2_{(1)}$  
$$
\| (-\Delta)^{\delta} (\lambda+A_1^*)^{-1} F_1\|_{L^2} 
\leq \dfrac{C_{n}}{|\lambda|^{1-\delta }}\{\|f_{11}\|_{L^2}+\|f_{12} \|_{L^2}\}
$$
for $0 \leq \delta \leq 1$. 
\end{cor}

\vspace{2ex}

\begin{cor}\label{resolvent2-cor-3-modify}
 Let $\lambda\in \mathcal R_{1,\pm}$  as in Proposition \ref{resolvent2}.   Then 
for $F_1=\trans(f_{11}, f_{12} )\in L^2_{(1)}$  
$$
\| (-\Delta)^{\delta}(\lambda +A_1^*)^{-1}F_1|
\leq \dfrac{C_{n}}{|\lambda|^{1/2+\epsilon-\delta }}\{\|f_{11}\|_{L^2}+\|f_{12} \|_{L^2}\}
$$
for $0 \leq \delta < \frac{1}{2}$ and any small $\epsilon>0$ independent of $c_0$ and $\lambda$. In addition, 
$$ 
\| (-\Delta)^{1/2}(\lambda +A_1^*)^{-1}F_1|
\leq C_{n}\{\|f_{11}\|_{L^2}+\|f_{12} \|_{L^2}\}. 
$$
\end{cor}

 \vspace{2ex}

Recall that 
\begin{eqnarray}
\left\{
\begin{array}{lll}
\partial_{t}u_1 +L_1 u_1=0,\\
u_1|_{t=0}=u_{01}
\label{cns-linear-low-per2}, 
\end{array}
\right.
\end{eqnarray}
where $L=A_1+B_1[u_{\omega}]$, for $u=\trans(\phi, w)$ and $B_1[u_\omega]$ is defined by \eqref{(b2-low)}. 
Here,  since $\lambda +L_1=(\lambda+A)(I + (\lambda+A)^{-1}B_1[u_{\omega}])$, we have the perturbed resolvent formula 
\begin{align}\label{L1-form}
(\lambda +L_1)^{-1}=\sum_{j=0}^{\infty}\{(\lambda+A_1)^{-1}B_1[u_{\omega}]\}^{j}(\lambda+A_1)^{-1}
\end{align}
or since 
$\lambda +L_1=(I + B_1[u_{\omega}](\lambda+A_1)^{-1})(\lambda+A_1)$, we also have 
\begin{align}\label{L1-form-2}
(\lambda +L_1)^{-1}=\sum_{j=0}^{\infty}(\lambda+A_1)^{-1}\{B_1[u_{\omega}](\lambda+A_1)^{-1}\}^{j}. 
\end{align}
Let $L^*_1= A^*_1+B_1^*[u_{\omega}]$ be the adjoint operator of $L^*_1$ with the adjoint $B_1^*[u_{\omega}]$ of $B_1[u_{\omega}]$.  We note that 
\begin{align}\label{L1-form-adjoint}
(\lambda +L^*_1)^{-1}=\sum_{j=0}^{\infty}\{(\lambda+A^*_1)^{-1}B^*_1[u_{\omega}]\}^{j}(\lambda+A^*_1)^{-1}
\end{align}
or 
\begin{align}\label{L1-form-adjoint-2}
(\lambda +L^*_1)^{-1}=\sum_{j=0}^{\infty}(\lambda+A^*_1)^{-1}\{B^*_1[u_{\omega}](\lambda+A^*_1)^{-1}\}^{j}. 
\end{align}

\vspace{2ex}

We prepare several fundamental estimates for the relation between $(\lambda+A_1)^{-1}$ and $B_1$. 

\vspace{2ex}

\begin{lem}\label{L1-B1}
Let $k$ satisfy $ \max\{0, 1-n/4\} < k/2 < \min\{n/4, 1\}$.  

{\rm (i)}  It holds that for $u \in L^2_{(1)}$ 
\begin{align*}
\begin{aligned}
&\|(-\Delta)^{\frac{k}{2}} (\lambda+A_1)^{-1}B_1[u_{\omega}] u \|_{L^2} 
\leq C \epsilon \|(-\Delta)^{\frac{k}{2}} u\|_{L^2},  \\
&\|(-\Delta)^{\frac{k}{2}} (\lambda+A^*_1)^{-1}B_1^*[u_{\omega}] u \|_{L^2} 
\leq C \epsilon \|(-\Delta)^{\frac{k}{2}} u\|_{L^2}. 
\end{aligned}
\end{align*}

{\rm (ii)}  There holds that for $u \in L^2_{(1)}$ 
\begin{align*}
\begin{aligned}
&\|(\lambda+A_1)^{-1}B_1[u_{\omega}] u \|_{L^2} 
\leq C \epsilon \|\nabla u\|_{L^2}, \\
&\| (\lambda+A^*_1)^{-1}B_1^*[u_{\omega}] u \|_{L^2} 
\leq C \epsilon \|\nabla u\|_{L^2}. 
\end{aligned}
\end{align*}

\end{lem}

\vspace{1ex}

\begin{proof}
We use the duality argument. We note that $(\lambda+A_1)^{-1}: L^2_{(1)} \rightarrow L^2_{(1)}$  and thus $(\lambda+A_1)^{-1}B_1[u_{\omega}]: L^2_{(1)} \rightarrow L^2_{(1)}$. 
Let $\varphi \in L^2_{(1)}$.  
Since 
$$
((-\Delta)^{\frac{k}{2}} (\lambda+A_1)^{-1}B_1[u_{\omega}] u, \varphi)_{L^2}
= (B_1[u_{\omega}] u, (\lambda+A^*_1)^{-1}(-\Delta)^{\frac{k}{2}}\varphi)_{L^2}, 
$$
the definition of $B_1$ and Lemma \ref{fgh-est} verify that for $ \max\{0, 1-n/4\} < k/2 < \min\{n/4, 1\}$, 
\begin{align}\label{perturbation-est}
|((-\Delta)^{\frac{k}{2}} (\lambda+A_1)^{-1}B_1[u_{\omega}] u, \varphi)_{L^2}|
\leq C\epsilon \|(-\Delta)^{\frac{k}{2}} u\|_{L^2} \|(-\Delta)^{1-\frac{k}{2}} (\lambda+A^*_1)^{-1}(-\Delta)^{\frac{k}{2}}\varphi\|_{L^2}. 
\end{align}
Indeed, let $\psi \in   L^{n,\infty}$ with $\nabla \psi \in L^{\frac{n}{2},\infty}$ satisfy $\|\psi\|_{L^{n,\infty}}+\|\nabla \psi\|_{L^{n/2,\infty}}\leq \epsilon$. 
It directly follows from Lemma \ref{fgh-est} {\rm (i)} and the property of the low frequency part, Lemma \ref{lemP_1} that for $3\leq n$ 
\begin{align}\label{est-perturbation-terms2-base}
|( \del_x  \psi   u, P_1(\lambda+A^*_1)^{-1}(-\Delta)^{\frac{k}{2}}\varphi)_{L^2} | 
\leq C\epsilon \|(-\Delta)^{\frac{k}{2}}u\|_{L^2}
\|(-\Delta)^{1-\frac{k}{2}}(\lambda+A^*_1)^{-1}(-\Delta)^{\frac{k}{2}}\varphi\|_{L^2}. 
\end{align}

When $1<k$,  {\rm i.e, } $\frac{1}{2} \leq \frac{k}{2} <  \max\{1, \frac{n}{4}\}$, Lemma \ref{fgh-est} {\rm (ii)} yields that 
$$
|( \psi \del_x u,  P_1 (\lambda+A^*_1)^{-1}(-\Delta)^{\frac{k}{2}}\varphi)_{L^2}|
\leq C\epsilon \|(-\Delta)^{\frac{k}{2}-\frac{1}{2}}\nabla  u\|_{L^2} \|(-\Delta)^{1-\frac{k}{2}} (\lambda+A^*_1)^{-1}(-\Delta)^{\frac{k}{2}}\varphi\|_{L^2}. 
$$

We note that $1-\frac{k}{2} < \max\{\frac{n}{4},1\}$. 
When $k<1$, {\rm i.e, } $0< \frac{k}{2}<\frac{1}{2}$, an integration by parts derive that 
\begin{align*}
( \psi \del_x u,  P_1 (\lambda+A^*_1)^{-1}(-\Delta)^{\frac{k}{2}}\varphi)_{L^2} &= -(\del_x \psi u,  P_1(\lambda+A^*_1)^{-1}(-\Delta)^{\frac{k}{2}}\varphi)_{L^2}\\
&\quad - (\psi u,  \del_x P_1 (\lambda+A^*_1)^{-1}(-\Delta)^{\frac{k}{2}}\varphi)_{L^2} \\
&= I_1+I_2. 
\end{align*}
We note that $(-\Delta)^{s}$ and $P_1$ are commute for $0\leq s \leq 1$, and $P_1=P_1^*$. 
$I_1$ is estimated by \eqref{est-perturbation-terms2-base}. Concerning $I_2$,  Lemma \ref{fgh-est} {\rm (ii)} also verifies that  
\begin{align}
\begin{aligned}\label{est-perturbation-terms2}
I_2  
&\leq C \epsilon \|(-\Delta)^{\frac{k}{2}} u\|_{L^2} \|(-\Delta)^{1-\frac{k}{2}-\frac{1}{2}} \nabla (\lambda+A^*_1)^{-1}(-\Delta)^{\frac{k}{2}}\varphi\|_{L^2}\\
&\leq C \epsilon \|(-\Delta)^{\frac{k}{2}} u\|_{L^2}\|(-\Delta)^{1-\frac{k}{2}}  (\lambda+A^*_1)^{-1}(-\Delta)^{\frac{k}{2}}\varphi\|_{L^2}. 
\end{aligned}
\end{align}
Substituting 
$\psi= \phi_\omega, w_\omega $ in $B_1[u_\omega]$ yields  \eqref{perturbation-est}. We can estimate $h^{(1)}(u_\omega)\phi_\omega \nabla \phi$ in $B_1[u_\omega]$ by the same argument because $\del_x(h^{(1)}(u_\omega))=\del_x h^{(1)}(u_\omega) \del_x \phi_\omega$ when we apply the integration by parts.  
We also can estimate $\phi_\omega \Delta w$ in $B_1[u_\omega]$ because 
$\|(-\Delta)^{k/2} \Delta w\|_{L^2_{(1)}} \leq C \|(-\Delta)^{k/2} \nabla  w\|_{L^2_{(1)}}$ by the low frequency property, Lemma \ref{lemP_1}. 

We see from Corollary \ref{resolvent2-cor} that 
$$
|((-\Delta)^{\frac{k}{2}} (\lambda+A_1)^{-1}B_1[u_{\omega}] u, \varphi)_{L^2}|
\leq C\epsilon \|(-\Delta)^{\frac{k}{2}} u\|_{L^2}. 
$$
This together with the duality argument implies (i) for $B_1[u_\omega]$. The estimate of $B_1^*[u_{\omega}]$ is analogous by the definition of $B_1[u_\omega]$. We note that 
since $L^2_{(1)}$ is the Hilbert space and $B_1[u_\omega]: L^2_{(1)}\rightarrow L^2_{(1)}$ is a bounded linear operator, $B_1^*[u_\omega]: L^2_{(1)}\rightarrow L^2_{(1)}$, as well as $ (\lambda+A^*_1)^{-1}B_1^*[u_{\omega}]:L^2_{(1)}\rightarrow L^2_{(1)}$.  
A similar argument to the proof of (i)  
derive (ii).    
\end{proof}

\vspace{2ex}

We have the following resolvent estimates of $(\lambda+L_1)^{-1}$. 

\vspace{2ex}

\begin{prop}\label{resolvent-est-L1}
Let $\lambda \in \mathbb{C}\setminus  \sigma (-{A}_1)$ as in Proposition {\rm \ref{resolvent2}} and $u_1=\trans(\phi_1,w_1)$ be a solution to \eqref{cns-linear-low-per2} for $F_1 \in L^2_{(1)}$.   Then it holds that 
\begin{align}
\begin{aligned}\label{resolvent-est-L1-1}
&\|(-\Delta)^{\frac{k}{2}} (\lambda+ L_1)^{-1} F_1\|_{L^2} \leq 
C|\lambda|^{\frac{k}{2}-1}\|F_1\|_{L^2},  \\ 
&\|(-\Delta)^{\frac{k}{2}} (\lambda+ L_1^*)^{-1} F_1\|_{L^2} \leq 
C|\lambda|^{\frac{k}{2}-1}\|F_1\|_{L^2},
\end{aligned}
\end{align}
for $\max\{0, 1-\dfrac{n}{4}\}< \dfrac{k}{2}< \min\{\dfrac{n}{4},1\}$ and 
\begin{align}
\begin{aligned}\label{resolvent-est-L1-2}
&\| (\lambda+ L_1)^{-1} F_1\|_{L^2} \leq 
C|\lambda|^{-1}\|F_1\|_{L^2},  \\
&\|(\lambda+ L_1^*)^{-1} F_1\|_{L^2} \leq 
C|\lambda|^{-1}\|F_1\|_{L^2}. 
\end{aligned}
\end{align}

\end{prop}

\vspace{2ex}

\eqref{resolvent-est-L1-1} directly follows from Proposition \ref{resolvent2}, Corollary \ref{resolvent2-cor-2}, Lemma \ref{L1-B1} {\rm (i)} with \eqref{L1-form} and \eqref{L1-form-adjoint} respectively. By analogy,  \eqref{resolvent-est-L1-2} follows from Proposition \ref{resolvent2}, Corollary \ref{resolvent2-cor-2}, Lemma \ref{L1-B1} {\rm (ii)} with \eqref{L1-form} and \eqref{L1-form-adjoint} respectively as 
$$
\|(\lambda+ L_1)^{-1} F_1\|_{L^2} \leq C|\lambda|^{-1}\|F_1\|_{L^2}
+ C|\lambda|^{-\frac{1}{2}}\sum_{j=1}^\infty \epsilon^j \|\nabla (\lambda+A_1)^{-1}F_1\|_{L^2}  \leq C |\lambda|^{-1}\|F_1\|_{L^2}. 
$$

\vspace{2ex}

By Corollary \ref{resolvent2-cor-2-modify} and \ref{resolvent2-cor-3-modify} and a similar argument to Proportion \ref{resolvent-est-L1}, we have 

\vspace{2ex}

\begin{cor}\label{resolvent-est-L1-modify}
{\rm (i)} Let $\lambda \in \mathcal R_{1,\pm}$ as in Proposition {\rm \ref{resolvent2}} and $u_1=\trans(\phi_1,w_1)$ be a solution to \eqref{cns-linear-low-per2} for $F_1 \in L^2_{(1)}$.   Then it holds that for any small $\epsilon>0$ independent of $\lambda$ and $c_0$ 
\begin{align}
\begin{aligned}\label{resolvent-est-L1-1-modify}
&\|(-\Delta)^{\frac{k}{2}} (\lambda+ L_1)^{-1} F_1\|_{L^2} 
 \leq 
C|\lambda|^{\frac{k}{2}-\frac{1}{2}-\epsilon}\|F_1\|_{L^2}  
\end{aligned}
\end{align}
for $\max\{0, 1-\dfrac{n}{4}\}< \dfrac{k}{2}< \dfrac{1}{2}$ and 
\begin{align}
\begin{aligned}\label{resolvent-est-L1-2-modify}
&\|(-\Delta)^{\frac{1}{2}} (\lambda+ L_1)^{-1}  F_1\|_{L^2}\leq 
C\|F_1\|_{L^2}. 
\end{aligned}
\end{align}

\end{cor}

\vspace{2ex}

We state the low frequency property of the resolvent $(\lambda+ L_1)^{-1}$. 

\vspace{2ex}

\begin{lem}\label{low-frequency-property-L1}
For $F_1 \in L^p_{(1)}$ with $1 < p \leq 2$ it holds that 
$\supp \mathcal{F}\{(\lambda+ L_1)^{-1} F_1\} \subset 
Q_{r'_\infty}$. 
\end{lem}

\vspace{2ex}

\begin{proof}

Let $u_1 = (\lambda+ L_1)^{-1} F_1$. 
We note that a similar argument to the proof of Lemma \ref{lemP_1} yields $F_1 \in L^2_{(1)}$ and thus $u_1$ is bounded in $L^2$ by \eqref{low-frequency-property-L1}. 
Let $R_{\lambda,k}:= \sum_{j=0}^{k}\{(\lambda+A_1)^{-1}B_1[u_{\omega}]\}^{j}(\lambda+A_1)^{-1}$.  Since $\mathcal{F}(\lambda+A_1)^{-1}=(\lambda+\hat{A}_1)^{-1}$ and $B_1[u_{\omega}]$ has the operator $P_1$, we see that $R_{\lambda,k}F_1 \in L^2_{(1)}$ and 
$$
R_{\lambda,k}F_1 \rightarrow (\lambda+ L_1)^{-1} F_1 \ \ \mbox{in}  \ \ L^2
$$
as $k$ goes to  $\infty$. Hence 
$$
\mathcal{F} R_{\lambda,k}F_1 \rightarrow \mathcal{F}[(\lambda+ L_1)^{-1} F_1] \ \ \mbox{in}  \ \ L^2
$$
as $k$ goes to  $\infty$. Let any bounded continuous function $\hat{\chi}_{0}$ satisfying $\hat{\chi}_0 =0$ on $Q_{r'_\infty}$. Since $\hat{\chi}_{0}$ is bounded on $L^2$, we obtain that 
$$
\hat{\chi}_{0}\mathcal{F} R_{\lambda,k}F_1 \rightarrow \hat{\chi}_{0}\mathcal{F}[(\lambda+ L_1)^{-1} F_1] \ \ \mbox{in}  \ \ L^2
$$
as $k$ goes to  $\infty$. Therefore, 
$$
\|\hat{\chi}_{0}\mathcal{F}[(\lambda+ L_1)^{-1} F_1]\|_{L^2} 
= \lim_{k \rightarrow \infty} \|\hat{\chi}_{0}\mathcal{F} R_{\lambda,k}F_1\|=0, 
$$
and thus $\hat{\chi}_{0}\mathcal{F}[(\lambda+ L_1)^{-1} F_1]=0$. 
This implies that $\supp \hat{u}_1 \subset 
Q_{r'_\infty}.$ 
\end{proof}

\vspace{2ex}

Let 
$\Sigma_{\epsilon}=\{z \in \mathbb{C}\setminus\{0\}; |\mbox{arg}z|\leq \pi-\epsilon\}$. Let $\Gamma$ be a contour in $\Sigma_{\epsilon}$ satisfying that 
for $t>0$, $\Gamma=\Gamma_+ \cap \Gamma_0 \cap \Gamma_-$ with   
\begin{align}
\begin{aligned}\label{integral-curve}
&\Gamma_\pm=\{\lambda \in \mathbb{C}; \lambda =  -r^2 \pm\Big(r + \frac{1}{t}\Big)i, \, r> 0\}, \\
&\Gamma_0 = \{\lambda \in \mathbb{C};\lambda=\frac{1}{t}e^{-i\theta}, -\frac{\pi}{2} \leq \theta \leq \frac{\pi}{2}\}.   
\end{aligned}
\end{align}
Then the Cauchy formula yields  the semigroup generated by $L_1$, namely,    
\begin{align}\label{cauchy-formula}
S_1(t)=
e^{-t L_1}=
\frac{1}{2\pi i}\displaystyle\int_{\Gamma}(\lambda+L_1)^{-1}e^{\lambda t}d\lambda. 
\end{align}

We next consider time decay estimates of the semigroup of the perturbed system in the low frequency part. 

\vspace{2ex}

\vspace{2ex}
\begin{prop}\label{analytic semigroup2}
For $u_1=\trans(\phi_1,w_1)=S_1(t)u_{01}=e^{-t L_1}u_{01}$ with $u_{01} \in L^2_{(1)}$, we have the decay estimate  
$$
\| u_1\|_{L^q} \leq C(1+t)^{-\delta+\epsilon}\|u_{01}\|_{L^2} , 
$$
where $\delta=\dfrac{n}{2}\Big(\dfrac{1}{2}-\dfrac{1}{q}\Big)$ for $2 < q$ with 
$\max\{0, 1-n/4\} < \delta < \min\{n/4,1\}$ and any small $\epsilon>0$ independent of $t$ and $\delta$ or $q=2$ and $\delta=\epsilon=0$.\BLACK   

In addition,  it holds that 
$$
\| \nabla u_1\|_{L^q} \leq C(1+t)^{-\delta+\epsilon-\frac{1}{2}}\|u_{01}\|_{L^2} .
$$
\end{prop}

\vspace{1ex}

\begin{proof}
Set 
$\mathcal L:=\mathcal L(L^q;L^2_{(1)})$.  We see from the Sobolev embedding with the fractional Laplacian that 
$$
\|(\lambda+L_1)^{-1}\|_{\mathcal L(L^2_{(1)};L^q)} 
\leq C\|(-\Delta)^{\delta}(\lambda+L_1)^{-1}\|_{\mathcal L(L^2_{(1)};L^2_{(1)})}. 
$$ 

We first consider the case $\max\{0,1-n/4\} <\delta<1/2$.  
Applying Proposition \ref{resolvent-est-L1} and Corollary \ref{resolvent-est-L1-modify} to the Cauchy formula \eqref{cauchy-formula} verifies that 

\begin{align}\label{e-tL}
    \|e^{-t L_1}\|_{\mathcal L} \leq C\int_{1/t}^\infty |e^{-\lambda t}|
    |\lambda|^{\delta-1/2-\epsilon}
    |\dlambda|\, + \, C\int_{\partial B_{1/t}} |\lambda|^{\delta-1} |\dlambda|. 
\end{align}
If $t<1$, then $|\lambda|>1$ in \eqref{e-tL}, and we get that
\begin{align}\label{e-tL-2}
\|e^{-t L_1}\|_{\mathcal L} \leq C\int_{1/t}^\infty e^{-c r^2 t} r^{\delta-1} \sqrt{1+4r^2}\dr + C\frac1t t^{1-\delta} \leq C(t^{-\delta}e^{-1/t}+ t^{-\delta}) \leq Ct^{-\delta}.  
  \end{align} 
\BLACK
Moreover, if $t>1$, we split the first integral in \eqref{e-tL-2} into two parts and obtain that for $0< \epsilon< \delta<1$  
\begin{align}\label{e-tL-3}
\begin{aligned}
    \|e^{-t L_1}\|_{\mathcal L} 
    & \leq C\int_{1/t}^1 e^{-c r^2 t} r^{2(\delta-1/2-\epsilon)} \sqrt{1+4r^2} \dr + C\int_1^\infty e^{-c r^2 t} r^{\delta-1} \sqrt{1+4r^2}\dr + \frac1t t^{1-\delta} \\
    & \leq C (t^{-\delta+\epsilon} + e^{-t} + t^{-\delta})\\
     & \leq C t^{-\delta+\epsilon}.
     \end{aligned}
\end{align} 
\BLACK 
Summarizing \eqref{e-tL}-\eqref{e-tL-3}, 
 we get that 
$$
\| u_1\|_{L^q} \leq Ct^{-\delta+\epsilon}\|u_{01}\|_{L^2}. 
$$
The case $1/2< \delta$ is derived by 
the formula \eqref{L1-form} and the semigroup property. 

Since $L_1$ is a bounded linear operator on $L^2_{(1)}$ due to Lemma \ref{lemP_1} and satisfies 
$$
\|L_1 u_1\|_{L^2} \leq C\|u_1\|_{L^2}, 
$$
$e^{-t L_1}$ is a uniformly continuous semigroup on $L^2_{(1)}$. This implies the case $q=2$ and $\delta=\epsilon=0$. 

The boundedness at $t=0$ is obtained by the property of the low frequency part, {\em i.e., } 
\begin{align*}
    \| e^{-t L_1} u_{01} \|_{L^q} 
    & \leq C\|(-\Delta)^\delta  e^{-t L} u_{01} \|_{L^2_{(1)}}\\[1ex]
    & \leq C\|e^{-t L} u_{01}  \|_{L^2_{(1)}} \leq C\|u_{01}\|_{L^2}. 
\end{align*}
Finally, the gradient estimate is obtained by analogy with the semigroup property.  
\end{proof}

\vspace{2ex}

In the proof of Proposition \ref{analytic semigroup2} we also have 

\vspace{2ex}

\begin{cor}\label{analytic semigroup2-cor}
For $u_1=\trans(\phi_1,w_1)=S_1(t)u_{01}=e^{-t L_1}u_{01}$ with $u_{01} \in L^2_{(1)}$, we get the decay estimate  
$$
\|(-\Delta)^\delta u_1\|_{L^2} \leq C(1+t)^{-\delta+\epsilon}\|u_{01}\|_{L^2} , 
$$
where 
$\max\{0, 1-n/4\} < \delta < \min\{n/4,1\}$ and any small $\epsilon>0$ independent of $t$ and $\delta$.  
\end{cor}

\vspace{1ex}

By virtue of the resolvent estimates of the adjoint operator $L_1^*$ in Proposition \ref{resolvent-est-L1}, the same argument as that in Proposition \ref{analytic semigroup2} implies 

\vspace{2ex}

\begin{prop}\label{analytic semigroup2-adjoint}
{\rm (i)} For $u_1=\trans(\phi_1,w_1)=S_1^*(t)u_{01}=e^{-t L_1^*}u_{01}$ with $u_{01} \in L^2_{(1)}$, we get the decay estimates  
$$
\| u_1\|_{L^q} \leq C(1+t)^{-\delta+\epsilon}\|u_{01}\|_{L^2} , 
$$
where $\delta=\frac{n}{2}\Big(\frac{1}{2}-\frac{1}{q}\Big)$ for $2 < q$ with 
$\max\{0, 1-n/4\} < \delta < \min\{n/4,1\}$  and any small $\epsilon>0$ independent of $t$ and $\delta$ or $q=2$ and $\delta=0$.  

Moreover,  it holds that 
$$
\| \nabla u_1\|_{L^q} \leq C(1+t)^{-\delta+\epsilon-\frac{1}{2}}\|u_{01}\|_{L^2} .
$$
\vspace{1ex}
{\rm (ii)}
We have the decay estimate  
$$
\|(-\Delta)^\delta u_1\|_{L^2} \leq C(1+t)^{-\delta+\epsilon}\|u_{01}\|_{L^2} , 
$$
where 
$\max\{0, 1-n/4\} < \delta < \min\{n/4,1\}$ and  any small $\epsilon>0$ independent of $t$ and $\delta$.   
\end{prop}

\vspace{2ex}

The duality argument together with Proposition \ref{analytic semigroup2-adjoint}  verifies the following $L^2$-$L^p$ estimates of the semigroup. 

\begin{prop}\label{analytic semigroup3}
For $u_1=\trans(\phi_1,w_1)=S_1 (t)u_{01}=e^{-t L_1}u_{01}$ with $u_{01} \in L^{p}_{(1)}$ and $1<p<2$, we have the  decay estimate 
$$
\| u_1\|_{L^2} \leq C(1+t)^{-\delta_2+\epsilon}\|u_{01}\|_{L^p},
$$ 
where $\delta_2=\frac{n}{2}\Big(\frac{1}{p}-\frac{1}{2}\Big)$ for $1< p< 2$ with 
$\max\{0, 1-n/4\} < \delta < \min\{n/4,1\}$  and any small $\epsilon>0$ independent of $t$ and $\delta_2$ or $q=2$ and $\delta_2=0$.  Furthermore, we have  that 
$$
\| \nabla u_1\|_{L^2} \leq C(1+t)^{-\delta_2+\epsilon-\frac{1}{2}}\|u_{01}\|_{L^p}. 
$$
\end{prop}

\vspace{2ex}

\begin{proof}
We see from Lemma \ref{low-frequency-property-L1} that $\supp \hat{u}_1 \subset Q_{r'_\infty}$ and thus $u_1 \in L^2_{(1)}$. 
Hence for $\varphi \in L^2_{(1)}$, 
$$
(u_1, \varphi) = (u_{01}, e^{-t L_1^*} \varphi). 
$$
Proposition \ref{analytic semigroup2-adjoint} can be applied to  $e^{-t L_1^*} \varphi$ and we get that 
$$
\|e^{-t L_1^*}\varphi\|_{L^{p'}} \leq  C(1+t)^{-\delta_2+\epsilon}\|\varphi\|_{L^2}. 
$$
Hence the duality argument yield that
$$
\| u_1\|_{L^2} \leq C(1+t)^{-\delta_2+\epsilon}\|u_{01}\|_{L^p}. 
$$ 
The estimate $\| \nabla u_1\|_{L^2}$ is obtained by the semigroup property. 
\end{proof}
\vspace{2ex}

By the duality argument with Proposition \ref{analytic semigroup2-adjoint}, we also have the following decay estimate which is useful for  nonlinear terms having divergence form. 

 \vspace{2ex}
 
 \begin{cor}\label{analytic semigroup3-divergence-form}
Let $u_1=\trans(\phi_1,w_1)=S_1 (t) \del_x u_{01}=e^{-t L_1}u_{01}$, where 
$u_{01} \in L^p_{(1)}$ for $1<p<2$ close to $2$ satisfying that $ \delta_3+\frac{1}{2}< \max\{\frac{n}{4},1\}$ with $\delta_3=\frac{n}{2}\Big(\frac{1}{p}-\frac{1}{2}\Big)$ and any small $\epsilon>0$ independent of $t$ and $\delta_3$.  Then  we have that 
$$
\| u_1\|_{L^2} \leq C(1+t)^{-\delta_3+\epsilon-\frac{1}{2}}\|u_{01}\|_{L^p}. 
$$
\end{cor}

\vspace{2ex}

The duality argument shows that for $\varphi \in L^2_{(1)}$, 
$$
(e^{-t L_1}\del_x u_{01}, \varphi)= -( u_{01}, \del_x e^{-t L^*_1}\varphi), 
$$
where 
\begin{align}\label{cauchy-formula-adjoint}
e^{-t L^*_1}=
\frac{1}{2\pi i}\displaystyle\int_{\Gamma}(\lambda+L^*_1)^{-1}e^{\lambda t}d\lambda. 
\end{align}
Substituting Proposition \ref{analytic semigroup2-adjoint} derives Corollary \ref{analytic semigroup3-divergence-form}.

\end{section}

\begin{section}{Estimates for the high frequency part}

On the high frequency part, we can apply the $L^2$ energy estimate.  \eqref{cns-linear-high-per} is concretely written  as follows. 

\vspace{2ex}

\begin{eqnarray}
\left\{\hspace*{11pt}
\begin{array}{lll}
 \partial_{t}\phi_\infty +\gamma \div  w_\infty+ \gamma P_\infty\tilde{w}\cdot \nabla \phi_{\infty} + \gamma P_\infty ( b_{11}[u_\omega]u_\infty+P_\infty b_{12}[u_\omega]u_\infty)
\hspace*{-4pt}  &=  f_{\infty 1},\\[1ex]
 \partial_{t}w_\infty- \alpha \Delta w_\infty- \beta \nabla\div w_\infty+\gamma  \nabla \phi_\infty+ \gamma P_\infty b_{21}[u_\omega]u_\infty\\
 \qquad + \gamma P_\infty b_{22}[u_\omega]u_\infty
+ \gamma P_\infty b_{23}[u_\omega]u_\infty
\hspace*{-4pt}  &= f_{\infty 2},\\[1ex]
 \phi_\infty|_{t=0}=\phi_{0\infty}, \ \ w_\infty|_{t=0}=w_{0\infty},& 
\label{cns-linear-high-energy}
\end{array}
\right.
\end{eqnarray}

\vspace{2ex}

Concerning the local existence of solutions, 
we first  obtain  local existence of solutions to the system with  $P_\infty (w_\omega \cdot \nabla \phi_\infty)$ in the right hand side by the same argument as that in \cite[Proposition 6.4]{Kagei-Tsuda}. 
Then  the direct computations based on Lemmata \ref{lem2.1.} and \ref{composition-est} yield 
$$
\| b_{11}(u_\omega, u_\infty)\|_{L^2([0,T']; H^s)}+ \|b_{12}(u_\omega, \nabla u_\infty)-w_\omega \cdot \nabla \phi_\infty\|_{L^2([0,T']; H^s)} \leq C\epsilon\|(\phi_\infty, w_\infty)\|_{L^2([0,T']; H^s \times H^{s+1})}
$$
and 
\begin{align*}
&\|b_{21}(u_\omega, u_\infty)\|_{L^2([0,T']; H^{s-1})} + \|b_{22}(u_\omega, \nabla u_\infty)\|_{L^2([0,T']; H^{s-1})} \\
& \quad + \|b_{23}(u_\omega, \nabla^2 u_\infty)\|_{L^2([0,T']; H^{s-1})}\\
&\qquad \leq  C\epsilon\|(\phi_\infty, w_\infty)\|_{L^2([0,T'];H^s \times  H^{s+1})}. 
\end{align*}
Hence \cite[Proposition 6.4]{Kagei-Tsuda} with the standard iteration argument verify the local existence of solutions to \eqref{cns-linear-high-energy} satisfying that 
\begin{eqnarray*}
\phi_{\infty}\in C([0,T'];H^{s}_{(\infty)}),\ 
w_{\infty}\in C([0,T'];H^{s}_{(\infty)})\cap L^{2}(0,T';H^{s+1}_{(\infty)})
\end{eqnarray*}
for all $T'>0$. 

We state an a priori estimate for the high frequency part as follows. 

\vspace{2ex}

\begin{prop}\label{energyest}
Let $s \geq 0$ be an integer with $s\geq [n/2]+2$. Suppose that 
\begin{eqnarray*}
u_{0\infty}=\trans(\phi_{0\infty},v_{0\infty})\in H^{s}_{(\infty)}\\
F_\infty =\trans(f_{\infty 1}, f_{\infty 2})\in L^2(0,T';  H^{s-1}_{(\infty)})
\end{eqnarray*}
for all $T'>0$. We suppose that $u_{\infty}=\trans(\phi_{\infty},w_{\infty})$ satisfies \eqref{cns-linear-high-per}  
and 
\begin{eqnarray*}
\phi_{\infty}\in C([0,T'];H^{s}_{(\infty)}),\ 
w_{\infty}\in C([0,T'];H^{s}_{(\infty)})\cap L^{2}(0,T';H^{s+1}_{(\infty)})
\end{eqnarray*}
for all $T'>0$.  
Then there eixst a positive constant $a_0$ and  an energy functional ${\cal E}[u_{\infty}]$ 
such that if 
$$
\|{w}\|_{C([0,T]; H^s)} \leq a_0, 
$$
then 
\begin{eqnarray}\label{energy}
\frac{d}{dt}{\cal E}[u_{\infty}](t)
+{d}(\|\nabla \phi_{\infty}(t)\|_{H^{s}}^{2}+\|\nabla w_{\infty}(t)\|_{H^{s}}^{2})
 \leq
C\|F_{\infty}(t)\|_{H^s\times H^{s-1}}^2
\end{eqnarray} 
on $(0,T')$ for all $T'>0$. 
Here $d$ is a positive constant,   
$C$ is a positive constant independent of  $T'$,   
${\cal E}[u_{\infty}]$ is equivalent to $\|u_{\infty}\|_{H^s}^2$, i.e, 
$$
C^{-1}\|u_{\infty}\|_{ H^s}^2
\leq {\cal E}[u_{\infty}]
\leq C\|u_{\infty}\|_{H^s}^2,  
$$
and ${\cal E}[u_{\infty}](t)$ is absolutely continuous in $t\in [0,T']$ for all $T'>0$.
\end{prop}

\vspace{2ex}

To prove \eqref{energy} 
for $0 \leq |k| \leq s$, 
we consider the $L^2$ inner products 
$\del_x^k \eqref{cns-linear-high-energy}_1 \times \del_x^k \phi_\infty$ and $\del_x^k \eqref{cns-linear-high-energy}_2 \times \del_x^k w_\infty$ respectively. Then we obtain that 
\begin{align}
\begin{aligned}\label{energy-1}
&\dfrac{d}{dt}\|\del_x^k \phi_\infty\|_{L^2}^2 
+  \gamma (\del_x^k\div w_\infty, \del_x^k \phi_\infty) 
=\sum_{j=1}^3 J_{1j,k}, \\
&\dfrac{d}{dt}\|\del_x^k w_\infty\|_{L^2}^2 
+ \alpha\|\del_x^k \nabla w_\infty\|_{L^2}^2
+ \beta \|\del_x^k \div w_\infty\|_{L^2}^2
+\gamma ( \del_x^k\nabla \phi_\infty, \del_x^k w_\infty) =\sum_{j=1}^2 J_{2j,k}, 
\end{aligned}
\end{align}
where 
\begin{align*}
J_{11,k}&=-\gamma (P_\infty \del_x^k(\tilde{w}\cdot \nabla \phi_\infty), \del_x^k\phi_\infty ), \\
J_{12,k}&= -\gamma(P_\infty (\del_x^k( b_{11}[u_\omega]u_\infty+ b_{12}[u_\omega]u_\infty), \del_x^k\phi_\infty), \\
J_{13,k}&=(\del_x^k f_{\infty 1}, \del_x^k \phi_\infty), \\
J_{21,k}&=-\gamma(P_\infty (\del_x^k( b_{21}[u_\omega]u_\infty+ b_{22}[u_\omega]u_\infty+ b_{23}[u_\omega]u_\infty), \del_x^k w_\infty),  \\
J_{22,k}&=(\del_x^k f_{\infty 2}, \del_x^k w_\infty). 
\end{align*}
In addition, for $0 \leq |k| \leq s-1$ by taking the the $L^2$ inner products 
$\del_x^k \eqref{cns-linear-high-energy}_2 \times \del_x^k \nabla \phi_\infty$ we have that 
\begin{align}
\begin{aligned}\label{energy-2}
\dfrac{d}{dt}(\del_x^k w_\infty, \del_x^k \nabla \phi_\infty) 
+\gamma \|\del_x^k \nabla \phi_\infty\|_{L^2}^2=\sum_{j=1}^3 J_{3j,k},   
\end{aligned}
\end{align}
where 
\begin{align*}
J_{31,k}&=-\alpha(\del_x^k \Delta w_\infty, \del_x^k \nabla \phi_\infty )  
- \beta (\del_x^k \div w_\infty, \del_x^k \nabla \phi_\infty ), \\
J_{32,k} &= - \gamma(P_\infty (\del_x^k( b_{21}[u_\omega]u_\infty+ b_{22}[u_\omega]u_\infty+ b_{33}[u_\omega]u_\infty), \del_x^k \nabla \phi_\infty) \\
J_{33,k}&=(\del_x^k f_{\infty 2}, \del_x^k \nabla \phi_\infty)
.
\end{align*}
We prepare estimates of $J_{ij,k}$ as follows. 

\vspace{2ex}

\begin{lem}\label{estimate-energy-reminder1}
It holds that 
\begin{align*}
    |J_{11,k}| &\leq C\|\div \tilde{w}\|_{H^{s-1}}\| \phi_\infty\|_{H^{s-1}}^2,  \ 0 \leq |k| \leq s, \\
    |J_{13,k}| & \leq \epsilon_1 \|\phi_\infty\|_{H^s}^2, + C\|f_{\infty1}\|_{H^2}^2, \ 0 \leq |k| \leq s, \\
    |J_{22,k}| & \leq \epsilon_1 \|\nabla w_\infty\|_{H^s}^2+ C\|f_{\infty2}\|_{H^{s-1}}^2,  \ 0 \leq |k| \leq s, \\
    |J_{33,k}| & \leq \epsilon_1\|\nabla \phi_\infty\|_{H^{s-1}}^2 +C\|f_{\infty2}\|_{H^{s-1}}^2, \ 0 \leq |k| \leq s-1. 
\end{align*}

\end{lem}

\vspace{2ex}

Lemma \ref{estimate-energy-reminder1} is derived directly by 
the Cauchy–Schwarz inequality,  integration by parts, and  Lemmata \ref{lem2.1.} and \ref{composition-est}. We omit details, see also \cite{Matsumura-Nishida2}. We note that since $[n/2]+ 2 \leq s $,  
the Sobolev embedding, Lemma \ref{lem2.1.} derives the norm $\|\div \tilde{w}\|_{H^{s-1}}$.  

\vspace{2ex}

On the perturbation terms including $u_\omega$, we have by the Hardy inequality 

\vspace{2ex}

\begin{lem}\label{estimate-energy-reminder2} 
Under \eqref{u-omega-esimate} there holds that for $0 \leq |k| \leq s $ on 
$J_{12,k}$ and $J_{21,k}$, and for $0 \leq |k| \leq s-1 $ on $J_{32,k}$  
\begin{align}
|J_{12,k}| + |J_{21,k}|+ |J_{32,k}| 
\leq C\epsilon(\|\nabla w_\infty\|_{H^s}(\|\nabla w_\infty\|_{H^s}+\|\nabla \phi_\infty\|_{H^{s-1}}) + \|\nabla \phi_\infty\|_{H^{s-1}}^2). 
\end{align}
\end{lem}

\vspace{2ex}

\begin{proof}
We estimate $|J_{22,k}|$. On the case $k=0$,  integration by parts,  Lemmata \ref{lem2.1.} and 
the property of the high frequency part, Lemma \ref{lemPinfty} 
derive that 
\begin{align*}
|(P_\infty(\phi_\omega\div w_{\infty}, w_\infty))_{L^2}| &\leq C\|\phi_\omega\|_{L^\infty}\|\nabla w_{\infty}\|_{L^2}^2 \leq C\|\nabla \phi_\omega\|_{H^{s-1}}\|\nabla w_{\infty}\|_{L^2}^2, \\
|(P_\infty(w_\infty\cdot \nabla \phi_\omega, w_\infty))_{L^2}| &\leq C\|\nabla \phi_\omega\|_{L^\infty}\|\nabla w_{\infty}\|_{L^2}^2 \leq C\|\nabla \phi_\omega\|_{H^{s-1}}\|\nabla w_{\infty}\|_{L^2}^2, \\
|(P_\infty(\phi_\infty\div w_{\omega}, w_\infty))_{L^2} |&\leq C\|\nabla w_\omega\|_{L^\infty}\|\nabla w_{\infty}\|_{L^2}\|\nabla \phi_{\infty}\|_{L^2}\leq C\|\nabla w_\omega\|_{H^{s-1}}\|\nabla w_{\infty}\|_{L^2}\|\nabla \phi_{\infty}\|_{L^2}, \\
|(P_\infty(w_\omega\cdot \nabla \phi_\infty, w_\infty))_{L^2} |&\leq C\|w_\omega\|_{L^\infty}\|\nabla w_{\infty}\|_{L^2}\|\nabla \phi_{\infty}\|_{L^2}, \\
|(P_\infty(w_\omega\cdot \nabla  w_{\infty}, w_\infty))_{L^2}| &\leq C\|w_\omega\|_{L^\infty}\|\nabla w_{\infty}\|_{L^2}^2, \\
|(P_\infty(w_\infty\cdot \nabla w_\omega, w_\infty))_{L^2}| &\leq C\|\nabla w_\omega\|_{L^\infty}\|\nabla w_{\infty}\|_{L^2}^2, \\
|(P_\infty(\phi_\omega \nabla  \phi_{\infty}, w_\infty))_{L^2}| &\leq C\|\phi_\omega\|_{L^\infty}
\|\nabla \phi_\infty\|_{L^2}\|\nabla w_\infty\|_{L^2}, \\
|(P_\infty(\phi_\infty \nabla  \phi_{\omega}, w_\infty))_{L^2}| &\leq C\|\nabla \phi_\omega\|_{L^\infty}
\|\nabla \phi_\infty\|_{L^2}\|\nabla w_\infty\|_{L^2}, \\
|(P_\infty(\phi_\omega\Delta  w_{\infty}), w_\infty)_{L^2}| 
&=|- (\phi_\omega, |\nabla w_\infty|^2) - (\nabla \phi_\omega \cdot w_\infty, \nabla w_\infty)|\\
& \leq C(\|\nabla \phi_\omega\|_{L^\infty}+\|\phi_\omega\|_{L^\infty})\|\nabla w_{\infty}\|_{L^2}^2, \\
|(P_\infty(\phi_\infty\Delta  w_\omega) , w_\infty)_{L^2}| 
& \leq C\|\Delta w_\omega\|_{L^\infty}\|\nabla w_{\infty}\|_{L^2}\|\nabla \phi_\infty\|_{L^2}\leq C\|\Delta w_\omega\|_{H^{s-1}}\|\nabla w_{\infty}\|_{L^2}\|\nabla \phi_\infty\|_{L^2}, 
\end{align*} 
Furthermore, Lemma \ref{composition-est} verifies that for $1 \leq |k| \leq s$ with $s\geq [n/2]+2$ 
\begin{align*}
|(P_\infty(\del_x^k (\phi_\omega\div w_{\infty}), \del_x^k w_\infty))|&\leq |( \del_x^k\div w_{\infty} \phi_\omega, \del_x^k w_\infty)| + |\del_x^k(\div w_{\infty} \phi_\omega - \del_x^k \div w_{\infty} \phi_\omega, \del_x^k w_\infty)|  \\
&\leq C\|\nabla w_\infty\|_{H^{s}}^2\|\nabla \phi_\omega\|_{H^{s-1}},\\
|(P_\infty(\del_x^k (w_\infty \cdot \nabla \phi_{\omega}), \del_x^k w_\infty))| &\leq |(  w_{\infty}\cdot \del_x^k\nabla  \phi_\omega, \del_x^k w_\infty)| + |\del_x^k(w_{\infty}\cdot \nabla  \phi_\omega -  w_{\infty}\cdot \del_x^k\nabla  \phi_\omega, \del_x^k w_\infty)|  \\
&\leq C\|\nabla w_\infty\|_{H^{s-1}}^2\|\nabla \phi_\omega\|_{H^{s-1}}
).
\end{align*}
Similarly we obtain that 
\begin{align*}
    |(P_\infty\del_x^k(\phi_\infty\div w_{\omega}), \del_x^kw_\infty)_{L^2}| &\leq C\|\nabla w_\omega\|_{H^{s}}\|\nabla w_{\infty}\|_{H^s}\|\nabla \phi_{\infty}\|_{H^{s-1}}, \\
|(P_\infty\del_x^k(w_\omega\cdot \nabla  w_{\infty}), \del_x^kw_\infty)_{L^2}| &\leq C\|\nabla w_\omega\|_{H^{s}}\|\nabla w_{\infty}\|_{H^s}^2, \\
|(P_\infty\del_x^k(w_\infty\cdot \nabla w_\omega)|, \del_x^kw_\infty)_{L^2} &\leq C\|\nabla w_\omega\|_{H^s}\|\nabla w_{\infty}\|_{H^s}^2, 
|(P_\infty\del_x^k(\phi_\infty \nabla \phi_\omega), \del_x^k w_\infty)_{L^2}|\\ &\leq C\|\nabla \phi_\omega\|_{H^s}\|\nabla \phi_{\infty}\|_{H^{s-1}}^2,  
\end{align*}
and an integration by parts derives that 
\begin{align*}
    |(P_\infty\del_x^k(w_\omega\cdot \nabla \phi_\infty), \del_x^kw_\infty)_{L^2}| &
    =|- (\del_x^k(w_\omega\cdot  \phi_\infty), \del_x^k \div w_\infty)_{L^2} - (\del_x^k(\nabla w_\omega\cdot \phi_\infty), \del_x^k w_\infty)_{L^2} |\\
    &\leq C\|\nabla w_\omega \|_{H^s}\|\nabla w_\infty\|_{H^s}\|\nabla \phi_{\infty}\|_{H^{s-1}}, \\
|(P_\infty\del_x^k(\phi_\omega\cdot \nabla \phi_\infty), \del_x^k w_\infty)_{L^2} |&
    =|- (\del_x^k(\phi _\omega\cdot  \phi_\infty), \del_x^k \div w_\infty)_{L^2} - (\del_x^k(\nabla \phi_\omega  \phi_\infty), \del_x^k w_\infty)_{L^2} |\\
    &\leq C\|\nabla \phi_\omega\|_{H^s}\|\nabla w_\infty\|_{H^s}\|\nabla \phi_{\infty}\|_{H^{s-1}}, \\
|(P_\infty\del_x^k(\phi_\omega\Delta  w_{\infty}), \del_x^kw_\infty)_{L^2}| 
&=|- (\del_x^k(\nabla \phi_\omega\cdot \nabla  w_\infty), \del_x^k  w_\infty)_{L^2} - (\del_x^k(\phi_\omega\nabla w_\infty), \del_x^k \nabla w_\infty)_{L^2} |\\
&\leq C\|\nabla w_\infty\|_{H^s}^2\|\nabla \phi_{\omega}\|_{H^{s}}, \\
|(P_\infty\del_x^k(\phi_\infty\Delta  w_\omega) , \del_x^k w_\infty)_{L^2}|
&= |- (\del_x^k(\nabla \phi_\infty \cdot \nabla  w_\omega), \del_x^k  w_\infty)_{L^2} - (\del_x^k(\phi_\infty\nabla w_\omega), \del_x^k \nabla w_\infty)_{L^2}|\\
&\leq  |(\del_x^k \phi_\infty \cdot \nabla^2  w_\omega), \del_x^k  w_\infty)_{L^2} 
+ (\del_x^k \phi_\infty \cdot  w_\omega, \del_x^k  \nabla w_\infty)_{L^2}|\\
&\quad 
+| (\del_x^k (\nabla \phi_\infty \cdot \nabla  w_\omega)| +| \del_x^k \nabla \phi_\infty \cdot \nabla  w_\omega, \del_x^k  \nabla w_\infty)_{L^2}|\\
&\quad 
+| (\del_x^k(\phi_\infty\nabla w_\omega), \del_x^k \nabla w_\infty)_{L^2}|\\
&\leq C\|\nabla w_\infty\|_{H^s}\|\nabla w_\omega\|_{H^s}\|\nabla \phi_{\infty}\|_{H^{s-1}}. 
\end{align*}
Since 
$$
\|h^{(1)}(\phi_\omega)\|_{L^\infty}
+ \|\nabla h^{(1)}(\phi_\omega)\|_{{H}^{s-1}} 
\leq C\|\nabla \phi_\omega\|_{H^{s-1}}, 
$$
we obtain the desired estimate of $J_{22,k}$. 
Estimates of $J_{12,k}$ and $J_{32,k}$ are analogous to  $J_{22,k}$. 
\end{proof} 

\vspace{2ex}
Let 
$$
\|\tilde{w}\|_{C([0,T; H^s])} \leq a. 
$$
Summarizing \eqref{energy-1} and \eqref{energy-2} on $k$, substituting 
Lemmata \ref{estimate-energy-reminder1} and \ref{estimate-energy-reminder2}, and taking $a$, $\epsilon_1$ and $\epsilon$ small,  
we see that 
\begin{eqnarray}\label{energy2}
\frac{d}{dt}{\cal E}[u_{\infty}](t)
+{d}(\| \phi_{\infty}(t)\|_{H^{s}}^{2}+\| w_{\infty}(t)\|_{H^{s+1}}^{2})
 \leq
C(\|f_{\infty 1}(t)\|_{H^{s}}^2+\|f_{\infty 2}(t)\|_{H^{s-1}}^2), 
\end{eqnarray} 
where we used  the Poincar\'e type inequality Lemma \ref{lemPinfty} so that 
$$
\epsilon\|u_\infty\|_{H^2} \leq C \epsilon\|\nabla u_\infty\|_{H^{s-1}} 
$$
and thus $C\epsilon\|\nabla u_\infty\|_{H^{s-1}}$ is absorbed to the left hand side. 

\end{section}

\begin{section}{Existence of global solutions with decay rates}\label{global sol}

On the nonlinear estimates, we apply the Duhamel formula further with the decay of the solution operator in the low frequency part, while the energy estimates are applied to the high frequency part. 

Suppose that $\|u_0\|_{L^p \cap H^s} \leq E_0$ for $1<p < 2$,  
$$
1< \frac{n}{2p},  \ \ 1< \frac{n}{2}\big(\frac{1}{p}-\frac{1}{2}\big)+\dfrac{1}{2}.  
$$
In addition, $\dfrac{1}{p}=\dfrac{1}{2}+\dfrac{1}{p_0}$ and $k=\dfrac{n}{2}\big(\dfrac{1}{2}-\dfrac{1}{p_0}\big)$ for $2<p_0$ close to $2$ with $2k\leq s$.  For $2<q$ we take $\delta_1=\min\{\dfrac{n}{2}\Big(\dfrac{1}{p}-\dfrac{1}{q}\Big)-\epsilon, 1+ \dfrac{n}{2}\Big(\dfrac{1}{p}-\dfrac{1}{2}\Big)-\epsilon\}$ for any small $\epsilon>0$.  
We can take a large $q$ satisfying that $\delta_1 >1$ and $\dfrac{n}{2q}< \dfrac{n}{4}-\dfrac{1}{2}$.  
We also set  $\delta_2 = \dfrac{n}{2}\Big(\dfrac{1}{p}-\dfrac{1}{2}\Big)-\epsilon$. 
We note that $\delta_1 > \delta_2 + \dfrac{1}{2}$, 
$1< \delta_2+ \dfrac{1}{2}$ and $\max\{0, 1-\dfrac{n}{4}\} < \delta_2 < \min\{\dfrac{n}{4},1\}$.  

For any $0<T< + \infty$ we define a time weighted function space $Z^s_{a}(0,T)$ by 
\begin{eqnarray*}
\lefteqn{Z^s_{a}(0,T)=\{u=\trans(\phi,w); u=u_1+u_\infty,}\\
&& u \in C([0,T]; L^2_{(1)}), \ \ P_\infty u=\trans(\phi_\infty, w_{\infty}), \\
&&\phi_\infty \in C([0,T]; H^{s}_{(\infty)}), \ \ w_\infty \in C([0,T]; H^{s}_{(\infty)})\cap L^2(0,T; H^{s+1}_{(\infty)}), \\
&&\|u_1,  u_\infty\|_{Z^s(0,T)} \leq a\}
\end{eqnarray*}
with the norm $\|u_1,  u_\infty\|_{Z^s(0,T)}$  defined by 
\begin{eqnarray*}
\|u_1, u_{\infty}\|_{Z^s(0,T)}&=&
\sup_{0\leq t \leq T}(1+t)^{\delta_{1}}\|  u_1\|_{L^\infty}+
\sup_{0\leq t \leq T}
\sum_{j=0}^{1}(1+t)^{\delta_2+\frac{j}{2}}\|\nabla^j  u_1\|_{L^2}\\
&&\quad 
\sup_{0\leq t \leq T}(1+t)^{\delta_1}\| \nabla  u_1\|_{L^n}
+\sup_{0\leq t \leq T}(1+t)^{\delta_1}\| u_\infty\|_{H^{s}}\\
&&\quad + \sup_{0\leq t \leq T}(1+t)^{\delta_1}\Big(\displaystyle\int_0^t e^{-C_2(t-\tau)}(1+\tau)^{-2\delta_2-1}\|\nabla u_\infty\|_{H^{s-1}\times H^{s}}^2d\tau\Big)^{\frac{1}{2}}, 
\end{eqnarray*}
where 
$a$ and $C_2$ are positive constants independent of $k$ and $T$ and are defined below respectively.  Note that the space $Z^s_{a}(0,T)$ is complete with the norm $\|u_1, u_{\infty}\|_{Z^s(0,T)}$. 

We first consider estimates of a solution $(u_1, u_\infty)$ to \eqref{equation-low} and \eqref{equation-high}:   
\begin{eqnarray}
\delt u_1+A u_1 +B_1[u_\omega]u_1= \mathcal{G}_1(u,u_\omega),
\end{eqnarray}
and 
\begin{eqnarray}
\delt u_\infty+A u_\infty +B_\infty[u_\omega,u^{(1)}]u_\infty= \mathcal{G}_\infty(u,u_\omega). 
\end{eqnarray}

On the low frequency part, let 

$$
u_1(t)=S_1(t)u_{01} +\displaystyle\int_0^t T (t-\tau)\mathcal{G}_1(u,u_\omega)d\tau, 
$$ 
where $S_1(t)$ is the solution operator in  the low frequency part. 
Using the Sobolev embedding, Lemma \ref{lemP_1}, Propositions \ref{analytic semigroup2} and \ref{analytic semigroup3} we see that for  any large $q$ satisfying $n <q$ 
\begin{align*}
\|u_1\|_{L^\infty} &\leq C\|u_1\|_{W^{1,q}}  
\leq C\|u_1\|_{L^q}\\
&\leq  C(1+t)^{-\delta_1}\|u_{01}\|_{L^p} +\displaystyle\int_0^t (1+t-\tau)^{-\delta_1} \|\mathcal{G}_1(u,u_\omega)\|_{L^p_{(1)}} d\tau. 
\end{align*}
Here we used 
the semigroup property as 
$$
\|e^{-tA} f\|_{L^q}  
\leq  C(1+t)^{-\frac{n}{2}\Big(\frac{1}{2}-\frac{1}{q}\Big)+\epsilon} \| e^{-t/2 L_1}f\|_{L^2}. 
$$
We note that for any large $q$ 
$$
1-\dfrac{n}{4}< \frac{n}{2}\Big(\frac{1}{2}-\frac{1}{q}\Big) < \dfrac{n}{4}. 
$$
In addition, Proposition \ref{analytic semigroup2} verifies that 
\begin{align}\label{nonlinear-est-low1}
\|u_1\|_{L^2}  \leq  C(1+t)^{-\delta_2}\|u_{01}\|_{L^p} +\displaystyle\int_0^t (1+t-\tau)^{-\delta_2} \|\mathcal{G}_1(u,u_\omega)\|_{L^p_{(1)}} d\tau, 
\end{align}
and 
\begin{align}\label{nonlinear-est-low2}
\|\nabla u_1\|_{L^2} \leq C(1+t)^{-\delta_2-\frac{1}{2}}\|u_{01}\|_{L^q} +\displaystyle\int_0^t (1+t-\tau)^{-\delta_2-\frac{1}{2}} \|\mathcal{G}_1(u,u_\omega)\|_{L^p_{(1)}} d\tau.  
\end{align}

Furthermore, for $\|\nabla u_1\|_{L^n}$, 
we see from Lemmata \ref{lemP_1} and  \ref{Hardy-Rellich} {\rm (ii)} that 
$$
\|\nabla u_1\|_{L^n} \leq C\|(-\Delta)^{\frac{n}{2}\Big(\frac{1}{2}-\frac{1}{n}\Big)+\frac{1}{2}}u_1 \|_{L^2} \leq C\|(-\Delta)^{\delta_4}u_1 \|_{L^2}, 
$$
where $\delta_4=\min\{\frac{n}{4}-\frac{n}{2q}-\epsilon, 1-\epsilon\}$. 
For any large $q$, it holds that 
$1-\frac{n}{4}<\frac{n}{4}-\frac{n}{2q}< \frac{n}{4}$. 
These together with Corollary  \ref{analytic semigroup2-cor}, the semigroup property and Proposition \ref{analytic semigroup3} verify that 
\begin{align}\label{nonlinear-est-low3}
\|\nabla u_1\|_{L^n} \leq C(1+t)^{-\delta_1}\|u_{01}\|_{L^q} +\displaystyle\int_0^t (1+t-\tau)^{-\delta_1} \|\mathcal{G}_1(u,u_\omega)\|_{L^p_{(1)}} d\tau.  
\end{align}

 On the nonlinear terms  $w_\omega \cdot \nabla w_\infty$, 
 we rewrite to 
 \begin{eqnarray}\label{derivative-form}
w_\omega \cdot \nabla w_\infty = 
\nabla (w_\omega w_\infty) - w_\infty\cdot \nabla w_\omega. 
 \end{eqnarray}
 On the first term we apply Corollary \ref{analytic semigroup3-divergence-form}, {\rm i.e., } we get the second terms of the right hand sides in \eqref{nonlinear-est-low1} and \eqref{nonlinear-est-low2}  as 
 \begin{align}\label{one point-estimate}
\displaystyle\int_0^t (1+t-\tau)^{-\delta_3-\frac{1+m}{2}} \|w_\omega w_\infty\|_{L^{p_2}} d\tau 
\end{align}
with  $m=0$ for $\|u_1\|_{L^2}$ and $m=1$ for $\|\nabla u_1\|_{L^2}$, respectively,  
where $\delta_3=\frac{n}{2}\Big(\frac{1}{p_2}-\frac{1}{2}\Big)-\epsilon$, $1< p_2 <2$ satisfying that 
$$
\dfrac{1}{p_2}= \frac{1}{2}-\dfrac{1}{n} +\dfrac{1}{p_3}.  
$$
Here we take $2< p_3 <n$ close to $2$ which satisfies that $\delta_3 + \frac{1}{2}= \frac{n}{2p_3}> \delta_2 $. The Sobolev embedding implies that 
$$
\|w_\omega w_\infty\|_{L^{p_2}} \leq 
C\|\nabla w_\omega\|_{L^2}\|w_\infty\|_{H^s}\leq C\epsilon (1+t)^{-\delta_1}\|u_1, u_\infty\|_{Z^s(0,T)}. 
$$
The nonlinear terms $\|\mathcal{G}_1(u,u_\omega)\|_{L^p_{(1)}}$ in \eqref{nonlinear-est-low1} and \eqref{nonlinear-est-low1} can be estimated by the properties of the low and high frequency parts,  Lemma \ref{lemP_1} and Lemma \ref{lemPinfty}, as 
\begin{align}\label{low-nonlinear-1}
\begin{aligned}
\|w \cdot \nabla w\|_{L^p} \leq C\|\nabla w\|_{L^2}\|w\|_{L^{p_0}} 
& \leq C\|\nabla w\|_{L^2}(\|(-\Delta)^k w_1\|_{L^2}+\|(-\Delta)^k w_\infty\|_{L^{2}})\\
& \leq  C\|\nabla w\|_{L^2}( \|w_1\|_{L^2}+\| w_\infty\|_{H^{s}}) \\
&\leq C(1+t)^{-2\delta_2-1/2}\|u_1, u_\infty\|_{Z^s(0,T)}^2, 
\end{aligned}
\end{align}  
 where $2<p_0$ satisfying  $1/p=1/2+1/p_0$ and $k=\frac{n}{2}\big(\frac{1}{2}-\frac{1}{p_0}\big)$ with $2k \leq s$.    
We see from the same argument as that in \eqref{low-nonlinear-1} that the second term in \eqref{derivative-form} is estimated by 
\begin{align}\label{low-nonlinear-2}
\begin{aligned}
\|w_\infty \cdot \nabla w_\omega\|_{L^p} &\leq C\|\nabla w_\omega\|_{L^2}\| w_\infty\|_{H^{s}} \\
&\leq C\epsilon (1+t)^{-\delta_1}\|u_1, u_\infty\|_{Z^s(0,T)}.  
\end{aligned}
\end{align}  
Concerning the nonlinear term $h^{(2)}(\phi,\phi_\omega)\phi \phi_\omega \nabla \phi_\omega$, 
The H\"older inequality 
yields 
\begin{align}
\begin{aligned}
\|h^{(2)}(\phi,\phi_\omega)\phi \phi_\omega \nabla \phi_\omega\|_{L^{p}} &\leq 
C\|h^{(2)}(\phi,\phi_\omega)\phi_\omega\|_{L^{p_0}}
\|\nabla \phi_\omega\|_{L^2}\|\phi\|_{L^\infty}\\
& \leq C \|h^{(2)}(\phi,\phi_\omega)\|_{L^{p_{4}}}\|\phi_\omega\|_{L^{2^*}}
\|\nabla \phi_\omega\|_{L^2}\|\phi\|_{L^\infty}\\
& \leq C \|h^{(2)}(\phi,\phi_\omega)\|_{L^{p_{4}}}\|\nabla \phi_\omega\|_{L^{2}}
\|\nabla \phi_\omega\|_{L^2}\|\phi\|_{L^\infty}\\
&\leq C\epsilon \|h^{(2)}(\phi,\phi_\omega)\|_{L^{p_{4}}}\|\phi\|_{L^\infty}, 
\end{aligned}
\end{align}
where $1/p=1/2+1/p_0=1-1/n + 1/p_4$. 
When $n=3$, since $1/p< n/4=3/4$, $1/p_4 < 1/12$, {\rm i.e., } $p_4>12$. 
When $4\leq n$, since $1/p <1$, $1/p_4< 1/n$ {\rm i.e., } $p_4>n$. We note that 
$\frac{1}{2}< \frac{n}{2}\Big(\frac{1}{2}-\frac{1}{q_4}\Big)<\frac{n}{4}$. 
Hence the Sobolev embedding, Lemma \ref{Hardy-Rellich} and the interpolation imply that for 
$2+[n/2] \leq s$ and some $0<\theta<1$, 
$$
\|h^{(2)}(\phi,\phi_\omega)\|_{L^{p_{4}}}
\leq C\|\nabla h^{(2)}(\phi,\phi_\omega)\|_{L^2}^{1-\theta}\|\nabla h^{(2)}(\phi,\phi_\omega)\|_{H^{s-1}}^{\theta}
\leq C(\|\nabla \phi_\omega\|_{H^{s-1}}+\|\nabla \phi\|_{H^{s-1}}). 
$$
Hence we obtain that 
\begin{align}
\begin{aligned}
\|h^{(2)}(\phi,\phi_\omega)\phi \phi_\omega \nabla \phi_\omega\|_{L^{p,2}} \leq C\epsilon (1+t)^{-\delta_1}\|u_1, u_\infty\|_{Z^s(0,T)}. 
\end{aligned}
\end{align}
 
Since other terms in $\mathcal{G}_1(u,u_\omega)$ can be estimated similarly,  we obtain that 

\vspace{2ex}

\begin{prop}\label{est-low-result}
Let $$
u_1(t)=S_1 (t)u_{01} +\displaystyle\int_0^t S_1 (t-\tau)\mathcal{G}_1(u,u_\omega)d\tau. 
$$ 
Then it holds that 
\begin{align*}
   &\sup_{0\leq t \leq T}(1+t)^{\delta_1} (\|  u_1\|_{L^q}+\|\nabla u_1\|_{L^n})+
\sup_{0\leq t \leq T}
\sum_{j=0}^{1}(1+t)^{\delta_2+\frac{j}{2}}\|\nabla^j  u_1\|_{L^2}
\\
& \leq C_0 E_0 +C_1 \|u_1, u_\infty\|_{Z^s(0,T)}^2 
+C\epsilon \|u_1, u_\infty\|_{Z^s(0,T)}, 
\end{align*}
where $C_0$ and $C_1$ are independent of $T$. 
\end{prop}

\vspace{2ex}

Concerning $u_\infty$, we use the estimate \eqref{energy2}.  
The following estimate which is related to the estimate of the nonlinearity $\mathcal{G}_\infty(u,u_\omega)$ is obtained by 
Lemmata \ref{lemP_1},  \ref{lemPinfty} and \ref{composition-est} as follows.

\vspace{2ex} 

\begin{lem}\label{nonlinearest-high}
It holds that for $t\in [0,T]$ and $\|u_1, u_\infty\|_{Z^s(0,T)} \leq a_0$ 
\begin{align*}
&\lefteqn{\|\mathcal{G}_\infty(u,u_\omega)\|_{H^{s}\times H^{s-1}}}\\
&\leq C(1+t)^{-2\delta_2-\frac{1}{2}}\|u_1, u_\infty\|_{Z^s(0,T)}^2
+ C\epsilon(1+t)^{-\delta_1}\|u_1, u_\infty\|_{Z^s(0,T)}\\
&\quad +C(1+t)^{-\delta_2-\frac{1}{2}}\|u_1, u_\infty\|_{Z^s(0,T)}\|\nabla u_\infty\|_{H^{s-1}\times H^s}.
\end{align*}
\end{lem}

\vspace{2ex} 

\noindent\textbf{Proof.} We write details of  estimates of  $P_\infty (\phi \Delta w)$, $P_\infty (h^{(1)}(\phi_\omega)  \phi_\omega \nabla \phi_1 )$, $P_\infty  ( w_1 \cdot \nabla \phi_\omega )$, $P_\infty (w \cdot \nabla w)$ 
$P_\infty(\phi \phi_\omega \nabla \phi_\omega)$ and $P_\infty(\phi \phi_\omega \Delta w_\omega)$. For $0 \leq |\alpha|\leq s-1$ and $1 \leq |\beta| \leq s-1$ we see from Lemmata \ref{lem2.1.}- \ref{composition-est} that  
\begin{eqnarray*}
\|\del_x^{\alpha }P_\infty (\phi \Delta w)\|_{L^2} 
&\leq& C\|\nabla \phi\|_{H^{s-1}}\|\Delta w\|_{H^{s-1}}\\
&\leq &   C(1+t)^{-\delta_2-\frac{1}{2}}\|u_1, u_\infty\|_{Z^s(0,T)}\|\nabla u_\infty\|_{H^{s-1}\times H^s}, \\
\|P_\infty (h^{(1)}(\phi_\omega)  \phi_\omega \nabla \phi_1 )\|_{L^2} 
&\leq& C\|\phi_\omega\|_{L^{2^*}} \|\nabla\phi_1\|_{L^n } \\
& \leq& C\|\nabla \phi_\omega\|_{L^{2}} \|\nabla\phi_1\|_{L^n} \leq C\epsilon(1+t)^{-\delta_1}\|u_1, u_\infty\|_{Z^s(0,T)}, \\
\|\del_x^{\beta} P_\infty (h^{(1)}(\phi_\omega)  \phi_\omega \nabla \phi_1 )\|_{L^2} 
&\leq& C\|\phi_\omega\|_{L^{2^*}} \|\nabla\phi_1\|_{L^n }+ C\|\nabla \phi_\omega\|_{H^{s-1}} \|\nabla \phi_1\|_{W^{s-1, \infty}}
\\
& \leq& C\|\nabla \phi_\omega\|_{L^{2}} \|\nabla\phi_1\|_{L^n} 
+C\|\nabla \phi_\omega\|_{H^{s-1}}\|\nabla\phi_1\|_{L^n}
\\
&\leq& C\epsilon(1+t)^{-\delta_1}\|u_1, u_\infty\|_{Z^s(0,T)}, \\
\|\del_x^{\alpha }P_\infty (\phi_1\cdot \nabla w_\omega )\|_{L^2} 
&\leq& C\|\nabla w_\omega\|_{H^{s-1}}\|\phi_1\|_{W^{s-1, \infty}}\\
& \leq& 
C\epsilon (\|\phi_1\|_{L^\infty}+\|\nabla \phi_1\|_{L^n})\\
&\leq& C\epsilon(1+t)^{-\delta_1}\|u_1, u_\infty\|_{Z^s(0,T)}, \\
\|\del_x^{\alpha }P_\infty (w \cdot \nabla w )\|_{L^2} 
&\leq& C \|\nabla w\|_{H^{s-1}} \|\nabla w\|_{H^{s-1}}\leq  C(1+t)^{-2\delta_2-\frac{1}{2}}\|u_1, u_\infty\|_{Z^s(0,T)}^2, \\
\|\del_x^{\alpha }P_\infty (\phi \phi_\omega  \nabla \phi_\omega) \|_{L^2} 
&\leq& C (\|\phi_1\|_{W^{s-1, \infty}}+\|\nabla \phi_\infty\|_{H^{s-1}})\|\nabla \phi_\omega \|_{H^{s-1}}\\
& \leq& C (\|\phi_1\|_{L^\infty}+\|\nabla \phi_1\|_{L^n}+\|\nabla \phi_\infty\|_{H^{s-1}})\|\nabla \phi_\omega \|_{H^{s-1}}\\
&\leq&  C\epsilon (1+t)^{-\delta_1 }\|u_1, u_\infty\|_{Z^s(0,T)}, \\
\|\del_x^{\alpha }P_\infty (\phi \phi_\omega  \Delta w_\omega )\|_{L^2} 
&\leq& C (\|\phi_1\|_{W^{s-1, \infty}}+\|\nabla \phi_\infty\|_{H^{s-1}})\|\nabla w_\omega \|_{H^{s}}\\
&\leq& C( \|\phi_1\|_{L^\infty}+\|\nabla \phi_1\|_{L^n}+\|\nabla \phi_\infty\|_{H^{s-1}} )\|\nabla w_\omega \|_{H^{s}}\\
& \leq&  C\epsilon (1+t)^{-\delta_1}\|u_1, u_\infty\|_{Z^s(0,T)}, 
\end{eqnarray*}
where $1/2^* = 1/2 - 1/n$. 
Since other nonlinear terms can be estimated similarly, we get Lemma \ref{nonlinearest-high}. 
$\hfill\square$

\vspace{2ex}

Let $D[u_\infty]=\|\nabla \phi_{\infty}(t)\|_{H^{s-1}}^{2}+\|\nabla w_{\infty}(t)\|_{H^{s}}^{2}$. 
By \eqref{energy2} and Lemma \ref{nonlinearest-high}, 
there exists a positive constant $C_2$ such that  for $t \in [0,T]$ and $\|u_1, u_\infty\|_{Z^s(0,T)} \leq a_0$
\begin{align}
\begin{aligned}\label{estenergy}
\lefteqn{{\cal E}[u_{\infty}](t)+d \displaystyle\int_{0}^{t}e^{-C_2(t-\tau)}D[u_\infty](\tau)d\tau}\\
&\quad \leq e^{-C_2 t}{\cal E}[u_{\infty}](0) \\
&\qquad  +C\|u_1, u_\infty\|_{Z^s(0,T)}^4 \displaystyle\int_{0}^{t}
e^{-C_2(t-\tau)}(1+\tau)^{-4\delta_2 -1}d\tau\\
&\qquad  +C\epsilon^2\|u_1, u_\infty\|_{Z^s(0,T)}^2 \displaystyle\int_{0}^{t}
e^{-C_2(t-\tau)}(1+\tau)^{-2\delta_2 -1}
d\tau \\
&\qquad 
+C\|u_1, u_\infty\|_{Z^s(0,T)}^2 \displaystyle\int_{0}^{t}
e^{-C_2(t-\tau)}(1+\tau)^{-2\delta_2 -1}D[u_\infty](\tau)
d\tau\\
&\quad \leq e^{-C_2 t}{\cal E}[u_{\infty}](0)\\
&\qquad  +C\|u_1, u_\infty\|_{Z^s(0,T)}^4 (1+t)^{-4\delta_2-1} \\
&\qquad  +C\epsilon^2\|u_1, u_\infty\|_{Z^s(0,T)}^2 (1+t)^{-2\delta_1} \\
&\qquad 
+C\|u_1, u_\infty\|_{Z^s(0,T)}^2 \displaystyle\int_{0}^{t}
e^{-C_2(t-\tau)}(1+\tau)^{-2\delta_2 -1}D[u_\infty](\tau)
d\tau. 
\end{aligned}
\end{align}

Define ${\cal D}[u_\infty]$ and $\tilde{\cal E}[u_\infty]$ by 
\begin{eqnarray*}
{\cal D}[u_\infty](t)&=&(1+t)^{2\delta_1}\displaystyle\int_{0}^{t}
e^{-C_2(t-\tau)}(1+\tau)^{-2\delta_2-1}D[u_\infty](\tau)d\tau, \\
\tilde{\cal E}[u_\infty](t) &=&\sup_{0\leq \tau \leq t}(1+\tau)^{2\delta_1}{\cal E}[u_{\infty}](\tau)
\end{eqnarray*}
We see from \eqref{estenergy} that 
\begin{eqnarray}
\tilde{\cal E}[u_\infty](t)+d{\cal D}[u_\infty](t) &\leq& C(\tilde{\cal E}[u_\infty](0)+\|u_1, u_\infty\|_{Z^s(0,T)}^4+C\|u_1, u_\infty\|_{Z^s(0,T)}^2{\cal D}[u_\infty](t))\nonumber\\
&&\quad + C\epsilon^2\|u_1, u_\infty\|_{Z^s(0,T)}^2\nonumber\\
&\leq & C_0^2E_0^2+ C^2\|u_1, u_\infty\|_{Z^s(0,T)}^4 
 + C\epsilon^2\|u_1, u_\infty\|_{Z^s(0,T)}^2,\label{est-high-result}
\end{eqnarray}
where the constant $C_3$ is independent of $T$. 
\vspace{2ex}

To apply Banach's fixed point theorem we use the iteration in 
$$ U= (u_1,u_\infty)= (u_1, u_\infty) \in Z^s_{a}(0,T)$$ 
defined by 
$U^{(1)}\mapsto U^{(2)} = H[U^{(1)}] $ for any $U^{(1)} \in 
Z^s_{a}(0,T)$, 
where the mapping $H$ is a solution operator of \eqref{equation-low}-\eqref{equation-high}, {\rm i.e., }
\begin{eqnarray}
\delt u^{(2)}_1+A u^{(2)}_1 +B_1[u_\omega]u^{(2)}_1= \mathcal{G}_1(u^{(1)},u_\omega),
\end{eqnarray}
and 
\begin{eqnarray}
\delt u^{(2)}_\infty+A u^{(2)}_\infty +B_\infty[u_\omega,u^{(1)}]u^{(2)}_\infty= \mathcal{G}_\infty(u^{(1)},u_\omega). 
\end{eqnarray}

We see from Proposition \ref{est-low-result} in the low frequency part,  Lemma \eqref{nonlinearest-high} with \eqref{est-high-result} in the high frequency part that  for $a \leq a_0$  
\begin{align}\label{U1U2}
\begin{aligned}
\|U^{(2)}\|_{Z^s_{a}(0,T)} & \leq C_0 E_0+C_1\|U^{(1)}\|^2_{Z^s_{a}(0,T)} + C_2\epsilon^2 \|U^{(1)}\|_{Z^s_{a}(0,T)}, 
\end{aligned}
\end{align}
where $C_0$, $C_1$ and $C_2$ are independent of $E_0$, $\epsilon$, $a$ and $T$.  
We note that if $a \leq \min\{a_0, 1, 1/2C_1\}$,  
and 
$E_0, \epsilon $ are small such that 
$$
c_0 E_0+ C_1 \epsilon \leq \frac{1}{2}, 
$$
then 
$\|U^{(2)}\|_{Z^s_{a}(0,T)} \leq a$. 
Moreover, given any $U^{(1)}, V^{(1)} \in Z^s_{a}(0,T)$ with solution  $U^{(2)}, V^{(2)} \in Z^s_{a}(0,T)$, respectively, we get the estimate
\begin{align}\label{U1-U2}
\|U^{(2)} - V^{(2)} \|_{Z^s_{a}(0,T)}  
& \leq C_3\big(\|U^{(1)}\|_{Z^s_{a}(0,T)}  + \|V^{(1)}\|_{Z^s_{a}(0,T)} + \epsilon \big)
   \|U^{(1)} - V^{(1)} \|_{Z^s_{a}(0,T)}, 
\end{align}
where the constant $C_3$ is independent of $T$ and $a$. 
By \eqref{U1U2}, \eqref{U1-U2} there exists a small $a$ such that the map $H[\cdot]$ is a strict contraction. Hence there exists a locally unique fixed point $U=(u_1, u_\infty)$ in $Z^s_{a}(0,T)$. Since constants are independent of $T$, letting $T$ go to $\infty$ and the definition of the time weighted space $Z^s_{a}(0,T)$, we also obtain the decay rates in Theorem \ref{stability} as $\delta=\delta_1$ with any large $q$. This completes  the proof.


\end{section}

\begin{section}{Application to the damped wave equations}\label{application}
\begin{subsection}{Formulation} 








To prove Theorem \ref{decay-DWE}, let us reformulate \eqref{eq-DWE}. We set $v=\del_t u$. Then \eqref{decay-DWE}  is rewritten to the following  equation 
\begin{eqnarray}\label{DEW-evolution}
    \del_t U +AU + \mathcal{B}[b]U =0,  \ \ U(0,x)=U_0, 
\end{eqnarray}
where $U=\trans(u, v)$, 
\begin{eqnarray*}
A&=&\begin{pmatrix}
0 &1 \\
-\mu\Delta  & -\mu'\Delta 
\end{pmatrix}
, \\[1ex]
\mathcal{B}[b]U &=&
\begin{pmatrix}
0    \\
 u \div b_1  + b_2 \cdot \nabla  u  + b_3 \Delta u
\end{pmatrix}
.
\end{eqnarray*}
Applying $P_1$ and $P_\infty$ to \eqref{DEW-evolution} and letting 
$U_1= P_1 U$ and $U_\infty = P_\infty U$, we obtain that 
\begin{eqnarray}\label{DEW-evolution-low}
    \del_t U_1 + L_1 U_1  = -P_1 \mathcal{B}[b]U_\infty,  \ \ U_1(0,x)=P_1 U_0 
\end{eqnarray}
and 
\begin{eqnarray}\label{DEW-evolution-high}
    \del_t U_\infty + L_\infty U_\infty  = -P_\infty \mathcal{B}[b]U_1,  \ \ U_\infty(0,x)=P_\infty U_0,  
\end{eqnarray}
where $L_1 = A_1 + P_1 \mathcal{B}[b]$, $A_1=A|_{L^2_{(1)}}$, $L_\infty =  A_\infty + P_\infty \mathcal{B}[b]$, $A_\infty= A|_{L^2_{(\infty)}}$. 
\end{subsection}

\begin{subsection}{Estimates of the low frequency part}

In this subsection. we estimate a solution to \eqref{DEW-evolution-low}. In this aim, we linearize \eqref{DEW-evolution-low} as 
\begin{eqnarray}\label{DEW-evolution-low-1}
    \del_t U_1 + L_1 U_1  = 0,  \ \ U_1(0,x)=U_{01}.  
\end{eqnarray}
We first consider the linear problem 
\begin{eqnarray}\label{DEW-evolution-low-1-A}
\del_t U_1 + A_1 U_1  = 0,  \ \ U_1(0,x)=U_{01}. 
\end{eqnarray}
The resolvent problem of \eqref{DEW-evolution-low-1-A} is given by 
\begin{eqnarray}\label{DEW-evolution-low-resolvent}
\left\{
\begin{array}{lll}
\lambda u_1 - v_1   =f_{11},\\
\lambda^2 u_1 - \mu \Delta \lambda u_1 - \mu' \Delta u_1  =f_{12}.  
\end{array}
\right.
\end{eqnarray}
 Applying the Fourier transform to \eqref{DEW-evolution-low-resolvent} yields the fundamental solution formula 
\begin{eqnarray}
\hat{u}_1 (t,\xi)&=&\displaystyle\frac{ \hat{f}_{12}}{(\lambda-\lambda_+)(\lambda-\lambda_-)},\nonumber\\
\hat{v}_1(t,\xi)&=&-\hat{f}_{11} 
+\displaystyle\frac{ \lambda \hat{f}_{12}}{(\lambda-\lambda_+)(\lambda-\lambda_-)}, 
\label{solution formula1-DEW}
\end{eqnarray} 
where $\lambda_{\pm}=-\dfrac{\mu}{2}|\xi|^2 \pm  \dfrac{1}{2}\sqrt{\mu^2 |\xi|^4- 4 \mu' |\xi|^2} \sim 
-\dfrac{\mu}{2}|\xi|^2 \pm \sqrt{\mu'}|\xi|i $ as $|\xi| \ll 1$ and $|\xi | \neq 0, \dfrac{2\mu'}{\sqrt{\mu}}$.  Hereafter we assume that $|\xi| \leq r_\infty < \dfrac{2\mu'}{\sqrt{\mu}}$. 
Let 
\begin{eqnarray}
\hat{A}_{1}=\begin{pmatrix}
0 &1\\
\mu |\xi|^2 & \mu' |\xi|^2
\end{pmatrix}
\ \ \ (\xi\in \mathbb{R}^n). 
\nonumber
\end{eqnarray}

\vspace{2ex}

\begin{prop}\label{resolvent2-each-point-DEW} 
Let $\lambda \in \mathbb{C}\setminus \sigma (-\hat{A}_1)$.  We suppose that $\mu \geq 5 \mu'$ with $\mu' \geq 1$.   Let a set $\mathcal R_1$ be defined as 
\begin{align*}
\mathcal{R}_{1, + } &=\{\lambda \in \mathbb{C}\setminus \sigma (-\hat{A}_1); 
\Re \lambda <0,  \   \lambda = -a^2 + (a + c_0) i \mbox{ with } a>0 
\}, \\
\mathcal{R}_{1, - } &=\{\lambda \in \mathbb{C}\setminus \sigma (-\hat{A}_1); 
\Re \lambda <0,  \   \lambda = -a^2 - (a + c_0) i \mbox{ with } a>0 
\}. 
\end{align*}
Further let the set $\mathcal R_2$ be defined as 
$$
\mathcal{R}_2 =\{\lambda \in \mathbb{C}\setminus \sigma (-\hat{A}_1); 
\Re \lambda >0
\}.
$$
 We consider the following two cases. 

{\rm Case 1}: $\lambda\in \mathcal{R}_{1, \pm}$,  \  \ {\rm Case 2}: $\lambda \in \mathcal R_2$. 
\BLACK





Then in each case the solution $\hat{U}_1(\xi)=\trans(\hat{u}_1, \hat{v}_1) \in L^2_{(1)}$ to \eqref{solution formula1-DEW} satisfies 
for $F_1=\trans(f_{11}, f_{12} )$ the estimate \BLACK
$$
|\hat{u}_1| \leq \dfrac{C}{|\lambda||\xi|^2} |\hat{f}_{12}|,  \ \ ||\xi|^2 \hat{u}_1| \leq \dfrac{C}{|\xi|}|\hat{f}_{12}|
$$
and 
$$
|\hat{v}_1| \leq C|\hat{f}_{11}|+\dfrac{C}{|\lambda|}|\hat{f}_{12}|,  \ \  
||\xi|^2 \hat{v}_1| \leq C  (|\hat{f}_{11}|+|\hat{f}_{12}|),  
$$
where positive constants $C$ are independent of $c_0$, $\lambda$ and $\xi$. \BLACK  
\end{prop}

\vspace{2ex}

Since 
$$
\Big|\dfrac{1}{ (\lambda-\lambda_+)(\lambda-\lambda_-)}\Big| \sim  \Big|\dfrac{1}{|\xi|}\frac{ |\xi|}{(\lambda-\lambda_+)(\lambda-\lambda_-)}\Big|
$$
as $|\xi| \rightarrow 0$, Proposition \ref{resolvent2-each-point-DEW} is directly obtained by the proof of Proposoition \ref{resolvent2}. 

\vspace{2ex}

For $5 \leq n$, we note that $1/|\xi|$ and $1/|\xi|^2$ are $L^2$ integrable on $|\xi| \leq r_\infty \leq 1$. Therefore, we have 

\vspace{2ex}

\begin{prop}\label{resolvent2-DEW} 
Let $\lambda \in \mathbb{C}\setminus \sigma (-{A}_1)$ as in Proposition \ref{resolvent2-each-point-DEW}. 
Then in each case the solution $U_1 \in L^2_{(1)}$ to \eqref{DEW-evolution-low-1-A} satisfies 
for $F_1=\trans(f_{11}, f_{12})$ the estimate \BLACK
$$
\|U_1\|_{L^2_{(1)}} \leq C(\|f_{11}\|+\frac{1}{|\lambda|}\|f_{12}\|_{L^2})
$$
and 
$$
\|\nabla^2 U_1\|_{L^2_{(1)}} \leq C(\|f_{11}\|+\|f_{12}\|_{L^2}), 
$$
where positive constants $C$ are independent of $c_0$ and  $\lambda$. \BLACK  
\end{prop}

\vspace{2ex}

By a similar argument to Corollary \ref{resovent-estimate-modify-0}, we have the following resolvent estimates in the Fourier space. 

\vspace{2ex}

\begin{cor}\label{resovent-estimate-DWE-modify-0}
Let $\lambda\in \mathcal R_{1,\pm}\cup \mathcal R_2$ as in Proposition \ref{resolvent2-each-point-DEW}. 
 
Then  the following estimates hold for any small epsilon $\epsilon$ independent of $c_0$, $\lambda$ and $\xi$. 
$$
|\hat{u}_1| + |\hat{v}_1| \leq C_n|\hat{f}_{11}|+\dfrac{C}{|\lambda|^{1/2 +\epsilon}|\xi|^{5/2-\epsilon}}|\hat{f}_{12}|
$$
and 
$$
||\xi|\hat{u}_1| + ||\xi| \hat{v}_2|\leq C(|\hat{f}_{11}|+|\hat{f}_{12}|), 
$$
where positive constants $C$ are independent of $c_0$, $\lambda$ and $\xi$.  
\end{cor}

\vspace{2ex}

As a result, we also obtain the following resolvent estimates in $L^2_{(1)}$.   

\vspace{2ex}

\begin{prop}\label{resolvent2-DEW-modify} 
Let $\lambda \in \mathbb{C}\setminus \sigma (-{A}_1)$ as in Proposition \ref{resolvent2-each-point-DEW}. 
Then in each case the solution $U_1 \in L^2_{(1)}$ to \eqref{DEW-evolution-low-1-A} satisfies 
for $F_1=\trans(f_{11}, f_{12})$ the estimate \BLACK
$$
\|U_1\|_{L^2_{(1)}} \leq C(\|f_{11}\|+\frac{1}{|\lambda|^{1/2+\epsilon}}\|f_{12}\|_{L^2})
$$
and 
$$
\|\nabla U_1\|_{L^2_{(1)}} \leq C(\|f_{11}\|+\|f_{12}\|_{L^2}), 
$$
where positive constants $C$ are independent of $c_0$ and $\lambda$. \BLACK  
\end{prop}

\vspace{2ex}

The Plancherel theorem and interpolation with Proposition \ref{resolvent2-DEW} directly yield the following corollary. 

\vspace{2ex}

\begin{cor}\label{resolvent2-cor-DWE}
Let $\lambda \in \mathbb{C}\setminus  \sigma (-A_1)$ as in Proposition \ref{resolvent2-DEW}.   Then 
for $F_1=\trans(f_{11}, f_{12})\in L^2_{(1)}$ 
$$
\| (-\Delta)^{\delta}(\lambda +A_1)^{-1}  F_1 \|_{L^2}\leq C_{n}\|f_{11}\|_{L^2}+C_n |\lambda|^{\delta-1}\|f_{12} \|_{L^2}, 
$$
for $0 \leq \delta \leq 1$ and 
$$
\| (-\Delta)^{\frac{k}{2}}(\lambda +A_1)^{-1} (-\Delta)^{1-\frac{k}{2}} F_1 \|_{L^2}\leq C_{n}\{\|f_{11}\|_{L^2}+\|f_{12} \|_{L^2}\}, 
$$
where $0 \leq k \leq 2$. 
\end{cor}

\vspace{2ex}

We also obtain $L^2$ estimates of the adjoint operator as follows. 

\vspace{2ex}

\begin{cor}\label{resolvent2-each-point-DEW-cor} 
Let $\lambda \in \mathbb{C}\setminus \sigma (-{A}_1)$ as in Proposition \ref{resolvent2-each-point-DEW}. 
Then it holds that 
$$
\|(\lambda+A_1^*)^{-1}F_1 \|_{L^2_{(1)}} \leq C(\|f_{11}\|+\frac{1}{|\lambda|}\|f_{12}\|_{L^2})
$$
and 
$$
\|\nabla^2 (\lambda+A_1^*)^{-1}F_1\|_{L^2_{(1)}} \leq C(\|f_{11}\|+\|f_{12}\|_{L^2}), 
$$
where positive constants $C$ are independent of $c_0$, $\lambda$ and $\xi$. 

\end{cor}

\vspace{2ex}

\begin{cor}\label{resolvent2-cor2-DWE}
Let $\lambda \in \mathbb{C}\setminus  \sigma (-{A}_1)$ as in Proposition \ref{resolvent2-DEW}.   Then 
for $F_1=\trans(f_{11}, f_{12})\in L^2_{(1)}$ 
$$
\| (-\Delta)^{\delta}(\lambda +A_1^*)^{-1}  F_1 \|_{L^2}\leq C_{n}\|f_{11}\|_{L^2}+C_n |\lambda|^{\delta-1}\|f_{12} \|_{L^2}, 
$$
for $0 \leq \delta \leq 1$ and 
$$
\| (-\Delta)^{\frac{k}{2}}(\lambda +A_1^*)^{-1} (-\Delta)^{1-\frac{k}{2}} F_1\|_{L^2}\leq C_{n}\{\|f_{11}\|_{L^2}+\|f_{12} \|_{L^2}\}, 
$$
where $0 \leq k \leq 2$.
\end{cor}

\vspace{2ex}

Corollary \ref{resolvent2-each-point-DEW-cor} and \ref{resolvent2-cor2-DWE} directly follow from Proposition \ref{resolvent2-each-point-DEW-cor}, Corollary \ref{resolvent2-cor-DWE}, the Plancherel theorem and a duality argument. 

\vspace{2ex}

Similarly, we see from Proposition \ref{resolvent2-DEW-modify}  the following specific resolvent estimates. 

\vspace{2ex}

\begin{cor}\label{resolvent2-cor-DWE-mpdify}
Let $\lambda \in \mathcal{R}_{1,\pm}$ as in Proposition \ref{resolvent2-DEW}.  

\vspace{1ex}

{\rm (i)} For $F_1=\trans(f_{11}, f_{12})\in L^2_{(1)}$ 
$$
\| (-\Delta)^{\delta}(\lambda +A_1)^{-1}  F_1 \|_{L^2}\leq C_{n}\|f_{11}\|_{L^2}+C_n |\lambda|^{\delta-1/2-\epsilon}\|f_{12} \|_{L^2}, 
$$
for $0 \leq \delta <1/2$ with any small $\epsilon>0$ independent of $\lambda$ and $\delta$  and 
$$
\| (-\Delta)^{1/2}(\lambda +A_1)^{-1}  F_1 \|_{L^2}\leq C_{n}\{\|f_{11}\|_{L^2}+\|f_{12} \|_{L^2}\}. 
$$ 

\vspace{1ex}

{\rm (ii)} It holds that 
$$
\| (-\Delta)^{\delta}(\lambda +A_1^*)^{-1}  F_1 \|_{L^2}\leq C_{n}\|f_{11}\|_{L^2}+C_n |\lambda|^{\delta-1/2-\epsilon}\|f_{12} \|_{L^2}, 
$$
for $0 \leq \delta <1/2$ with any small $\epsilon>0$ independent of $\lambda$ and $\delta$  and 
$$
\| (-\Delta)^{1/2}(\lambda +A_1^*)^{-1}  F_1 \|_{L^2}\leq C_{n}\{\|f_{11}\|_{L^2}+\|f_{12} \|_{L^2}\}, 
$$
where positive constants $C$ are independent of $c_0$ and $\lambda$. 
\end{cor}

\vspace{2ex}

We next consider the perturbed operator $L_1$. We prepare the following lemma. 

\vspace{2ex}

\begin{lem}\label{L1-B1-DWE}
Let $k$ satisfy $\max\{0, 1-n/4 \}< k/2 < \min\{1, n/4\}$.  

{\rm (i)}  It holds that for $u \in L^2_{(1)}$ 
\begin{align*}
\begin{aligned}
&\|(-\Delta)^{\frac{k}{2}} (\lambda+A_1)^{-1}B_1[b] u \|_{L^2} 
\leq C \epsilon \|(-\Delta)^{\frac{k}{2}} u\|_{L^2},  \\
&\|(-\Delta)^{\frac{k}{2}} (\lambda+A^*_1)^{-1}B_1^*[b] u \|_{L^2} 
\leq C \epsilon \|(-\Delta)^{\frac{k}{2}} u\|_{L^2}. 
\end{aligned}
\end{align*}

{\rm (ii)}  There holds that for $u \in L^2_{(1)}$ 
\begin{align*}
\begin{aligned}
&\|(\lambda+A_1)^{-1}B_1[b] u \|_{L^2} 
\leq C \epsilon \|\nabla u\|_{L^2}, \\
&\| (\lambda+A^*_1)^{-1}B_1^*[b] u \|_{L^2} 
\leq C \epsilon \|\nabla u\|_{L^2}. 
\end{aligned}
\end{align*}

\end{lem}
\vspace{2ex}
By virtue of Assumption \ref{ass-DWE} and Corollary \ref{resolvent2-each-point-DEW-cor} and \ref{resolvent2-cor2-DWE},  a similar argument to Lemma \ref{L1-B1} verifies Lemma \ref{L1-B1-DWE}. 
We note that by the property of the low frequency part,  Lemma \ref{lemP_1}, the perturbation terms $\mathcal{B}[b]u_1$ including $b_3 \Delta u_1$ satisfies the same estimates.

\vspace{2ex}

Substituting Lemma \ref{L1-B1-DWE} to the forms \begin{align}\label{L1-form-DWE}
(\lambda +L_1)^{-1}=\sum_{j=0}^{\infty}\{(\lambda+A_1)^{-1}B_1[b]\}^{j}(\lambda+A_1)^{-1}
\end{align}
and
\begin{align}\label{L1-form-adjoint-DWE}
(\lambda +L^*_1)^{-1}=\sum_{j=0}^{\infty}\{(\lambda+A^*_1)^{-1}B^*_1[b]\}^{j}(\lambda+A^*_1)^{-1}
\end{align}
respectively and taking $\epsilon$ small, we have 

\vspace{2ex}

\begin{prop}\label{resolvent-est-L1-DWE}
Let $\lambda \in \mathbb{C}\setminus  \sigma (-{A}_1)$ as in Proposition {\rm \ref{resolvent2}} and $u_1=\trans(\phi_1,w_1)$ be a solution to \eqref{cns-linear-low-per2} for $F_1=\trans(f_{11}, f_{12}) \in L^2_{(1)}$.  Then it holds that 
\begin{align}
\begin{aligned}\label{resolvent-est-L1-1-DWE}
&\|(-\Delta)^{\frac{k}{2}} (\lambda+ L_1)^{-1} F_1\|_{L^2} \leq C\|f_{11}\|_{L^2}+ 
C|\lambda|^{\frac{k}{2}-1}\|f_{12}\|_{L^2},  \\ 
&\|(-\Delta)^{\frac{k}{2}} (\lambda+ L_1^*)^{-1} F_1\|_{L^2} \leq C\|f_{11}\|_{L^2}+
C|\lambda|^{\frac{k}{2}-1}\|f_{12}\|_{L^2},
\end{aligned}
\end{align}
for $\max\{0, 1-n/4 \}< k/2 < \min\{1, n/4\}$ and 
\begin{align}
\begin{aligned}\label{resolvent-est-L1-2-DWE}
&\| (\lambda+ L_1)^{-1} F_1\|_{L^2} \leq C\|f_{11}\|_{L^2}+
C|\lambda|^{-1}\|f_{12}\|_{L^2},  \\
&\|(\lambda+ L_1^*)^{-1} F_1\|_{L^2} \leq C\|f_{11}\|_{L^2}+
C|\lambda|^{-1}\|f_{12}\|_{L^2}. 
\end{aligned}
\end{align}
\end{prop}

\vspace{2ex}

Let 
$\Sigma_{\epsilon}=\{z \in \mathbb{C}\setminus\{0\}; |\mbox{arg}z|\leq \pi-\epsilon\}$. Let $\Gamma$ be an integral curve in $\Sigma_{\epsilon}$ satisfying that 
for $t>0$, $\Gamma=\Gamma_+ \cap \Gamma_0 \cap \Gamma_-$ with   
\begin{align*}
&\Gamma_\pm=\{\lambda \in \mathbb{C}; \lambda =  -r^2 \pm\Big(r + \frac{1}{t}\Big)i, \, r> 0\}, \\
&\Gamma_0 = \{\lambda \in \mathbb{C};\lambda=\frac{1}{t}e^{-i\theta}, -\frac{\pi}{2} \leq \theta \leq \frac{\pi}{2}\}.   
\end{align*}
Then the Cauchy formula of the semigroup generated by $L_1$ yields  
\begin{align}\label{cauchy-formula}
S_1(t)=
e^{-t L_1}=
\frac{1}{2\pi i}\displaystyle\int_{\Gamma}(\lambda+L_1)^{-1}e^{\lambda t}d\lambda. 
\end{align}

The same arguments as those in Proposition \ref{analytic semigroup3} and Corollary \ref{analytic semigroup3-divergence-form} yield the followig estimates. 

\vspace{2ex}

\begin{prop}\label{analytic semigroup3-DWE}
{\rm (i)} For $u_1=\trans(\phi_1,w_1)=S_1 (t)u_{01}=e^{-t L_1}u_{01}$, we have the decay estimates  
$$
\| u_1\|_{L^2} \leq C(1+t)^{-\delta_2+\epsilon}\|u_{01}\|_{L^p},
$$ 
where $\delta_2=\dfrac{n}{2}\Big(\dfrac{1}{p}-\dfrac{1}{2}\Big)$ for $1< p< 2$ with 
$\max\{0,1-n/4\} < \delta_2 <\min\{1, n/4\}$ and any small $\epsilon>0$ independent of $c_0$, and $\delta_2$ or $p=2$ and $\delta_2=\epsilon=0$.  Furthermore, we have  that 
$$
\| \nabla u_1\|_{L^2} \leq C(1+t)^{-\delta_2-\frac{1}{2}}\|u_{01}\|_{L^p}. 
$$

{\rm (ii)} Let $u_1=\trans(\phi_1,w_1)=S_1 (t) \del_x u_{01}=e^{-t L_1}u_{01}$, where 
$u_{01} \in L^p_{(1)}$ for $1<p<2$ close to $2$ satisfying that $ \delta_3+\dfrac{1}{2}< \min\{1, \dfrac{n}{4}\}$ with $\delta_3=\dfrac{n}{2}\Big(\dfrac{1}{p}-\dfrac{1}{2}\Big)+\epsilon$, where $\epsilon>0$ is any small constant.  Then  we have that 
$$
\| u_1\|_{L^2} \leq C(1+t)^{-\delta_3+\epsilon-\frac{1}{2}}\|u_{01}\|_{L^p}. 
$$
\end{prop}

\end{subsection}

\begin{subsection}{Estimates of the high frequency part} 

On the high frequency part, we apply the $L^2$ energy estimates. We consider the following linearized system 

\begin{align}
\begin{aligned}\label{eq-DWE-high}
    \del_{tt} u_\infty - \mu \Delta \del_t u_\infty - \mu' \Delta u_\infty + B[b]u_\infty =f_\infty, \\
u_\infty(0,x)=u_{0\infty},  \ \ \del_t u_\infty(0,x)=v_{0\infty}. 
\end{aligned}
\end{align}

Taking the $L^2$ inner products \eqref{eq-DWE-high} with $u_\infty$ and $\del_t u_\infty$ respectively and applying the Cauchy-Schwartz inequality,  there exist positive constants $d_1 $ and $d_2$ such that 
\begin{align}\label{energy-high}
    \dfrac{d}{dt}\|\del_t u_\infty\|_{L^2}^2 
    + \mu\dfrac{d}{dt}\|\nabla  u_\infty\|_{L^2}^2 + \del_t(\del_tu_\infty, u_\infty) +d_1 \|\nabla u_\infty\|_{L^2}^2+ d_2 \|\nabla \del_t u_\infty\|_{L^2}^2 \leq C\|f_\infty\|_{L^2}^2. 
\end{align}
    Here we used Lemma \ref{lemPinfty} and thus  terms in the right hand side were absorbed  to the left hand sides in \eqref{energy-high}. 
    We also used that 
    $$
 (B[b]u_\infty, u_\infty) \leq 
 C\epsilon\|\nabla u_\infty\|_{L^2}^2,  \ \ 
  (B[b]u_\infty, \del_t u_\infty) \leq 
 C\epsilon\|\nabla u_\infty\|_{L^2}\|\nabla \del_t u_\infty\|_{L^2}.  
    $$
Indeed,  by Assumption \ref{ass-DWE}, integration by parts and  Lemma \ref{lemPinfty}, we get that 
\begin{align*}
| (u_\infty\div b_1, u_\infty)| 
&= |-(b_1, \nabla |u_\infty|^2)| 
\leq C \|b_1\|_{L^{n,\infty}}\|u_\infty\|_{L^{2^*, 2}}\|\nabla u_\infty\|_{L^2} \\
&\leq C\epsilon\|\nabla u_\infty\|_{L^2}^2, \\
| (b_2 \cdot \nabla u_\infty , u_\infty)| 
&\leq C \|b_2\|_{L^{n,\infty}}\|u_\infty\|_{L^{2^*, 2}}\|\nabla u_\infty\|_{L^2}\\
&\leq C\epsilon\|\nabla u_\infty\|_{L^2}^2,\\
| (b_3 \Delta u_\infty , u_\infty)| 
&\leq  |-(\nabla b_3 \cdot \nabla u_\infty, u_\infty)|
+ |-( b_3 \cdot \nabla u_\infty, \nabla u_\infty)|
\leq C \|b_3\|_{W^{1,\infty}}\|\nabla u_\infty\|_{L^2}^2\\
&\leq C\epsilon\|\nabla u_\infty\|_{L^2}^2.\\
\end{align*} 
We note that by Lemma \ref{lemPinfty} 
    $$
|(\del_tu_\infty, u_\infty)| \leq \frac{1}{2}\|\del_t u_\infty\|_{L^2}^2 + C_1 \|\nabla u_\infty\|_{L^2}^2.
    $$
Hence if 
\begin{align} \label{mu-condition}
\mu \geq 2C_1,
\end{align}
$C_1 \|\nabla u_\infty\|_{L^2}^2$ is also absorbed to $\mu \|\nabla u_\infty\|_{L^2}^2$. Therefore, let us define an energy functional ${\cal E}[u_\infty]$ by 
\begin{eqnarray*}
{\cal E}[u_\infty](t) =\|\nabla u_\infty\|_{L^2}^2+ \|\del_t u_\infty\|_{L^2}^2. 
\end{eqnarray*}
The we see from \eqref{energy-high} that 
\begin{align}\label{energy-high-DWE}
    \dfrac{d}{dt}{\cal E}[u_\infty]   +d_1 \|\nabla u_\infty\|_{L^2}^2+ d_2 \|\nabla \del_t u_\infty\|_{L^2}^2 \leq C\|f_\infty\|_{L^2}^2. 
\end{align}
    
\end{subsection}

\begin{subsection}{Deriving an a priori estimate} 
Finally, we derive an a priori estimate Theorem \ref{decay-DWE}. We consider the perturbation system \eqref{DEW-evolution-low}-\eqref{DEW-evolution-high}. For any $0<T< + \infty$ we define a time weighted function space $Z(0,T)$ by 
\begin{eqnarray*}
\lefteqn{Z(0,T)=\{ U =\trans(u,v) =P_1 U +P_\infty U,  \ \  P_1 U=: U_1  =\trans(u_1, v_1) \in C([0,T]; L^2_{(1)}), }\\
&& \qquad P_\infty U=: U_\infty = \trans(u_\infty, v_\infty)\in C([0,T]; H^1_{(\infty)}\times L^2_{(\infty)}),  \ \ \|U_1,  U_\infty\|_{Z(0,T)} <+\infty\}
\end{eqnarray*}
with the norm $\|U_1,  U_\infty\|_{Z(0,T)}$ defined by 
\begin{eqnarray*}
\|U_1, U_{\infty}\|_{Z(0,T)} &=&
\sup_{0\leq t \leq T}(1+t)^{\delta}\|U_1\|_{L^2}\\ &&+\sup_{0\leq t \leq T}(1+t)^{\delta+\frac{1}{2}}(\|\nabla U_1\|_{L^2}+ \|U_\infty\|_{H^1}+ \|\del_t U_\infty\|_{L^2}),  
\end{eqnarray*}
where $\delta=\frac{n}{2}\Big(\frac{1}{p}-\frac{1}{2}\Big)-\epsilon$ for any small $\epsilon>0$ satisfying that $\max\{0,1-\frac{n}{4}\}<\delta < \min\{1, \frac{n}{4}\}$ and $1<\delta+\frac{1}{2}$. 
Let 
$$
E_0 = \|
u_0\|_{L^p\cap H^1 }+ \|
v_0\|_{L^p \cap H^1}. 
$$
\vspace{2ex}

On the low frequency part, we have the following estimate. 

\vspace{2ex}

\begin{prop}\label{est-low-result-DWE}
Let $$
U_1(t)=S_1 (t)U_{01} -\displaystyle\int_0^t S_1 (t-\tau)P_1 \mathcal{B}[b]U_\infty d\tau. 
$$ 
Then it holds that 
\begin{align*}
   \sup_{0\leq t \leq T}\{(1+t)^{\delta}\|  U_1\|_{L^2} +(1+t)^{\delta+\frac{1}{2}}\|\nabla  U_1\|_{L^2} \}\leq C_0 E_0 +C_1 \epsilon \|U_1, U_\infty\|_{Z(0,T)}, 
\end{align*}
where $C_0$ and $C_1$ are independent of $T$. 
\end{prop}

\vspace{2ex}

Proposition \ref{est-low-result-DWE} follows from Proposition \ref{analytic semigroup3-DWE} and a similar argument to that in Proposition \ref{est-low-result}. We note that for $1/p= 1/p_0 +1/2$ 
$$
\| u_\infty \div b_1\|_{L^p}+\|b_2 \nabla u_\infty\|_{L^p} \leq C(\|\nabla b_1\|_{L^{p_0}}+\|b_2\|_{L^{p_0}})  \|\nabla u_\infty\|_{L^2}\leq 
C\epsilon \|\nabla u_\infty\|_{L^2}.  
$$
In addition, since 
$$
b_3 (x)\Delta u_\infty = \div (b_3 (x)\nabla u_\infty)- \nabla b_3(x) \cdot \nabla u_\infty, 
$$
hence we can apply Proposition \ref{analytic semigroup3-DWE} {\rm (ii)} for $b_3 (x)\Delta u_\infty$. We also note that 
$$
\|b_3 \nabla u_\infty\|_{L^p}+\|\nabla b_3\cdot  \nabla u_\infty\|_{L^p} \leq C(\|b_3\|_{L^{p_0}}+\|\nabla b_3\|_{L^{p_0}})  \|\nabla u_\infty\|_{L^2}\leq 
C\epsilon \|\nabla u_\infty\|_{L^2}.  
$$

On the other hand, for the high frequency part, let $D[u_\infty]=\|\nabla u_\infty\|_{L^2}^2+  \|\nabla \del_t u_\infty\|_{L^2}^2 $.  
We apply \eqref{energy-high-DWE} 
for $f_\infty=P_\infty \mathcal{B}[b]U_1$. We note that Lemma \ref{lemP_1} and the Sobolev iequality,  \ref{Hardy-Rellich} verify that 
$$
\|\mathcal{B}[b]U_1\|_{L^2} \leq C\epsilon \|\nabla U_1\|_{L^2}.
$$
In fact, 
\begin{align*}
\begin{aligned}
\|u_1 \div b_1\|_{L^2} 
&\leq C\|\nabla b_1\|_{L^{n/2, \infty}}\|u_1\|_{L^{q,2}} 
\leq C \|\nabla b_1\|_{L^{n/2, \infty}}
\|(-\Delta)^{\frac{n}{2}\big(\frac{1}{2}-\frac{1}{q}\big)} u_1\|_{L^2}\\
& \leq C\epsilon \|\nabla u_1\|_{L^{2}},\\
\|b_2 \cdot \nabla u_1\|_{L^2}& \leq 
C\|b_2\|_{L^{n,\infty}}\|\nabla u_1\|_{L^{2^*,2}}  \leq C\|b_2\|_{L^{n,\infty}}\|\nabla^2 u_1\|_{L^{2}} \\
&\leq C\epsilon \|\nabla u_1\|_{L^{2}},\\
\|b_3 \Delta u_1\|_{L^2}& \leq 
C\|b_3\|_{L^{\infty}}\|\Delta u_1\|_{L^{2}}\\
&\leq C\epsilon \|\nabla u_1\|_{L^{2}},
\end{aligned}
\end{align*}
where $2<q$ satisfying that $1/2=2/n +1/q$, as well as $2^*<q$. 
Hence it follows from \eqref{energy-high-DWE} that there exists a positive constant $d$ such that   
\begin{align}
\begin{aligned}\label{estenergy-DWE}
\lefteqn{{\cal E}[U_{\infty}](t)+d \displaystyle\int_{0}^{t}e^{-C_2(t-\tau)}D[U_\infty](\tau)d\tau}\\
&\quad \leq e^{-C_2 t}{\cal E}[U_{\infty}](0) \\
&\qquad  +C\epsilon^2\|U_1, U_\infty\|_{Z(0,T)}^2 \displaystyle\int_{0}^{t}
e^{-C_2(t-\tau)}(1+\tau)^{-2\delta-1}d\tau\\
&\quad \leq e^{-C_2 t}{\cal E}[U_{\infty}](0) 
+ C\epsilon^2 \|U_1, U_\infty\|_{Z(0,T)}^2(1+t)^{-2\delta-1}. 
\end{aligned}
\end{align}
Define $\tilde{\cal E}[U_\infty]$ by  
\begin{eqnarray*}
\tilde{\cal E}[U_\infty](t) &=&\sup_{0\leq \tau \leq t}(1+\tau)^{2\delta+1}{\cal E}[U_{\infty}](\tau)
\end{eqnarray*}
We see from \eqref{estenergy-DWE} that 
\begin{eqnarray}
\tilde{\cal E}[U_\infty](t)+d{\cal D}[U_\infty](t)
\leq  C_0^2E_0^2+ C_3^2\epsilon^2\|U_1, U_\infty\|_{Z(0,T)}^2 
 \label{est-high-result-DWE}
\end{eqnarray}
where the constant $C_3$ is independent of $T$. By Proposition \ref{est-low-result-DWE} and \eqref{est-high-result-DWE}, taking $\epsilon \leq 1/2 \times \max\{C_1, C_3\}$, we obtain the a priori estimate Theorem \ref{decay-DWE}. This completes the proof. 

\vspace{2ex}

{\rm \begin{rem}\label{existence-sol-DWE} 
Concerning  existence of a solution $u \in C([0,T]; H^1) \cap L^2 (0,T; H^2)$ with $\del_t u \in C([0,T]; L^2)$  to \eqref{eq-DWE}, 
Under $(u_0,v_0) \in L^p \cap H^2$, 
\cite[Proposition 2.1]{ITY} implies 
existence of the semigroup of \eqref{eq-DWE} 
in $C([0,T];H^2)$ without the perturbation terms having $b_j$. Hence there exists a solution $U=\trans(u, \del_t u)$ in $C([0,T];H^2)$ for initial data $U_{0}$ without the perturbation terms having $b_j$.  

We show existence of a solution adding the perturbation terms. 
Let $U=U_1+U_\infty=P_1 U+P_\infty U$. 
By Proposition \ref{analytic semigroup3-DWE}, we can consider the semigroup $S_1(t)=
e^{-t L_1}$ for the low frequency part with the same decay rates.  Since $U_1$ satisfies that 
$$
U_1(t)=S_1 (t)U_{01} -\displaystyle\int_0^t S_1 (t-\tau)P_1 \mathcal{B}[b]U_\infty d\tau,  
$$ 
applying similar arguments to Proposition \ref{est-low-result-DWE} 
we see that 
$$
\|U_1(t)\|_{L^2} \leq C \|U_0\|_{L^2} +C\epsilon T\sup_{0<t \leq T}\|u_1,u_\infty\|_{L^2 \times H^1}
$$
where $C$ are independent of $t$ and $T$. Lemma \ref{lemP_1} yields 
$$
\|\nabla^2 U_1(t)\|_{L^2} \leq C \|U_0\|_{L^2} +C\epsilon T\sup_{0<t \leq T}\|u_1,u_\infty\|_{L^2 \times H^1}.  
$$
On the other hand, the energy estimate yields a similar estimate to \eqref{estenergy-DWE} in the high frequency part, {\rm i.e.,} 
\begin{align}
\begin{aligned}\label{estenergy-DWE-2}
\lefteqn{{\cal E}[U_{\infty}](t)+d \displaystyle\int_{0}^{t}e^{-C_2(t-\tau)}D[U_\infty](\tau)d\tau}\\
&\quad \leq e^{-C_2 t}{\cal E}[U_{\infty}](0) \\
&\qquad  +C\epsilon^2\|u_1, u_\infty\|_{C([0,T];L^2\times H^1)}^2 \displaystyle\int_{0}^{t}
e^{-C_2(t-\tau)}d\tau\\
&\quad \leq e^{-C_2 t}{\cal E}[U_{\infty}](0) 
+ C\epsilon^2 \|u_1, u_\infty\|_{C([0,T];L^2\times H^1)}^2. 
\end{aligned}
\end{align}
Therefore, the iteration argument derives that  if $\epsilon$ is sufficiently small, 
there exists a solution $U =\trans(u, \del_t u)$ with $u \in C([0,T]; H^1)$ and $\del_t u \in C([0,T]; L^2)$ to  \eqref{eq-DWE} with  the perturbation terms having $b_j$. 
By Lemma \ref{lemP_1} and the energy inequality, we also have that $\nabla^2 u \in L^2 (0,T; L^2)$.  
\end{rem}
}

\end{subsection}

\end{section}

\vspace{2ex}

\section{Appendix}
Our aim is to prove the following claims:\\

(i) 
$ v\in \widehat {B}^{1/2}_{2,\infty}(\IR^3)$ when $(1+|x|)v\in L^\infty$ and $(1+|x|)^2 \nabla v \in L^\infty.$\\

(ii) 
$L^{3,\infty} \not\subset \widehat{B}^{1/2}_{2,\infty}$ on $\mathbb{R}^3$.

\vspace{2ex}

\noindent
{\it Proof of (i).} 
We assume that  
$v=v(x)$ satisfies
\begin{equation}\label{assump}
M_1:=\|(1+|x|)\,v\|_{L^\infty}<\infty,\quad
M_2:=\|(1+|x|)^2\,\nabla v\|_{L^\infty}<\infty.
\end{equation}
By the finite-difference characterization of homogeneous Besov spaces, see \cite[Proposition 2.36]{BCD}, one has
\BLACK for $0<s<1$ and $1\le p\le\infty$  
\begin{equation}\label{eq:diff}
\|f\|_{\widehat {B}^{s}_{p,\infty}}\cong  \sup_{0<|h|} |h|^{-s}\,\|f(\cdot+h)-f(\cdot)\|_{L^p},
\end{equation}
for tempered distributions with at most polynomial growth. In particular, for $s=\tfrac12$, $p=2$,
\begin{equation}\label{eq:diff-half}
\|f\|_{\widehat {B}^{1/2}_{2,\infty}}\cong \sup_{0<|h|} |h|^{-1/2}\,\|f(\cdot+h)-f(\cdot)\|_{L^2}.
\end{equation}

\vspace{2ex}


\begin{prop}
Under \eqref{assump}, $v\in \widehat {B}^{1/2}_{2,\infty}(\IR^3)$ and
\begin{equation}\label{eq:main-bound}
\|v\|_{\widehat {B}^{1/2}_{2,\infty}}\;\lesssim\; \|\nabla v\|_{L^2}+M_1+M_2.
\end{equation}
\end{prop}

\vspace{2ex}

\begin{proof}[Proof] 
By \eqref{eq:diff-half}, it is sufficient to show a  bound 
\begin{align}\label{eq:main-bound2}
\begin{aligned} 
\sup_{0<|h|} |h|^{-1/2}\,& \|v(\cdot+h)-v(\cdot)\|_{L^2} \\
&\leq \sup_{0<|h|\leq 1} |h|^{-1/2}\,\|v(\cdot+h)-v(\cdot)\|_{L^2}
+ \sup_{1 \leq |h|} |h|^{-1/2}\,\|v(\cdot+h)-v(\cdot)\|_{L^2} \\
& =: I_1 +I_2.\BLACK 
\end{aligned}
\end{align}
Concerning $I_1$, fix $0<|h|\le 1$ and set $R:=|h|^{-1}\ge 1$. Take a smooth radial cutoff function $\chi\in C_0^\infty(\IR^3)$ with $\chi\equiv 1$ on $B_1(0)$, $\supp\chi\subset B_2(0)$, and define $\chi_R(x):=\chi(x/R)$. Split
\[
v=v_{\mathrm{loc}}+v_{\mathrm{far}},\qquad v_{\mathrm{loc}}:=v\,\chi_R,\quad v_{\mathrm{far}}:=v\,(1-\chi_R).
\]

\smallskip
\noindent\emph{1) Estimate of $v_{\rm loc}$.
} 
By the fundamental theorem of calculus,
\[
v_{\mathrm{loc}}(x+h)-v_{\mathrm{loc}}(x)=\int_0^1 h\cdot\nabla v_{\mathrm{loc}}(x+\theta h)\,d\theta.
\]
Taking $L^2$ norms and using $\nabla(v\chi_R)=\chi_R\nabla v + v\,\nabla\chi_R$,
\begin{equation}\label{eq:loc-1}
\|v_{\mathrm{loc}}(\cdot+h)-v_{\mathrm{loc}}(\cdot)\|_{L^2(B_{3R}(0))}
\le |h|\Big(\|\nabla v\|_{L^2(B_{3R}(0))}+\|v\,\nabla\chi_R\|_{L^2}\Big).
\end{equation}
Clearly $\|\nabla v\|_{L^2(B_{3R}(0))}\le \|\nabla v\|_{L^2}$. Moreover, $\nabla\chi_R(x)=R^{-1}(\nabla\chi)(x/R)$ is supported in the annulus $A_R:=\{R\le |x|\le 2R\}$ and $\|\nabla\chi_R\|_{L^\infty}\lesssim R^{-1}$. Using $|v(x)|\le M_1(1+|x|)^{-1}$ (from \eqref{assump}) we obtain
\begin{align}
\|v\,\nabla\chi_R\|_{L^2}
&\lesssim R^{-1}\Big(\int_{A_R} |v(x)|^2\,dx\Big)^{1/2}
\le R^{-1} M_1\Big(\int_{A_R} (1+|x|)^{-2}\,dx\Big)^{1/2}\nonumber\\
&\lesssim R^{-1} M_1\Big(\int_{R}^{2R} r^{-2}\cdot r^2\,dr\Big)^{1/2}
\lesssim R^{-1/2} M_1.\label{eq:loc-2}
\end{align}
Combining \eqref{eq:loc-1} and \eqref{eq:loc-2},
\begin{equation}\label{eq:loc-final}
\|v_{\mathrm{loc}}(\cdot+h)-v_{\mathrm{loc}}(\cdot)\|_{L^2}
\;\lesssim\; |h|\,\|\nabla v\|_{L^2} + |h|^{3/2}\,M_1.
\end{equation}

\smallskip
\noindent\emph{2) Estimate of $v_{\rm far}$.}
For $x$ with $|x|\ge R$ and $\theta\in[0,1]$, we have $|x+\theta h|\ge |x|-|h|\ge |x|-1$, 
hence $(1+|x+\theta h|)\gtrsim |x|$. \BLACK  By the Taylor inequality  and \eqref{assump},
\begin{align*}
|v_{\mathrm{far}}(x+h)-v_{\mathrm{far}}(x)|
& \lesssim  |h|\,\|(1+|x|)^2\nabla v\|_{L^\infty}\,|x|^{-2} 
\leq \BLACK |h|\,M_2\,|x|^{-2}.
\end{align*}
Therefore,
\begin{equation}\label{eq:far-final}
\|v_{\mathrm{far}}(\cdot+h)-v_{\mathrm{far}}(\cdot)\|_{L^2}
\lesssim |h|\,M_2\Big(\int_{|x|\ge R} |x|^{-4}\,dx\Big)^{1/2}
\lesssim |h|^{3/2} M_2,
\end{equation}
since $\int_{|x|\ge R}|x|^{-4}\,dx \simeq \int_{R}^\infty r^{-2}\,dr \simeq R^{-1}=|h|$.

\smallskip
\noindent\emph{3)}
Adding \eqref{eq:loc-final} and \eqref{eq:far-final} and dividing by $|h|^{1/2}$ 
\BLACK yields
\[
|h|^{-1/2}\,\|v(\cdot+h)-v(\cdot)\|_{L^2}
\;\lesssim\; |h|^{1/2}\,\|\nabla v\|_{L^2} + |h|\,(M_1+M_2).
\]
Taking the supremum over $0<|h|\le 1$ proves the bound of $I_1$.  

Concerning $I_2$, set $R:=|h|\ge 1$ and define  $\chi_R(x):=\chi(x/2R)$.  We also define $v_{\mathrm{loc}}$ and $v_{\mathrm{far}}$ similarly to $I_1$. Then a direct computation yields
\begin{align}
\begin{aligned}\label{eq:loc-final2}
\|v_{\mathrm{loc}}(\cdot+h)-v_{\mathrm{loc}}(\cdot)\|_{L^2}
& \lesssim  \|v_{\mathrm{loc}}\|_{L^2} \\
& \lesssim M_1 \Big(\displaystyle\int_{|x|\leq 4R} 
\dfrac{1}{(1+|x|)^2} \, dx \Big)^\frac{1}{2}\\
& \lesssim M_1 R^{\frac{1}{2}}
\simeq M_1 |h|^{\frac{1}{2}}.
\end{aligned}
\end{align}
As for $v_{\mathrm{far}}$, since 
$|x+\theta h|\ge |x|-|h|\ge \frac{|x|}{2}$ for $|x| \geq 2R$, 
{\em i.e.,} $1+|x+\theta h| \ge \frac{|x|}{2}$, a similar argument to \eqref{eq:far-final} shows that 
\begin{equation}\label{eq:far-final2}
\|v_{\mathrm{far}}(\cdot+h)-v_{\mathrm{far}}(\cdot)\|_{L^2}
\lesssim |h|\,M_2\Big(\int_{|x|\ge 2R} |x|^{-4}\,dx\Big)^{1/2}
\lesssim |h|^{1/2} M_2.
\end{equation}
Hence \eqref{eq:loc-final2} and \eqref{eq:far-final2} derive that 
\[
|h|^{-1/2}\,\|v(\cdot+h)-v(\cdot)\|_{L^2}
\;\lesssim\; M_1+M_2.
\]
Taking the supremum over $1\leq |h|$  we  prove the bound of $I_2$.

By \eqref{eq:diff-half} this is equivalent to $v\in \widehat {B}^{1/2}_{2,\infty}(\IR^3)$.  
\end{proof}

\vspace{1ex}



\begin{proof}[Proof of (ii)]
Let $\psi \in C_c^\infty(\mathbb{R}^3)$, $\psi \not\equiv 0$, satisfy $\int_{\mathbb{R}^3} \psi(x)\,dx = 0$. Fix $L \gg 1$ and define points $x_k := kLe_1$ so that the supports of $\psi(\cdot - x_k)$ are pairwise disjoint. Set amplitudes $a_k := k^{-1/3}$ and for $N\in\mathbb{N}$ define
\begin{equation*}
    f_N(x) := \sum_{k=1}^N a_k\, \psi(x - x_k).
\end{equation*}
Let $S_k := \mathrm{supp}\,\psi(\cdot - x_k)$. Since the supports are pairwise disjoint, on each $S_k$ we have $f_N (x)=a_k \psi(x-x_k)$ and thus $|f_N| = a_k |\psi(\cdot - x_k)|$.  
For any $\lambda > 0$,
\begin{equation*}
    \{ |f_N| > \lambda \} = \bigcup_{k: a_k M > \lambda} \{ |\psi(\cdot - x_k)| > \lambda/a_k \},
\end{equation*}
where $M := \|\psi\|_{L^\infty}$. Because $\psi$ is fixed and compactly supported,
\[
    |\{ |\psi(\cdot - x_k)| > \lambda / a_k \}| \le |\mathrm{supp}(\psi)| =: C_\psi,
\]
independent of $k$. Since the supports $S_k$ are  pairwise \BLACK disjoint,
\[
    |\{|f_N| > \lambda\}| \le C_\psi\, \#\{k : a_k M > \lambda\}.
\]
Hence we obtain that 
\begin{equation*}
    |\{ |f_N| > \lambda \}| \lesssim \#\{k: a_k \gtrsim \lambda\} \lesssim \lambda^{-3},
\end{equation*}
which implies
\begin{equation*}
    \|f_N\|_{L^{3,\infty}} \le C, \quad \text{independent of $N$}.
\end{equation*}

Choose a dyadic index $j_0$ with $\|\dot{\Delta}_{j_0} \psi\|_{L^2} > 0$ in the homogeneous Besov norm. 
Here $\dot{\Delta}_{j_0}\psi = 
\mathcal{F}^{-1}(\hat{\phi}_j \hat{\psi})$, where $\hat{\phi}_j(\xi) = \hat{\phi}(2^{-j}\xi)$ and $\phi$ defines a homogeneous Littlewood-Paley dyadic decomposition; see \cite[Section 2]{BCD} for details. 
Then
\begin{equation*}
    \dot{\Delta}_{j_0} f_N = \sum_{k=1}^N a_k\, \dot{\Delta}_{j_0}\psi(\cdot - x_k).
\end{equation*}
Because the supports are far apart and $\dot{\Delta}_{j_0}\psi$ is in the Schwartz class and $x_k= kLe_1$, the  cross terms in the inner products are summable  and small by taking $L \geq 1$  for $k\neq k'$. To be more precise, \BLACK 
\begin{align}\label{cross-term}
|(\dot{\Delta}_{j_0} \psi(\cdot-x_k ), \dot{\Delta}_{j_0} \psi(\cdot-x_{k'}))| \leq  C\frac{1}{L^2|k-k'|^2} \|\dot{\Delta}_{j_0} \psi\|_{L^2}^{2}. 
\end{align}
Indeed, let $g=\dot{\Delta}_{j_0} \psi$. 
By the definition, 
$$
(g(\cdot-x_k), g(\cdot-x_{k'}))=
g \ast \tilde{g}(x_{k'}-x_k) = g \ast \tilde{g}((k'-k)Le_1),
$$
where $\tilde{g}(x)=\overline{g(-x)}$. Moreover, 
we note by the Plancherel theorem that 
$$
(g(\cdot-x_k), g(\cdot-x_{k'})) 
=(2\pi)^{-n}(\mathcal{F}g(\cdot-x_k ), \mathcal{F}g(\cdot-x_{k'})).
$$
Since $\mathcal{F}g(\cdot-x_k)(\xi) =  \hat g(\xi) \BLACK e^{-i x_k \cdot\xi}$, we obtain  that 
$$
(g(\cdot-x_k), g(\cdot-x_{k'})) 
= (2\pi)^{-n} \displaystyle\int_{\mathbb{R}^n}
|\hat{g}|^2 e^{-i(x_{k'}-x_k)\cdot \xi} d\xi= (2\pi)^{-n} \displaystyle\int_{\mathbb{R}^n}
|\hat{g}|^2 e^{-i Le_1(k'-k) \cdot \xi} d\xi. 
$$
We set 
$h(z)= g \ast \tilde{g}(z)= (2\pi)^{-n} 
\int_{\mathbb{R}^n}
|\hat{g}|^2 e^{-i z \cdot \xi} d\xi. $
Since $\hat{g}\in C^{\infty}_c$, an integration by parts with respect to $\xi_1$ twice, we have that
$$
|h(z)| \leq C\frac{1}{|z_1|^2},
$$
where $z_1= z \cdot e_1$. Hence with 
$z=Le_1(k'-k)$ and $L \geq 1$,  
$$
|(g(\cdot-x_k ), g(\cdot-x_{k'}))|
\leq C\frac{1}{L^2|k-k'|^2} 
\|g\|_{L^2}^2. 
$$
\BLACK

We note that 
$$
\Big\|\sum_{k=1}^N a_k  \dot{\Delta}_{j_0}\psi(\cdot - x_k)\Big\|_{L^2}^2 = \Big(\sum_{k=1}^N a_k^2\Big)\, \|\dot{\Delta}_{j_0}\psi\|_{L^2}^2 + (\mbox{small cross terms}).
$$ 
When $k\not =k'$, {\rm i.e.} $|k-k'|\geq 1$,  we see from \eqref{cross-term} that the cross terms are bounded by 
$$
\sum_{k, k'=1,\, k\neq k'}^N |a_k a_{k'}| \frac{1}{L^2 |k-k'|^2}  \|\dot{\Delta}_{j_0} \psi\|_{L^2}^2 \leq C L^{-2} \Big(\sum_{k=1}^N |a_k|^2\Big) \|\dot{\Delta}_{j_0} \psi\|_{L^2}^2. 
$$ 
Here 
we used the H\"older and Young (Minkowski) inequalities on sequences \cite[Theorem (20.18)]{Hewitt-Ross}, 
{\rm i.e.,} for $\mathfrak{a}=(b_k)_{k\in \mathbb{Z}}$, where $b_k=a_k$ on $k \in \mathbb{N}$ and otherwise $b_k=0$,  and 
$\mathfrak{w}=(w_m)_{m \in \mathbb{Z}}$, where $w_m= |m|^{-2}$ for $m\not=0$ and $w_m=0$ for $m=0$,  
\begin{align*}
|(\mathfrak{a}, \mathfrak{w}\ast \mathfrak{a})_{\ell^2} | & \leq 
C\|\mathfrak{a}\|_{\ell^2}
\|\mathfrak{w}\ast \mathfrak{a}\|_{\ell^2}
 \leq 
C\|\mathfrak{a}\|_{\ell^2} \|\mathfrak{w}\|_{\ell^1} \|\mathfrak{a}\|_{\ell^2}\\
& = C\zeta(2)\, \|\mathfrak{a}\|_{\ell^2}^2 =
C\frac{\pi^2}{6}\|\mathfrak{a}\|_{\ell^2}^2.
\end{align*}
\BLACK Therefore, 
$$
\sum_{k, k'=1,\, k\neq k'}^N |a_k a_{k'}| \frac{1}{L^{2}|k-k'|^2}  \|\dot{\Delta}_{j_0} \psi\|_{L^2}^2 
\leq C L^{-2}\|\mathfrak{a}\|_{\ell^2}^2 
\leq C L^{-2} \Big(\sum_{k=1}^N |a_k|^2\Big) \|\dot{\Delta}_{j_0} \psi\|_{L^2}^2. 
$$ 


Hence we have that for a large $L$, 
\begin{equation*}
    \|\dot{\Delta}_{j_0} f_N\|_{L^2}^2 \gtrsim \sum_{k=1}^N a_k^2\, \|\dot{\Delta}_{j_0}\psi\|_{L^2}^2 = c \sum_{k=1}^N k^{-2/3}.
\end{equation*}
Therefore, 
\begin{equation*}
    \|f_N\|_{\widehat{B}^{1/2}_{2,\infty}} \ge 2^{j_0/2}\, \|\dot{\Delta}_{j_0} f_N\|_{L^2} \gtrsim N^{1/6} \to \infty \quad (N\to\infty).
\end{equation*}
Thus we have constructed a sequence $(f_N)$ bounded in $L^{3,\infty}$ but unbounded in $\widehat {B}^{1/2}_{2,\infty}$, proving that $L^{3,\infty} \not\subset \widehat {B}^{1/2}_{2,\infty}$. 
\end{proof} 

\BLACK 

\BLACK


\vspace{2ex}

\noindent {\bf Acknowledgements.}
I would like to express my appreciation to Professor
 Reinhard Farwig for his valuable advices and suggestions.  I also would like to express my appreciation to Professor
 Yoshiyuki Kagei and Takayuki Kobayashi for their kindly advices and comments when I was a post doctoral fellow. 
 
 The author is supported by JSPS grant whose number is 22K13946.\vspace{1ex} 

\noindent {\bf Statements and Declarations}\vspace{1ex}
 
\noindent {\bf Conflicts of interest statement.}  There is no conflict of interest. \vspace{1ex} 

\noindent {\bf Data Availability statement.} No datasets were generated or analysed during the current study.

%

\end{document}